\documentclass{gtart}

\def\ifplaintex{\expandafter\ifx\csname documentclass\endcsname\relax}

\ifplaintex 
\else
\fi

\expandafter\ifx\csname beginpicture\endcsname\relax
\expandafter\ifx\csname documentclass\endcsname\relax
\input pictex \else
\input prepictex \input pictex \input postpictex \fi\fi

\def\gt{{\mathsurround=0pt\it $\cal G\mskip-2mu$eometry \&\ 
$\cal T\!\!$opology}}        

\def\gtp{{\mathsurround=0pt\it $\cal G\mskip-2mu$eometry \&\ 
$\cal T\!\!$opology $\cal P\!$ublications}}  

\def\lognumber#1{\def\thelognumber{#1}}
\def\volumenumber#1{\def\thevolumenumber{#1}}
\def\papernumber#1{\def\thepapernumber{#1}}
\def\volumeyear#1{\def\thevolumeyear{#1}}

\def\pagenumbers#1#2{\def\startpage{#1}\def\finishpage{#2}}
\def\published#1{\def\publishdate{#1}}
\def\proposed#1{\def\theproposer{#1}}
\def\seconded#1{\def\theseconders{#1}}
\def\received#1{\def\receiveddate{#1}}
\def\revised#1{\def\reviseddate{#1}}
\def\accepted#1{\def\accepteddate{#1}}

\def\asciiaddress#1{\def\theasciiaddress{#1}}
\def\asciiemail#1{\def\theasciiemail{#1}}
\long\def\asciiabstract#1{\long\def\theasciiabstract{#1}}
\def\asciikeywords#1{\def\theasciikeywords{#1}}

\def\shorttitle#1{\def\theshorttitle{#1}}

\let\\\par\let\thelognumber\relax
\let\thevolumenumber\relax\let\thepapernumber\relax
\let\thevolumeyear\relax\let\thesamplenumber\relax\let\startpage\relax
\let\finishpage\relax\let\publishdate\relax\let\receiveddate\relax
\let\reviseddate\relax\let\accepteddate\relax\let\theasciititle\relax
\let\theasciiauthors\relax\let\theasciiaddress\relax
\let\theasciiabstract\relax\let\theasciikeywords\relax
\let\theasciiemail\relax\let\theshortauthors\relax\let\theshorttitle\relax

\long\def\maketitlep{   

\count0=\startpage

\gt\hfill      
\beginpicture
\setcoordinatesystem units <0.33truein, 0.33truein> point at 2.2 0.9
\setplotsymbol ({$\cal G$})
\plotsymbolspacing=9truept
\circulararc 315 degrees from 0 1 center at 0 0
\setplotsymbol ({$\cal T$})
\circulararc 315 degrees from 1 -1 center at 1 0
\endpicture
\break
{\small\ifx\thesamplenumber\relax 
Volume \else Sample
\fi\thevolumenumber\ (\thevolumeyear)
\startpage--\finishpage\nl
Published: \publishdate}
\vglue 0.5truein plus 0.4fil minus 0.1truein

{\parskip=0pt\leftskip 0pt plus 1fil\def\\{\par\smallskip}{\ifplaintex\large
\else\Large\fi\bf\thetitle}\par\medskip}   

\vglue 0pt plus 0.1fil 

{\parskip=0pt\leftskip 0pt plus 1fil\def\\{\par}{\sc\theauthors}
\par\medskip}

\vglue 0pt plus 0.1fil 

{\small\parskip=0pt\let\newline\\
{\leftskip 0pt plus 1fil\def\\{\par}{\sl\theaddress}\par}
\expandafter\ifx\theemail\relax    
\relax\else\vglue 5pt plus 0.02fil minus 2pt\def\\{\stdspace{\rm 
and}\stdspace} 
\cl{Email:\stdspace\tt\theemail}\fi
\ifx\theurl\relax                  
\relax\else\vglue 5pt plus 0.02fil minus 2pt\def\\{\stdspace{\rm 
and}\stdspace}
\cl{URL:\stdspace\tt\theurl}\fi\par}

\vglue 7pt plus 0.3fil minus 3pt

{\bf Abstract}
\vglue 5pt plus 0.1fil minus 2pt

\theabstract

\vglue 7pt plus 0.3fil minus 3pt

{\bf AMS Classification numbers}\quad Primary:\quad \theprimaryclass

Secondary:\quad \thesecondaryclass

\vglue 5pt plus 0.3fil minus 2pt

{\bf Keywords}\quad \thekeywords

\vglue 10pt plus 0.5fil minus 5pt

{\small  Proposed: \theproposer\hfill Received: \receiveddate\nl
Seconded: \theseconders\hfill 
\ifx\reviseddate\relax                         
Accepted: \accepteddate                        
\else
Revised: \reviseddate                          
\fi}
\eject
}       

\let\maketitlepage\maketitlep
\let\maketitle\maketitlepage

\font\phead=cmsl9 scaled 950
\font=cmsl9 scaled 1050
\font\pnum=cmbx10 scaled 913
\font\lnum=cmbx10 
\font\pfoot=cmsl9 scaled 950
\font=cmsl9 scaled 1050
\ifplaintex
\headline{\vbox to 0pt{\vskip -4.5mm\line{\small\phead\ifnum
\count0=\startpage ISSN 1364-0380 (on line)
1465-3060 (printed) \hfill {\pnum\folio}\else\ifodd\count0\def\\{ }%
\ifx\theshorttitle\relax\thetitle\else\theshorttitle\fi\hfill{\pnum\folio}
\else\def\\{ and }{\pnum\folio}\hfill\ifx\theshortauthors\relax\theauthors
\else\theshortauthors\fi\fi\fi}\vss}}
\footline{\vbox to 0pt{\vglue 0mm\line{\small\pfoot\ifnum\count0=\startpage
\copyright\ \gtp\hfill\else
\gt, Volume \thevolumenumber\ (\thevolumeyear)\hfill\fi}\vss
}}
\else
\makeatletter
\def\@oddhead{{\small\lhead\ifnum\count0=\startpage ISSN 1364-0380 (on line)
1465-3060 (printed) \hfill {\lnum\number\count0}\else\ifodd\count0
\def\\{ }\ifx\theshorttitle\relax \thetitle \else\theshorttitle\fi\hfill
{\lnum\number\count0}\else\def\\{ and }{\lnum\number\count0}
\hfill\ifx\theshortauthors\relax 
\theauthors\else\theshortauthors\fi\fi\fi}}\def\@evenhead{\@oddhead}
\def\@oddfoot{\small\lfoot\ifnum\count0=\startpage\copyright\ \gtp\hfill\else
\gt, Volume \thevolumenumber\ (\thevolumeyear)\hfill\fi}
\def\@evenfoot{\@oddfoot}
\makeatother
\fi

\newwrite\gtoutfile
\long\gdef\makeheadfile{  
{\def\\{, }\def\s{ }
\immediate\openout\gtoutfile head.xxx
\immediate\write\gtoutfile{To: math@arxiv.org}
\immediate\write\gtoutfile{Subject: put or rep NNNNN:pppp}
\immediate\write\gtoutfile{--text follows this line--}
\immediate\write\gtoutfile{Proxy-for: \ifx\theasciiauthors\relax
\theauthors\else\theasciiauthors\fi\s<\ifx\theasciiemail\relax\theemail\else\theasciiemail\fi>}
\immediate\write\gtoutfile{\noexpand\\}
\immediate\write\gtoutfile{Authors: \ifx\theasciiauthors\relax
\theauthors\else\theasciiauthors\fi}
\immediate\write\gtoutfile{Title: \ifx\theasciititle\relax
\thetitle\else\theasciititle\fi}
\immediate\write\gtoutfile{Subj-class: GT or SG or MG etc}
\immediate\write\gtoutfile{MSC-class: \theprimaryclass\ifx\thesecondaryclass\relax\else, \thesecondaryclass\fi}
\immediate\write\gtoutfile{Journal-ref: Geom. Topol. \thevolumenumber
(\thevolumeyear) \startpage-\finishpage}
\immediate\write\gtoutfile{Comments: Published by Geometry and Topology at}
\immediate\write\gtoutfile{\s\s http://www.maths.warwick.ac.uk/gt/GTVol\thevolumenumber/paper\thepapernumber.abs.html}
\immediate\write\gtoutfile{\noexpand\\}
\immediate\write\gtoutfile{}
\ifx\theasciiabstract\relax
\immediate\write\gtoutfile{\theabstract}\else
\immediate\write\gtoutfile{\theasciiabstract}\fi
\immediate\write\gtoutfile{}
\immediate\write\gtoutfile{\noexpand\\}
\immediate\write\gtoutfile{}
\immediate\closeout\gtoutfile}}  

\def\maketitlepage{\maketitlep\makeheadfile}
\let\maketitle\maketitlepage

\def\ifplaintex{\expandafter\ifx\csname documentclass\endcsname\relax}

\ifplaintex 
\else
\fi

\expandafter\ifx\csname beginpicture\endcsname\relax
\expandafter\ifx\csname documentclass\endcsname\relax
\input pictex \else
\input prepictex \input pictex \input postpictex \fi\fi

\def\gt{{\mathsurround=0pt\it $\cal G\mskip-2mu$eometry \&\ 
$\cal T\!\!$opology}}        

\def\gtp{{\mathsurround=0pt\it $\cal G\mskip-2mu$eometry \&\ 
$\cal T\!\!$opology $\cal P\!$ublications}}  

\def\lognumber#1{\def\thelognumber{#1}}
\def\volumenumber#1{\def\thevolumenumber{#1}}
\def\papernumber#1{\def\thepapernumber{#1}}
\def\volumeyear#1{\def\thevolumeyear{#1}}

\def\pagenumbers#1#2{\def\startpage{#1}\def\finishpage{#2}}
\def\published#1{\def\publishdate{#1}}
\def\proposed#1{\def\theproposer{#1}}
\def\seconded#1{\def\theseconders{#1}}
\def\received#1{\def\receiveddate{#1}}
\def\revised#1{\def\reviseddate{#1}}
\def\accepted#1{\def\accepteddate{#1}}

\def\asciiaddress#1{\def\theasciiaddress{#1}}
\def\asciiemail#1{\def\theasciiemail{#1}}
\long\def\asciiabstract#1{\long\def\theasciiabstract{#1}}
\def\asciikeywords#1{\def\theasciikeywords{#1}}

\def\shorttitle#1{\def\theshorttitle{#1}}

\let\\\par\let\thelognumber\relax
\let\thevolumenumber\relax\let\thepapernumber\relax
\let\thevolumeyear\relax\let\thesamplenumber\relax\let\startpage\relax
\let\finishpage\relax\let\publishdate\relax\let\receiveddate\relax
\let\reviseddate\relax\let\accepteddate\relax\let\theasciititle\relax
\let\theasciiauthors\relax\let\theasciiaddress\relax
\let\theasciiabstract\relax\let\theasciikeywords\relax
\let\theasciiemail\relax\let\theshortauthors\relax\let\theshorttitle\relax

\long\def\maketitlep{   

\count0=\startpage

\gt\hfill      
\beginpicture
\setcoordinatesystem units <0.33truein, 0.33truein> point at 2.2 0.9
\setplotsymbol ({$\cal G$})
\plotsymbolspacing=9truept
\circulararc 315 degrees from 0 1 center at 0 0
\setplotsymbol ({$\cal T$})
\circulararc 315 degrees from 1 -1 center at 1 0
\endpicture
\break
{\small\ifx\thesamplenumber\relax 
Volume \else Sample
\fi\thevolumenumber\ (\thevolumeyear)
\startpage--\finishpage\nl
Published: \publishdate}
\vglue 0.5truein plus 0.4fil minus 0.1truein

{\parskip=0pt\leftskip 0pt plus 1fil\def\\{\par\smallskip}{\ifplaintex\large
\else\Large\fi\bf\thetitle}\par\medskip}   

\vglue 0pt plus 0.1fil 

{\parskip=0pt\leftskip 0pt plus 1fil\def\\{\par}{\sc\theauthors}
\par\medskip}

\vglue 0pt plus 0.1fil 

{\small\parskip=0pt\let\newline\\
{\leftskip 0pt plus 1fil\def\\{\par}{\sl\theaddress}\par}
\expandafter\ifx\theemail\relax    
\relax\else\vglue 5pt plus 0.02fil minus 2pt\def\\{\stdspace{\rm 
and}\stdspace} 
\cl{Email:\stdspace\tt\theemail}\fi
\ifx\theurl\relax                  
\relax\else\vglue 5pt plus 0.02fil minus 2pt\def\\{\stdspace{\rm 
and}\stdspace}
\cl{URL:\stdspace\tt\theurl}\fi\par}

\vglue 7pt plus 0.3fil minus 3pt

{\bf Abstract}
\vglue 5pt plus 0.1fil minus 2pt

\theabstract

\vglue 7pt plus 0.3fil minus 3pt

{\bf AMS Classification numbers}\quad Primary:\quad \theprimaryclass

Secondary:\quad \thesecondaryclass

\vglue 5pt plus 0.3fil minus 2pt

{\bf Keywords}\quad \thekeywords

\vglue 10pt plus 0.5fil minus 5pt

{\small  Proposed: \theproposer\hfill Received: \receiveddate\nl
Seconded: \theseconders\hfill 
\ifx\reviseddate\relax                         
Accepted: \accepteddate                        
\else
Revised: \reviseddate                          
\fi}
\eject
}       

\let\maketitlepage\maketitlep
\let\maketitle\maketitlepage

\font\phead=cmsl9 scaled 950
\font=cmsl9 scaled 1050
\font\pnum=cmbx10 scaled 913
\font\lnum=cmbx10 
\font\pfoot=cmsl9 scaled 950
\font=cmsl9 scaled 1050
\ifplaintex
\headline{\vbox to 0pt{\vskip -4.5mm\line{\small\phead\ifnum
\count0=\startpage ISSN 1364-0380 (on line)
1465-3060 (printed) \hfill {\pnum\folio}\else\ifodd\count0\def\\{ }%
\ifx\theshorttitle\relax\thetitle\else\theshorttitle\fi\hfill{\pnum\folio}
\else\def\\{ and }{\pnum\folio}\hfill\ifx\theshortauthors\relax\theauthors
\else\theshortauthors\fi\fi\fi}\vss}}
\footline{\vbox to 0pt{\vglue 0mm\line{\small\pfoot\ifnum\count0=\startpage
\copyright\ \gtp\hfill\else
\gt, Volume \thevolumenumber\ (\thevolumeyear)\hfill\fi}\vss
}}
\else
\makeatletter
\def\@oddhead{{\small\lhead\ifnum\count0=\startpage ISSN 1364-0380 (on line)
1465-3060 (printed) \hfill {\lnum\number\count0}\else\ifodd\count0
\def\\{ }\ifx\theshorttitle\relax \thetitle \else\theshorttitle\fi\hfill
{\lnum\number\count0}\else\def\\{ and }{\lnum\number\count0}
\hfill\ifx\theshortauthors\relax 
\theauthors\else\theshortauthors\fi\fi\fi}}\def\@evenhead{\@oddhead}
\def\@oddfoot{\small\lfoot\ifnum\count0=\startpage\copyright\ \gtp\hfill\else
\gt, Volume \thevolumenumber\ (\thevolumeyear)\hfill\fi}
\def\@evenfoot{\@oddfoot}
\makeatother
\fi

\newwrite\gtoutfile
\long\gdef\makeheadfile{  
{\def\\{, }\def\s{ }
\immediate\openout\gtoutfile head.xxx
\immediate\write\gtoutfile{To: math@arxiv.org}
\immediate\write\gtoutfile{Subject: put or rep NNNNN:pppp}
\immediate\write\gtoutfile{--text follows this line--}
\immediate\write\gtoutfile{Proxy-for: \ifx\theasciiauthors\relax
\theauthors\else\theasciiauthors\fi\s<\ifx\theasciiemail\relax\theemail\else\theasciiemail\fi>}
\immediate\write\gtoutfile{\noexpand\\}
\immediate\write\gtoutfile{Authors: \ifx\theasciiauthors\relax
\theauthors\else\theasciiauthors\fi}
\immediate\write\gtoutfile{Title: \ifx\theasciititle\relax
\thetitle\else\theasciititle\fi}
\immediate\write\gtoutfile{Subj-class: GT or SG or MG etc}
\immediate\write\gtoutfile{MSC-class: \theprimaryclass\ifx\thesecondaryclass\relax\else, \thesecondaryclass\fi}
\immediate\write\gtoutfile{Journal-ref: Geom. Topol. \thevolumenumber
(\thevolumeyear) \startpage-\finishpage}
\immediate\write\gtoutfile{Comments: Published by Geometry and Topology at}
\immediate\write\gtoutfile{\s\s http://www.maths.warwick.ac.uk/gt/GTVol\thevolumenumber/paper\thepapernumber.abs.html}
\immediate\write\gtoutfile{\noexpand\\}
\immediate\write\gtoutfile{}
\ifx\theasciiabstract\relax
\immediate\write\gtoutfile{\theabstract}\else
\immediate\write\gtoutfile{\theasciiabstract}\fi
\immediate\write\gtoutfile{}
\immediate\write\gtoutfile{\noexpand\\}
\immediate\write\gtoutfile{}
\immediate\closeout\gtoutfile}}  

\def\maketitlepage{\maketitlep\makeheadfile}
\let\maketitle\maketitlepage

\lognumber{134}
\volumenumber{5}\papernumber{2}\volumeyear{2001}
\pagenumbers{7}{74}
\accepted{31 January 2001}
\proposed{Walter Neumann}
\seconded{Steve Ferry, Ralph Cohen}
\received{1 September 2000}
\revised{13 December 2000}
\published{02 February 2001}

\usepackage[leqno]{amsmath} 
\usepackage{amssymb} 
\usepackage{amstext}
\usepackage{amscd}

\newcommand{\myequation}[1]{$\displaystyle #1$\vskip 1.25ex plus 1ex 
minus .2ex}

\renewcommand\theenumi{\arabic{enumi}}

\makeatletter
            
  \renewcommand{\subsection}[1]{\@startsection{subsection}{2}{\z@}%
    {3.25ex\@plus1ex \@minus.2ex}%
    {-1.5ex plus .2ex}%
    {\bf}{\hspace{-.25em}#1\hspace{.5em}\ignorespaces}}%
  
  \renewcommand{\subsubsection}[1]{\def\@tempd{#1}\ifx\@tempd\@empty
	\@startsection{subsubsection}{3}{\z@}%
    {0.75ex \@plus1ex \@minus.2ex}%
    {-1ex plus .2ex}%
    {\bf}{\hspace{-.5em}\ignorespaces}
	\else
	\@startsection{subsubsection}{3}{\z@}%
    {0.75ex \@plus1ex \@minus.2ex}%
    {-1ex plus .2ex}%
    {\bf}{\hspace{-.25em}{\rm(\sl#1\rm)}\hspace{.5em}\ignorespaces}\fi}

\makeatother

\DeclareMathOperator{\I}{\bf I} 
\DeclareMathOperator{\II}{\bf II}
\DeclareMathOperator{\III}{\bf III} 
\DeclareMathOperator{\IV}{\bf IV}
\DeclareMathOperator{\V}{\bf V} 
\DeclareMathOperator{\J}{\bf J}
 
\DeclareMathOperator{\im}{Im}
\DeclareMathOperator{\Ker}{Ker} 
\DeclareMathOperator{\Link}{Link}
\DeclareMathOperator{\flag}{Flag}

\newcommand{\norm}[1]{\lVert#1\rVert} 
\def\b{\operatorname{\beta}}
\newcommand\BH{\mathbb{H}} 
\newcommand\BQ{\mathbb{Q}}
\newcommand\BE{\mathbb{E}} 
\newcommand\BZ{\mathbb{Z}}
\newcommand\BR{\mathbb{R}} 
\newcommand\BN{\mathbb{N}}
\newcommand\BS{\mathbb{S}} 
\newcommand\BX{\mathbb{X}}

\newtheorem{Theorem}[subsubsection]{Theorem}

\newtheorem{Conjecture}[subsubsection]{Conjecture}
\newtheorem{Lemma}[subsubsection]{Lemma}
\newtheorem{Corollary}[subsubsection]{Corollary}
\newtheorem{Proposition}[subsubsection]{Proposition}

\newtheorem{The Euler Characteristic Conjecture}[subsection]{The Euler
  Characteristic Conjecture} 
\newtheorem{The Flag Complex Conjecture}[subsection]{The Flag Complex
  Conjecture} 
\newtheorem{The Vanishing Conjecture}[subsection]{Singer's Conjecture}
\newtheorem{uConjecture}[subsection]{Conjecture}
\newtheorem*{uTheorem}{Theorem} 
\newtheorem*{I(n)}{$\I(n)$}
\newtheorem*{II(n)}{$\II(n)$} 
\newtheorem*{III(2k+1)}{$\III(2k+1)$}
\newtheorem*{III'(2k+1)}{$\III'(2k+1)$} 
\newtheorem*{IV(n)}{$\IV(n)$}
\newtheorem*{IV(2k+1)}{$\IV(2k+1)$} 
\newtheorem*{V(n)}{$\V(n)$}

\def\H{\operatorname{\mathfrak h}}
\DeclareMathOperator{\E}{\mathcal E}
\DeclareMathOperator{\h}{\mathcal H} 
\DeclareMathOperator{\Supp}{Supp}
\DeclareMathOperator{\st}{St} 
\DeclareMathOperator{\tr}{tr}
\DeclareMathOperator{\Hom}{Hom} 
\DeclareMathOperator{\rend}{End}
\newcommand\ab{\text{ab}} 
\DeclareMathOperator{\Ind}{Ind}
\def\l{\operatorname{\ell}}
\newcommand\orb{\text{orb}}

\begin{document}
\abovedisplayskip=6pt plus3pt minus3pt
\belowdisplayskip=6pt plus3pt minus3pt

\title{Vanishing theorems and conjectures for the\\$\ell^2$--homology
of right-angled Coxeter groups}

\shorttitle{Vanishing theorems and conjectures for homology of Coxeter
groups}

\author{Michael W Davis\\Boris Okun} 
\address{Department of Mathematics, The Ohio State 
University\\Columbus, OH 43210, USA\\\smallskip
\\Department of Mathematics, Vanderbilt 
University\\Nashville, TN 37240, USA}
\asciiaddress{Department of Mathematics, The Ohio State 
University\\Columbus, OH 43210, USA
\\Department of Mathematics, Vanderbilt 
University\\Nashville, TN 37240, USA}

\email{mdavis@math.ohio-state.edu}
\secondemail{okun@math.vanderbilt.edu}
\asciiemail{mdavis@math.ohio-state.edu, okun@math.vanderbilt.edu}

\begin{abstract}
Associated to any finite flag complex $L$ there is a right-angled
Coxeter group $W_L$ and a cubical complex $\Sigma_L$ on which $W_L$
acts properly and cocompactly.  Its two most salient features are that
(1) the link of each vertex of $\Sigma_L$ is $L$ and (2) $\Sigma_L$ is
contractible.  It follows that if $L$ is a triangulation of
$\BS^{n-1}$, then $\Sigma_L$ is a contractible $n$--manifold.  We
describe a program for proving the Singer Conjecture (on the vanishing
of the reduced $\l^2$--homology except in the middle dimension) in the
case of $\Sigma _L$ where $L$ is a triangulation of $\BS^{n-1}$.  The
program succeeds when $n\le 4$. This implies the Charney--Davis
Conjecture on flag triangulations of $\BS^3$.  It also implies the
following special case of the Hopf--Chern Conjecture: every closed
4--manifold with a nonpositively curved, piecewise Euclidean, cubical
structure has nonnegative Euler characteristic.  Our methods suggest
the following generalization of the Singer Conjecture.

{Conjecture:}  If a discrete group $G$ acts properly on a contractible
$n$--manifold, then its $\l^2$--Betti numbers $b_i^{(2)} (G)$ vanish for
$i>n/2$.
\end{abstract}
\asciiabstract{
Associated to any finite flag complex L there is a right-angled
Coxeter group W_L and a cubical complex \Sigma_L on which W_L
acts properly and cocompactly.  Its two most salient features are that
(1) the link of each vertex of \Sigma_L is L and (2) \Sigma_L is
contractible.  It follows that if L is a triangulation of
S^{n-1}, then \Sigma_L is a contractible n-manifold.  We
describe a program for proving the Singer Conjecture (on the vanishing
of the reduced L^2-homology except in the middle dimension) in the
case of \Sigma _L where L is a triangulation of S^{n-1}.  The
program succeeds when n < 5. This implies the Charney-Davis
Conjecture on flag triangulations of S^3.  It also implies the
following special case of the Hopf-Chern Conjecture: every closed
4-manifold with a nonpositively curved, piecewise Euclidean, cubical
structure has nonnegative Euler characteristic.  Our methods suggest
the following generalization of the Singer Conjecture.

{Conjecture:}  If a discrete group G acts properly on a contractible
n-manifold, then its L^2-Betti numbers b_i^{(2)} (G)$ vanish for
i>n/2.
}

\asciikeywords{Coxeter group, aspherical manifold, nonpositive curvature,
L^2-homology, L^2-Betti numbers}

\maketitlepage

\section{\label{0}Introduction}

\begin{The Euler Characteristic Conjecture}\label{0.1}
If $M^{2k}$ is a closed, aspherical manifold of dimension $2k$, then
its Euler characteristic, $\chi (M^{2k})$, satisfies:
$$(-1)^k\chi (M^{2k})\ge 0.$$
\end{The Euler Characteristic Conjecture}

In the special case of Riemannian manifolds of nonpositive sectional
curvature, this conjecture is usually attributed to H Hopf.
  (In this special case, in dimensions $2$ and $4$, the conjecture
follows from the Gauss--Bonnet Theorem. The proof in dimension 4 is
given in Chern's 1956 paper \cite{C}, where it is attributed to
Milnor.)  In the early 1970's, Thurston suggested that the conjecture
might  hold for all closed aspherical manifolds.

In \cite{CD}, R Charney and the first author discuss the Euler
Characteristic Conjecture in the context of piecewise Euclidean
manifolds which are nonpositively curved in the sense of Aleksandrov
and Gromov \cite{G1}.  The case where the manifold is cellulated by
regular Euclidean cubes is particularly easy to discuss. In this case, 
by a lemma of Gromov \cite{G1}, the nonpositive curvature condition
becomes a combinatorial statement:  the link of each vertex must be a
``flag complex''.  (A simplicial complex $L$ is a {\it flag complex}
if any finite nonempty set of vertices, which are pairwise connected
by edges, span a simplex of $L$.)

There is also a combinatorial version of the Gauss--Bonnet Theorem for
a piecewise Euclidean space $X$ (cf \cite{CMS}).  It states that
$\chi (X)$ is the sum over the vertices of $X$ of a local contribution
coming from the link $L$ of a vertex.  In the cubical case, the
formula for the local contribution $\kappa (L)$ coming from a link $L$
is simply, 
$$\kappa (L) = \sum^{\dim L}_{i=-1} \left(-\frac{1}{2}\right)^{i+1}
f_i(L), $$ where $f_i(L)$ denotes the number of $i$--simplices in $L$
and $f_{-1}(L)=1$.  Hence, for piecewise Euclidean cubical manifolds
of nonpositive curvature, the Euler Characteristic Conjecture is
implied by (and, in fact, is equivalent to) the following conjecture
of \cite{CD}.

\begin{The Flag Complex Conjecture}\label{0.2}\qua If $S$ is a flag
triangulation of a\break $(2k-1)$--sphere, then $$(-1)^k\kappa (S)\ge 0, $$
where $\kappa (S)$ is defined by the above formula.
\end{The Flag Complex Conjecture}

In 1976, in \cite{A}, Atiyah introduced the study of $\ell^2$--homology
(or cohomology) into topology.  Here one is interested in the
following situation:  $X$ is either a closed manifold or a finite
$CW$--complex, $\widetilde{X}$ is its universal cover and $\pi$ is its
fundamental group.  For each natural number $i$, one can then define a
Hilbert space, ${\h}_i(\widetilde{X})$, the ``reduced
$\ell^2$--homology'' of $\widetilde{X}$.  There are two methods for
defining this.  In the case where $X$ is a Riemannian manifold, one
lifts the metric to $\widetilde{X}$ and then defines (de Rham)
$\ell^2$--cohomology by using differential forms with square integrable
norms.  When $X$ is a finite $CW$--complex, one lifts the cell
structure to $\widetilde{X}$ and then defines ${\mathcal
H}_i(\widetilde{X})$ by using infinite cellular chains with square
summable coefficients.  In either case, the Hilbert space ${\mathcal
H}_i ( \widetilde{X})$ comes equipped with an orthogonal $\pi$--action.
When $X$ is a triangulated Riemannian manifold, the equivalence of the
two definitions was proved by Dodziuk in \cite{Do1}.  In this paper, 
we will deal only with the cellular version of $\ell^2$--homology.

A key feature of the $\ell^2$--theory is that, by using the
$\pi$--action, it is possible to attach to the Hilbert space ${\mathcal
H}_i (\widetilde{X})$ a nonnegative real number, called the
``$i^{\text{th}}~\ell^2$--Betti number''.  (This is explained in
Section 3.)  A formula of Atiyah \cite{A} states that the alternating
sum of these $\ell^2$--Betti numbers is the ordinary Euler
characteristic $\chi (X)$.  (The precise statement of Atiyah's Formula
can be found in Section~\ref{3.3} of this paper.)

Shortly after this formula became known, Dodziuk and Singer pointed
out that Atiyah's Formula shows that the Euler Characteristic
Conjecture follows if one can prove that the reduced $\ell^2$--homology
of the universal cover of any even dimensional, closed, aspherical
manifold vanishes except in the middle dimension.  (This is explained
in the introduction of \cite{Do2}.)  This led to the following
conjecture.

\begin{The Vanishing Conjecture}\label{0.3}  If $M^n$ is a closed
aspherical manifold, then 
$${\h}_i(\widetilde{M}^n)=0 \quad \text{  
for all\qua $i\not= \frac{n}{2}$. }$$
\end{The Vanishing Conjecture}

Singer's Conjecture holds for elementary reasons in dimensions $\le
2$.  In \cite{LL} Lott and L\"uck proved that it holds for those
aspherical $3$--manifolds for which Thurston's Geometrization
Conjecture is true.  It is also known to hold for (a) locally
symmetric spaces, (b) negatively curved K\"ahler manifolds (by
\cite{G2}), (c) Riemannian manifolds of sufficiently pinched negative
sectional curvature (by \cite{DX}), (d)  closed aspherical manifolds
with  fundamental group containing an infinite amenable normal
subgroup (by \cite{CG2}), and (e) manifolds
which fiber over $\BS^1$ (by \cite{L1}).

We note that the Euler Characteristic Conjecture and Singer's
Conjecture both make sense for closed aspherical orbifolds or, for
that matter, for virtual Poincar\'e duality groups. 

In several earlier papers (eg, \cite{CD}, \cite{D1}, \cite{D2}, 
\cite{D3}, or \cite{DM}), the first author has described a
construction which associates to any finite flag complex $L$, a
``right-angled'' Coxeter group $W_L$ and a cubical cell complex
$\Sigma_L$ on which $W_L$ acts properly and cocompactly.  (The details
of this construction will be given in Sections 5 and 6, below.)  Its
two most salient features are that (1) the link of each vertex of
$\Sigma_L$ is isomorphic to $L$ and (2) $\Sigma_L$ is contractible.

If $\Gamma$ is a torsion-free subgroup of finite index in $W_L$, then
$\Gamma$ acts freely on $\Sigma_L$ and $\Sigma_L/\Gamma$ is a finite
complex.  By (2), $\Sigma_L/\Gamma$ is aspherical.  If $L$ is
homeomorphic to the $(n-1)$--sphere, then by (1), $\Sigma_L$ is an
$n$--manifold.  Hence, this construction gives many examples of closed
aspherical manifolds.  Singer's Conjecture for such manifolds becomes
the following.

\begin{uConjecture}\label{0.4}  Suppose $S$ is a triangulation of the
$(n-1)$--sphere as a flag complex.  Then
$${\h}_i(\Sigma_S)=0 \quad \text{  for all\qua $i\not= \frac{n}{2}$. }$$
\end{uConjecture}

The purpose of this paper is to describe a partially successful
program for proving this conjecture by using standard techniques of
algebraic topology and induction on the dimension $n$.  Our main
result, Theorem~\ref{9.3.1}, is that the program succeeds in half the
cases:  if Conjecture~\ref{0.4} is true in some odd dimension $n$, 
then it is also true in dimension $n+1$. Moreover, in odd dimensions it 
is only necessary to establish a weak form of the conjecture. 

As we shall see in Section~\ref{10}, the Geometrization Conjecture is
true for the $3$--manifolds which we are considering.  Hence, the
Lott--L\"uck result implies that Conjecture~\ref{0.4} is true for $n=3$
and, therefore, also for $n=4$.  This gives the following
(Theorem~\ref{11.1.1} of Section~\ref{11}).

\begin{uTheorem}  Conjecture~\ref{0.4} is true for $n\le 4$.
\end{uTheorem}

Hence, $4$--manifolds of the form $\Sigma_S/\Gamma$ have nonnegative
Euler characteristic.  As explained in \cite{CD} and~\ref{6.3.4},
below,  this implies the following (Theorem~\ref{11.2.1}).

\begin{uTheorem}  The Flag Complex Conjecture is true in dimension $3$.  In
other words, if $S$ is a triangulation of a homology $3$--sphere as a
flag complex, then
$$\kappa (S)\ge 0.$$
\end{uTheorem}

The combinatorial Gauss--Bonnet Theorem then implies the next result
(Theorem~\ref{11.2.2}).

\begin{uTheorem}  The Euler Characteristic Conjecture holds true for all
nonpositively curved, piecewise Euclidean $4$--manifolds which are
cellulated by regular Euclidean cubes.  In other words, for any such
$4$--manifold $M^4$, 
$$\chi (M^4)\ge 0.$$
\end{uTheorem}

A surprising aspect of our analysis is that it turns out that
Conjecture~\ref{0.4} is equivalent to a statement about the vanishing
of ${\h}_i(\Sigma_L)$ for an arbitrary finite flag complex $L$
(not necessarily a sphere).  More precisely, we will show in Section
9, that Conjecture~\ref{0.4} is equivalent to the following.

\begin{uConjecture}\label{0.5}
  Suppose $L$ is a finite flag complex.  If $L$ can be embedded as a
full subcomplex of some flag triangulation of the $2k$--sphere, then
$${\h}_i(\Sigma_L)=0 \quad \text{  for all $i>k$. }$$
\end{uConjecture}

Let us say that an $n$--dimensional polyhedron $X$ has {\it spherical
links in codimensions $\le m$} if, for $i\le m$, the link of any
$(n-i)$--cell in $X$ is an $(i-1)$--sphere.  For example, if $m=1$, then
$X$ is a pseudomanifold, while if $m=n$, then $X$ is a manifold.  The
inductive arguments of Section 9 suggest the following generalization
of Singer's Conjecture.

\begin{uConjecture}\label{0.6}
Suppose  an $n$--dimensional aspherical polyhedron $X$ has \linebreak spherical
links in codimensions $\le 2l +1$, where $2l +1\le n$.  Then
$${\h}_{n-i} (\widetilde{X})=0 \quad \text{  for $i\le l$. }$$
\end{uConjecture}

When $l =0$ (so that $X$ is a pseudomanifold), this conjecture holds
for elementary reasons (as we explain in~\ref{2.6}).  For right-angled
Coxeter groups the conjecture reads as follows.

\begin{uConjecture}\label{0.7}
Suppose an $(n-1)$--dimensional flag complex $L$ has spherical links in
codimensions $\le 2l +1$, where $2l +1\le n$.  (If $n=2l +1$, we take
this to mean that $L$ is an $(n-1)$--sphere.)  Then, for $i\le l$,  
$${\h}_{n-i}(\Sigma_L)=0.$$
\end{uConjecture}

In Theorem~\ref{9.3.3}, we show that if Conjecture~\ref{0.4} is true
for  $n=2k +1$, then Conjecture~\ref{0.7} holds for $l=k$ and any
$n\ge 2k +1$.  In particular, since Conjecture~\ref{0.4} holds for
$n=3$ we get the following (Theorem~\ref{11.3.2} in Section~\ref{11}).

\begin{uTheorem}  Suppose $S$ is a flag triangulation of an $(n-1)$--sphere, 
$n\ge 3$.  Then
$${\h}_i(\Sigma_S)={\h}_{n-i}(\Sigma_S)=0 \quad \text{  for
$i=0, 1$. }$$
\end{uTheorem}

Conjecture~\ref{0.5}, taken together with recent work of Bestvina, Kapovich and
Kleiner~\cite{BKK}, suggests the following different generalization of Singer's
Conjecture (Conjecture~\ref{8.9.1} of Section 8).

\begin{uConjecture}\label{0.8}  Suppose that a discrete group $G$ acts properly
on a contractible $n$--manifold.  Then
$$b^{(2)}_i(G)=0 \quad \text{  for $i>\frac{n}{2}$. }$$  
\rm(See \ref{3.3.7} for the
definition of the $\ell^2$--Betti numbers  $b^{(2)}_i(G)$.)
\end{uConjecture}

In the last three sections (\ref{12}, \ref{13} and \ref{14}) we discuss some 
possible attacks on (a weak form of) Conjecture~\ref{0.4} in odd dimensions.

Some of the results of this paper appeared in a preliminary form in
\cite{DM}.\nl
We thank the referee for some useful suggestions.\nl
The authors were partially supported by NSF grants.

\section{\label{1}Group actions on $CW$ complexes}

\subsection{Geometric $\boldsymbol G$--complexes}\label{1.1}
Let $G$ be a discrete group.  A {\it $G$--complex} is a $CW$ complex
$X$ together with a cellular action of $G$ on $X$.  All $G$--complexes
in this paper will be {\it geometric}.  By this we mean that the
$G$--action is proper (ie, that each cell stabilizer is finite) and
cocompact (ie, that $X/G$ is compact).

\subsection{Regular complexes and orbihedra}\label{1.2}
 A $CW$ complex $X$ is {\it regular}
if the characteristic map of each cell is an embedding (so that the
boundary of each cell is an embedded sphere).  If $X$ is a geometric
$G$--complex and if it is regular, then $X/G$ is an {\it orbihedron} in
the sense of \cite{H}.  The structure of an orbihedron encodes not
only the topological space $X/G$, but also the isomorphism types of
the cell stabilizers for each $G$--orbit of cells.  If $H$ is a
subgroup of $G$, then the natural projection $X/H\to X/G$ is an {\it
orbihedral covering map}.

\subsection{The orbihedral Euler characteristic}\label{1.3}
Suppose $X$ is a geometric $G$--complex.  Then there are only a finite
number of $G$--orbits of cells in $X$ and the order of each cell
stabilizer is finite.  The orbihedral Euler characteristic of $X/G$, 
denoted $\chi^{\orb}(X/G)$, is the rational number defined by

\subsubsection{}\label{1.3.1}  
\myequation{ \chi^{\orb}(X/G ) =  \sum_\sigma
\frac{(-1)^{\dim\sigma}}{\vert G_\sigma\vert},}

 where the summation is over a set of representatives for the
$G$--orbits of cells and where $\vert G_\sigma\vert$ denotes the order
of the stabilizer $G_\sigma$ of $\sigma$.

\subsubsection{}\label{1.3.2} If $G$ acts freely on $X$, then 
$\chi^{\orb}(X/G )$
is the ordinary Euler characteristic of the finite $CW$ complex $X/G$.

\subsubsection{}\label{1.3.3} If $H$ is a subgroup of finite index $m$ in $G$, 
then it follows immediately from the definition that
$$\chi^{\orb}(X/H)=m\chi^{\orb}(X/G).$$

\subsection{Universal spaces for proper $\boldsymbol G$--actions}\label{1.4}  A
$G$--complex $X$ is a {\it universal space for proper $G$--actions} if
the action is proper and if the fixed point set $X^F$ is contractible
for each finite subgroup $F$ of $G$.  (In particular, taking $F$ to be
the trivial subgroup, this means that $X$ is contractible.)  Such
universal spaces always exist and are unique up to $G$--equivariant
homotopy equivalence.  It is often denoted by $\underline{EG}$.  If, 
in addition, the action is cocompact, then
$\chi^{\orb}(\underline{EG}/G)$ is defined and is an invariant of $G$.
It is the {\it Euler characteristic of $G$}, $\chi (G)$, in the sense
of \cite{W}.

\section{\label{2}$\ell^2$--homology}

We review some basic facts about the $\ell^2$--homology of geometric
$G$--com\-plexes.  References for this material include \cite{CG2}, 
\cite{Do1}, \cite{G3}, and \cite{E} (which is particularly easy to
read).

\subsection{Square summable functions}\label{2.1}  Suppose $G$ is a
countable discrete group.  Let $\ell^2 (G)$ denote the vector space of
real-valued, square-summable functions on $G$, ie, $\ell^2 (G)=\{
f\colon G\to\BR\mid\sum f(g)^2<\infty\}$.  It is a Hilbert space:  the
inner product is given by
$$\langle f_1, f_2\rangle = \sum_{g\in G} f_1(g)f_2(g).$$ The group
ring $\BR G$ can be identified with the dense subspace of $\ell^2 (G)$
consisting of the functions with finite support.

The action of $G$ on itself by left translation induces an orthogonal
(left) $G$--action on $\ell^2(G)$.  (There is also an orthogonal right
$G$--action on $\ell^2(G)$ induced by right translation.)

\subsection{Hilbert $\boldsymbol G$--modules}\label{2.2}  Given a natural 
number $n$, 
let $\ell^2(G)^n$ denote the direct sum of a $n$ copies of
$\ell^2(G)$, equipped with the diagonal (left) $G$--action.  A Hilbert
space $V$ with orthogonal $G$--action is a {\it Hilbert $G$--module} if
it is isomorphic to a closed, $G$--stable subspace of $\ell^2(G)^n$, 
for some $n\in\BN$.  (In the literature, this is sometimes called a
``finitely generated'' Hilbert $G$--module or a Hilbert $G$--module of
``finite type''.)

\subsubsection{}\label{2.2.1}  
If $F$ is a finite subgroup of $G$, then $\ell^2(G/F)$, 
the space of square summable functions on $G/F$, can be identified
with the subspace of $\ell^2(G)$ consisting of the square summable
functions on $G$ which are constant on each coset.  This subspace is
clearly closed and $G$--stable; hence, $\ell^2(G/F)$ is a Hilbert
$G$--module.

\subsubsection{}\label{2.2.2}  
A {\it map} of Hilbert $G$--modules is a $G$--equivariant, 
bounded linear map.

The complication which arises at this point is that the image of such
a map need not be a closed subspace.  This leads to the notions of a
``weakly'' exact sequence and a ``weak'' isomorphism, defined below.

\subsubsection{Weak exactness}\label{2.2.3}   A sequence
$U\xrightarrow{e} V\xrightarrow{f}  W$ of maps of Hilbert
$G$--mo\-dules is {\it weakly exact} at $V$ if the closure of the image
of $e$ (denoted $\overline{\im e}$) is the kernel of $f$ (denoted
$\Ker f$).  Similarly, $e\colon U\to V$ is {\it weakly surjective} if
$\overline{\im e} = V$ and it is a {\it weak isomorphism} if it is
injective and weakly surjective.

\subsubsection{}\label{2.2.4}  
If two Hilbert $G$--modules are weakly isomorphic, then
they are $G$--isometric (Lemma 2.5.3 in \cite{E}).

\subsubsection{Induced representations}\label{2.2.5}  Suppose $H$ is a
subgroup of $G$ and that $W$ is a Hilbert $H$--module.  The {\it
induced representation}, $\Ind^G_H(W)$, can be defined as the
$\ell^2$--completion of $\BR G\otimes_{\BR H}W$.  Alternatively, it is
the vector space of all square summable sections of the vector bundle
$G\times_HW\to G/H$.  (Here $G/H$ is discrete.)  The induced
representation is obviously a Hilbert space with orthogonal
$G$--action.  If $W$ is a closed subspace of $\ell^2(H)^n$, then
$\Ind^G_H(W)$ is a closed subspace of $\ell^2 (G)^n$.  (This follows
from the observation that $\Ind^G_H(\ell^2(H))$ can be identified with
$\ell^2 (G)$.)  Thus, $\Ind^G_H(W)$ is a Hilbert $G$--module. For
example, if $F$ is a finite subgroup of $G$ and $\BR$ denotes the
trivial $1$--dimensional representation of $F$, then $\Ind^G_F(\BR)$
can be identified with $\ell^2 (G/F)$.

\subsection{$\boldsymbol \ell^{\boldsymbol 2}$--homology and 
cohomology}\label{2.3}  
Given a geometric
$G$--complex $X$, let $K_\ast (X)$ denote the usual cellular chain
complex on $X$, regarded as a left $\BZ(G)$--module.  
(We use this notation since we want to reserve
$C_\ast (X)$ for the chain complex of $\ell^2$--chains on $X$.)

\subsubsection{$\ell^2$--chains}\label{2.3.1}  Set
$$C_i(X)= \ell^2(G) \otimes_{\BZ G} K_i(X)$$ 
where $\l^2(G)$ is regarded as a right $\BZ G$--module. 
An element of $C_i(X)$ is
an $\ell^2${\it --chain}; it is an infinite chain with square summable
coefficients.  The Hilbert space $C_i(X)$ can also be regarded as the
space of $\ell^2$--cochains on $X$.

\subsubsection{}\label{2.3.2}  
If $\sigma$ is an $i$--cell of $X$, then the space of
$\ell^2$--chains which are supported on the $G$--orbit of $\sigma$ can
be identified with $\ell^2(G/G_\sigma )$.  Since there are a finite
number of such orbits, $C_i(X)$ is the direct sum of a finite number
of such subspaces.  Hence, by~\ref{2.2.1}, $C_i(X)$ is a Hilbert
$G$--module.

\subsubsection{Unreduced and reduced $\ell^2$--homology}\label{2.3.3}  
We define the
boundary map $d_i\colon C_i(X)\to C_{i-1} (X)$ and the coboundary map 
$\delta^i\colon C_i(X)\to C_{i+1}(X)$
by the usual formulae. Then the
boundary and the coboundary maps are  $G$--equivariant, 
bounded linear maps.  The coboundary map $\delta^i$ can be identified 
with $d^\ast_{i+1}$ (the adjoint of
$d_{i+1}$).  Define subspaces of $C_i(X)$:
\begin{align*}
Z_i(X) &= \Ker d_i & Z^i(X)&=\Ker \delta^i\\ B_i(X) &= \im d_{i+1} &
B^i(X)&=\im \delta^{i-1} 
\end{align*}
the $\ell^2${\it --cycles, --cocycles, --boundaries} and {\it
--coboundaries}, respectively.  The corresponding quotient spaces
$$H^{(2)}_i(X) = Z_i(X)/B_i(X)$$ and
$$H^i_{(2)}(X) = Z^i(X)/B^i(X)$$ are the {\it unreduced} $\ell^2${\it
--homology} and {\it --cohomology groups}, respectively. (In other words, 
$H^{(2)}_i(X)$ is the ordinary equivariant homology of $X$ with 
coefficients in $\l^2(G)$, ie, $H^{(2)}_i(X)=H^{G}_i(X;\l^2(G) )$.)   
Since the 
subspaces $B_i(X)$ and $B^i(X)$ need not be closed, these quotient
spaces need not be isomorphic to Hilbert spaces.

Let $\overline{B_i}(X)$ (respectively, $\overline{B^i}(X)$) denote the
closure of $B_i(X)$ (respectively, $B^i(X)$).  The {\it reduced
$\ell^2$--homology} and {\it --cohomology groups} are defined by:
$${\h}_i(X) = Z_i(X)/\overline{B_i}(X)$$
$${\h}^i (X) = Z^i(X)/\overline{B^i}(X) .$$ They are Hilbert
$G$--modules (since each can be identified with the orthogonal
complement of a closed $G$--stable subspace in a closed $G$--stable
subspace of $C_i(X)$).

\subsubsection{Hodge decomposition}\label{2.3.4}  Since $\langle\delta^{i-1}
(x), y\rangle = \langle x, d_i(y)\rangle$ for all $x\in$\break $C_{i-1}(X)$ and
$y\in C_i(X)$, we have orthogonal direct sum decompositions:
\begin{align*}
C_i(X)=\overline{B_i}(X)\oplus Z^i(X) \\ \intertext{and}
C_i(X)=\overline{B^i}(X)\oplus Z_i(X).
\end{align*} 
 Since $\langle \delta^{i-1}(x), d_{i+1}(y)\rangle =\langle
x, d_id_{i+1}(y)\rangle =0$, the subspaces $\overline{B_i}(X)$ and
$\overline{B^i}(X)$ are orthogonal.  Hence, 
$$C_i(X)=\overline{B_i}(X)\oplus\overline{B^i}(X)\oplus(Z_i(X)\cap
Z^i(X)).$$ It follows that the reduced $\ell^2$--homology and
--cohomology groups can both be identified with the subspace
$Z_i(X)\cap Z^i(X)$.  We denote this intersection again by ${\mathcal
H}_i(X)$ and call it the subspace of {\it harmonic} $i$--cycles.  Thus, 
an $i$--chain is harmonic if and only if it is simultaneously a cycle
and a cocycle.

The {\it combinatorial Laplacian} $\Delta\colon C_i(X)\to C_i(X)$ is
defined by $\Delta =\delta^{i-1} d_i+d_{i+1}\delta^i$.  One checks
that ${\h}_i(X)=\Ker \Delta$.

\subsubsection{Relative groups}\label{2.3.5}  If $Y$ is a $G$--stable
subcomplex of $X$, then $(X, Y)$ is a {\it pair of geometric
$G$--complexes}.  The reduced $\ell^2$--homology (or --cohomology) groups
${\h}_i(X, Y)$ are then defined in the usual manner.

\subsection{Basic algebraic topology}\label{2.4}
Suppose $(X, Y)$ is a pair of geometric $G$--complexes.  Versions of
most of the Eilenberg--Steenrod homology theory Axioms hold for ${\h}_\ast
(X, Y)$.  We list some standard properties below.  (Of course, similar
results hold for the contravariant $\ell^2$--cohomology functor.)

\subsubsection{Functoriality}\label{2.4.1}  For $i=1, 2$, 
 suppose $(X_i, Y_i)$ is a pair of geometric $G$--complexes and that
$f\colon (X_1, Y_1)\to (X_2, Y_2)$ is a $G$--equivariant map (a $G${\it
--map} for short).  Then there is an induced map $f_\ast\colon{\mathcal
H}_i(X_1, Y_1)\to{\h}_i(X_2, Y_2)$ and this gives a functor from
pairs of $G$--complexes to Hilbert $G$--modules.  Moreover, if $f'
\colon (X_1, Y_1)\to (X_2, Y_2)$ is another $G$--map which is homotopic
to $f$ (not necessarily $G$--homotopic), then $f_\ast = f'_\ast$.

\subsubsection{Exact sequence of a pair}\label{2.4.2}  The sequence of a
pair $(X, Y)$, 
$$\to {\h}_i(Y)\to{\h}_i(X)\to{\h}_i(X, Y)\to$$
is weakly exact.

\subsubsection{Excision}\label{2.4.3}  Suppose that $(X, Y)$ is a pair of
geometric $G$--complexes and that $U$ is a $G$--stable subset of $Y$
such that $Y-U$ is a subcomplex.  Then the inclusion $(X-U, Y-U)\to
(X, Y)$ induces an isomorphism:
$${\h}_i(X-U, Y-U)\cong {\h}_i (X, Y) .$$

A standard consequence of the last two properties is the following.

\subsubsection{Mayer--Vietoris sequences}\label{2.4.4}  Suppose $X=X_1\cup
X_2$, where $X_1$ and $X_2$   are  $G$--stable subcomplexes of $X$.
Then $X_1\cap X_2$ is also $G$--stable and the Mayer--Vietoris sequence, 
$$\to {\h}_i(X_1\cap X_2)\to{\h}_i(X_1)\oplus{\mathcal
H}_i(X_2)\to {\h}_i(X)\to$$ is weakly exact.

\subsubsection{Twisted products and the induced representation}\label{2.4.5}
Suppose that $H$ is a subgroup of $G$ and that $Y$ is a space on which
$H$ acts.  The {\it twisted product}, $G\times_H Y$, is the quotient
space of $G \times Y$ by the $H$--action defined by
$h(g, y)=(g h^{-1}, h y)$.  It is a left $G$--space and a $G$--bundle over
$G/H$.  Since $G/H$ is discrete, $G \times_H Y$ is a disjoint union of
copies of $Y$, one for each element of $G/H$.  If $Y$ is a geometric
$H$--complex, then $G\times_H Y$ is a geometric $G$--complex and the
following formula obviously holds:
$${\h}_i(G\times_H Y)\cong \Ind^G_H({\h}_i(Y)) .$$

\subsubsection{K\"unneth Formula}\label{2.4.6}  Suppose $G=G_1\times G_2$ and
that for $j=1, 2$, $X_j$ is a geometric $G_j$--complex.  Then $X_1\times
X_2$ is a geometric $G$--complex and
$${\h}_k(X_1\times X_2)\cong \sum_{i+j=k}{\mathcal
H}_i(X_1)\widehat{\otimes} {\h}_j(X_2) , $$ where
$\widehat{\otimes}$ denotes the completed tensor product.

\subsection{Homology in dimension $\boldsymbol 0$}\label{2.5}
An element of $C_0(X)$ is an $\ell^2$ function on the set of vertices
of $X$; it is a $0$--cocycle if and only if it takes the same value on
the endpoints of each edge.  Hence, if $X$ is connected, any
$0$--cocycle is constant.  If, in addition, $G$ is infinite (so that
the $1$--skeleton of $X$ is infinite), then this constant must be $0$.
So, when $X$ is connected and $G$ is infinite, $H^0_{(2)}(X)={\mathcal
H}^0(X)=0$.  Hence, 
\subsubsection{}\label{2.5.1}
\myequation{{\h}_0(X) =0.}

\subsubsection{}\label{2.5.2}  
On the other hand, the unreduced homology $H^{(2)}_0
(X)$ need not be $0$.  For example, if $X=\BR$, cellulated as the
union of intervals $[n, n+1]$, and $G=\BZ$, then any vertex of $\BR$ is
an $\ell^2$--$0$--cycle which is not $\ell^2$--boundary.  (A vertex
bounds a half-line which can be thought of as an infinite $1$--chain
but this $1$--chain is not square summable.)  In fact, if $G$ is
infinite, then a theorem of Kesten~\cite{Ke} 
implies that $H^{(2)}_0(X)=0$ if and only if $G$ is not amenable.

\subsection{The top-dimensional homology of a
pseudomanifold}\label{2.6}$\phantom1$\kern-.3em
Suppose that an $n$--dimensional, regular $G$--complex $X$ is a
pseudomanifold.  This means that each $(n-1)$--cell is contained in
precisely two $n$--cells.  If a component of the complement of the
$(n-2)$--skeleton is not orientable, then it does not support a nonzero
$n$--cycle (with coefficients in $\BR$).  If such a component is
orientable, then any $n$--cycle supported on it is a constant multiple
of the $n$--cycle with all coefficients equal to $+1$.  If the component has
an infinite number of $n$--cells, then this $n$--cycle does not have
square summable coefficients.  Hence, if each component of the
complement of the $(n-2)$--skeleton is either infinite or
nonorientable, then $H^{(2)}_n(X)=0$.  In particular, if the
complement of the $(n-2)$--skeleton is connected and if $G$ is
infinite, then $H^{(2)}_n(X)=0$.

\subsection{Poincar\'e duality}\label{2.7}
Suppose $(X, \partial X)$ is a pair of geometric $G$--com\-plexes and that
$X$ is an $n$--dimensional manifold with boundary.  Then
\subsubsection{}\label{2.7.1}
\myequation{ {\h}_i(X, \partial X) \cong{\h}^{n-i}(X) \qquad
\text{and}}

\subsubsection{}\label{2.7.2}\myequation{ {\h}_i(X)\cong
{\h}^{n-i}(X, \partial X).}

In the case where $X$ is cellulated as a $PL$ manifold with boundary, 
these isomorphisms  are induced by the bijective
correspondence $\sigma\leftrightarrow D\sigma$ which associates to
each $i$--cell $\sigma$ its dual $(n-i)$--cell $D\sigma$.  A slight
elaboration of this argument also works in the case where $(X, \partial
X)$ is a polyhedral homology manifold with boundary; the only
complication being that the ``dual cells'' need not actually be cells, 
rather they are ``generalized homology disks'' as defined in
Section~\ref{4.3}, below.

\subsubsection{}\label{2.7.3} 
In fact, as is shown in \cite[Theorem 3.7.2]{E}, in order
to have the Poincar\'e duality isomorphisms of~\ref{2.7.1}, all one
need assume is that $(X, \partial X)$ is a ``virtual $PD^n$--pair''.
This means that there is a subgroup $H$ of finite index in $G$ so that
the chain complexes $K_\ast (X, \partial X)$ and $^n\!DK_\ast(X)$ are
chain homotopy equivalent, where $^n\!DK_i(X)$ is defined by
$^n\!DK_i(X)=\Hom_{\BZ H}(K_{n-i}(X), \BZ H)$.

\subsection{Extended $\boldsymbol \ell^{\boldsymbol 2}$--homology}\label{2.8}
In \cite{F}, Farber defines an ``extended $\ell^2$--(co)\-ho\-mo\-logy''
theory and demonstrates that this is the correct categorical framework
for $\ell^2$--homology.  An extended $\ell^2$--homology object is
isomorphic to the sum of its ``projective part'' and its ``torsion
part''.  The projective part is essentially the reduced
$\ell^2$--homology group while its torsion part contains information
such as Novikov--Shubin invariants.  Since we have nothing to say about
this torsion part, we shall stick to the simpler reduced
$\ell^2$--homology groups.

\section{\label{3}$\ell^2$--Betti numbers}

The feature which distinguishes $\ell^2$--homology from its brothers, 
the $\ell^p$--ho\-mology theories, is that one can associate to each Hilbert
$G$--module a nonnegative real number called its ``von Neumann
dimension''.

\subsection{von Neumann algebra}\label{3.1}
The {\it von Neumann algebra} ${\mathcal N} (G)$ associated to $G$ is
the algebra of all $G$--equivariant, bounded linear endomorphisms of
$\ell^2 (G)$.  Since $\ell^2(G)$ is also a right $\BR G$--module we see
that $\BR G\subset{\mathcal N}(G)$.  In fact, ${\mathcal N}(G)$ is the
weak closure of $\BR G$ in the space $\rend (\ell^2(G))$ of all
bounded linear endomorphisms of $\ell^2 (G)$.

For each $g\in G$, let $e_g$ denote the characteristic function of $\{
g\}$, ie, $e_g(h)=0$ if $h\not= g$ and $e_g(h)=1$ if $h=g$.  Then
$\{ e_g\}_{g \in G}$ is a basis for $\BR G$ and an orthonormal basis
for the Hilbert space $\ell^2(G)$.

Define a linear functional $\tr_G\colon{\mathcal N}(G)\to\BR$ by
\subsubsection{}\label{3.1.1}
\myequation{\tr_G(\varphi ) =\langle\varphi (e_1), e_1\rangle.}

(The restriction of $\tr_G$ to the subset $\BR G$ is the classical
Kaplansky trace.)

Next, suppose that $\varphi$ is a $G$--equivariant, bounded linear
endomorphism of $\ell^2(G)^n$, $n\in\BN$.  Then $\varphi$ can be
represented as an $n$ by $n$ matrix $(\varphi_{ij})$ with coefficients
in ${\mathcal N}(G)$.  Define
\subsubsection{}\label{3.1.2}
\myequation{\tr_G(\varphi ) = \sum^n_{i=1} \tr_G (\varphi_{ii}) .}

  The standard argument shows that $\tr_G(\varphi )$ depends only on
the conjugacy class of $\varphi$.

\subsection{von Neumann dimension}\label{3.2}
Let $V$ be a Hilbert $G$--module.  Choose an embedding of $V$ as a
closed $G$--stable subspace of $\ell^2 (G)^n$ for some $n\in\BN$.  Let
$p_V\colon \ell^2(G)^n\to\ell^2(G)^n$ denote orthogonal projection
onto $V$.  The {\it von Neumann dimension} of $V$, denoted by $\dim_G
(V)$, is defined by
\subsubsection{}\label{3.2.1}
\myequation{\dim_G (V) = \tr_G (p_V) .}

  Standard arguments (as in \cite{E}) show that this definition is
independent of the choice of embedding $V\to\ell^2(G)^n$.

We list some properties of $\dim_G (V)$.  Proofs can be found in
\cite{E}.

\subsubsection{}\label{3.2.2}  \myequation{\dim_G (V)\in [0, \infty ).}

\subsubsection{}\label{3.2.3}  
\myequation{\dim_G (V)=0 \text{ if and only if }V=0.}

\subsubsection{}\label{3.2.4}  If $G$ is the trivial group (so that the Hilbert
space $V$ is finite dimensional), then $\dim_G(V)=\dim (V)$.

\subsubsection{}\label{3.2.5} \myequation{ \dim_G (\ell^2(G))=1.}

\subsubsection{}\label{3.2.6}  
\myequation{\dim_G (V_1\oplus V_2)=\dim_G (V_1)+\dim_G(V_2).}

\subsubsection{}\label{3.2.7}  
If $f\colon V\to W$ is a map of Hilbert $G$--modules, 
then by~\ref{2.2.2} and~\ref{3.2.6}, 
$$\dim_G (V)=\dim_G (\Ker f)+\dim_G (\overline{\im f}) .$$

\subsubsection{}\label{3.2.8} 
 If $f\colon V\to W$ is a map of Hilbert $G$--modules and
$f^\ast\colon W\to V$ is its adjoint,  then $\Ker f$ and
$\overline{\im f^*}$ are orthogonal complements in $V$.  Hence, 
$$\dim_G(V)=\dim_G (\Ker f)+\dim_G (\overline{\im f^*}) .$$ So, 
by~\ref{3.2.7} $$\dim_G (\overline{\im f})=\dim_G (\overline{\im
f^*}).$$

\subsubsection{}\label{3.2.9}  
By~\ref{3.2.6} and~\ref{3.2.7}, if $0\to V_n\to\cdots\to V_0\to 0$
is a weakly exact sequence of Hilbert $G$--modules, then
$$\sum^n_{i=0} (-1)^i\dim_G (V_i)=0 .$$

\subsubsection{}\label{3.2.10}  
If $H$ is a subgroup of finite index $m$ in $G$, 
then
$$\dim_H(V)=m\dim_G (V) .$$

Combining \ref{3.2.10} with \ref{3.2.4} we get the following.

\subsubsection{}\label{3.2.11}  If $G$ is finite, then
$$\dim_G (V)= \frac{1}{\vert G\vert}\dim (V) .$$

\subsubsection{}\label{3.2.12}  
If $H$ is a subgroup of $G$ and $W$ is a Hilbert
$H$--module, then
$$\dim_G (\Ind^G_H(W))=\dim_H(W) .$$

\subsubsection{}\label{3.2.13}  
If $F$ is a finite subgroup of $G$, then by~\ref{2.2.5}
and \ref{3.2.12}, $$\dim_G (\ell^2(G/F))=\frac{1}{\vert F\vert}.$$

\subsubsection{}\label{3.2.14} 
Suppose $G=G_1\times G_2$ and that for $j=1, 2$, 
$V_j$ is a Hilbert $G_j$--module.  Then $V_1\widehat{\otimes} V_2$ is a
Hilbert $G$--module and
$$\dim_G (V_1\widehat{\otimes} V_2)= \dim_{G_1}(V_1)\dim_{G_2}(V_2) .$$

\subsection{$\boldsymbol \ell^{\boldsymbol 2}$--Betti numbers}\label{3.3}
Given a pair $(X, Y)$ of geometric $G$--complexes, its {\it $i^{\text{th}}$
$\ell^2$--Betti number}, $b^{(2)}_i(X, Y;G)$, is defined by
\subsubsection{}\label{3.3.1}
\myequation{ b^{(2)}_i(X, Y;G) = \dim_G ({\h}_i(X, Y)) .}

From the properties of von Neumann dimension in~\ref{3.2} and the
properties of reduced $\ell^2$--homology in Section 2, we get
properties of $\ell^2$--Betti numbers.  We list a few of these
properties below.

\subsubsection{}\label{3.3.2}  
\myequation{b^{(2)}_i(X, Y;G)=0 \text{ if and only if } {\h}_i(X, Y)=0 
\text{ (by~\ref{3.2.3})}.}

\subsubsection{}\label{3.3.3}  
If $H$ is a subgroup of finite index $m$ in $G$, 
then, by~\ref{3.2.10}, 
$$b^{(2)}_i(X, Y;H)=mb^{(2)}_i(X, Y;G) .$$

\subsubsection{}\label{3.3.4}  
By~\ref{2.4.5} and~\ref{3.2.12}, for any geometric $H$--complex
$Y$, with $H\subset G$, 
$$b^{(2)}_i(G\times_H Y;G)=b^{(2)}_i(Y;H) .$$

\subsubsection{K\"unneth Formula}\label{3.3.5}  If $G=G_1\times G_2$ and for
$j=1, 2$, $X_j$ is a geometric $G_j$--complex, then by~\ref{2.4.6}
and~\ref{3.2.14}, 
$$b^{(2)}_k(X_1\times
X_2;G)=\sum_{i+j=k}b^{(2)}_i(X_1;G_1)b_j(X_2;G_2) .$$

\subsubsection{Atiyah's Formula}\label{3.3.6}  
By~\ref{1.3.1}, \ref{2.2.1} and~\ref{3.2.13}, 
$$\chi^{\orb}(X/G) = \sum \frac{(-1)^{\dim\sigma}}{\vert
G_\sigma\vert} = \sum^{\dim X}_{i=0} (-1)^i\dim_G (C_i(X)) .$$ A
standard argument (given in \cite[Theorem 3.6.1]{E}) then proves
Atiyah's Formula:
$$\chi^{\orb}(X/G)= \sum^{\dim X}_{i=0} (-1)^i b^{(2)}_i(X;G) .$$

\subsubsection{$\ell^2$--Betti numbers of a group}\label{3.3.7}  
As in~\ref{1.4}, 
let $\underline{EG}$ denote the universal space for proper
$G$--actions.  Also, assume that $\underline{EG}/G$ is compact (so that
$\underline{EG}$ is a geometric $G$--complex).  Since any two
realizations of $\underline{EG}$ as a geometric $G$--complex are
$G$--equivariantly homotopy equivalent, the $\ell^2$--Betti number
$b^{(2)}_i(\underline{EG};G)$ is an invariant of the group.  We denote
this number by $b^{(2)}_i(G)$.

\subsubsection{Poincar\'e duality}\label{3.3.8}  Suppose $(X, \partial X)$
is a pair of geometric $G$--com\-plexes and also an $n$--dimensional
polyhedral homology manifold with boundary.  Then, by~\ref{2.7}, 
$$b^{(2)}_i(X;G)=b^{(2)}_{n-i}(X, \partial X;G) .$$

\section{\label{4}Simplicial complexes and flag complexes}

\subsection{Definitions and notation}\label{4.1}
Given a simplicial complex $L$, denote by ${\mathcal S}(L)$ the set of
simplices in $L$ together with the empty set $\emptyset$.  It is
partially ordered by inclusion.  ${\mathcal S}_i(L)$ denotes the
subset of ${\mathcal S}(L)$ consisting of the simplices of dimension
$i$.  (For notational purposes it will be convenient to regard
$\emptyset$ as an element of dimension $-1$ in ${\mathcal S}(L)$.)
${\mathcal S}_0(L)$ is the {\it vertex set} of $L$.

\subsubsection{Full subcomplexes}\label{4.1.1}  A subcomplex $A$ of $L$
is a {\it full subcomplex} if whenever $\sigma\in{\mathcal S}(L)$ is
such that the vertex set of $\sigma$ is contained in ${\mathcal
S}(A)$, then $\sigma\in {\mathcal S}(A)$.

\subsubsection{Joins}\label{4.1.2}  Suppose $L_1$ and $L_2$ are simplicial
complexes.  Define a partial order on ${\mathcal
S}(L_1)\times{\mathcal S}(L_2)$ by $(\sigma , \tau )\le (\sigma ', \tau
')$ if and only if $\sigma\le \sigma '$ and $\tau\le \tau '$.  For
example, if $\sigma$ and $\tau$ are simplices of dimension $i$ and
$j$, respectively, then ${\mathcal S}(\sigma )\times{\mathcal S} (\tau
)$ is isomorphic to the poset of faces of a simplex of dimension
$i+j+1$.  We denote this simplex by $\sigma\ast\tau$.  It follows that
there is a unique simplicial complex $L_1\ast L_2$, called the {\it
join} of $L_1$ and $L_2$, characterized by the property that
${\mathcal S}(L_1\ast L_2)$ is isomorphic to ${\mathcal
S}(L_1)\times{\mathcal S}(L_2)$.  The empty element of ${\mathcal
S}(L_1\ast L_2)$ corresponds to $(\emptyset , \emptyset )\in{\mathcal
S}(L_1)\times {\mathcal S} (L_2)$ and the vertex set of $L_1\ast L_2$
corresponds to $({\mathcal S}_0(L_1) \times\{\emptyset\})\cup
(\{\emptyset\}\times{\mathcal S}_0(L_2))$.

As is well known, the geometric realization of $L_1\ast L_2$ is
homeomorphic to the space formed from $L_1\times L_2\times [-1, 1]$ by
identifying points of the form $(x_1, x_2, -1)$ with $(x'_1, x_2, -1)$ and
those of the form $(x_1, x_2, +1)$ with $(x_1, x'_2, +1)$.

\subsubsection{Cones}\label{4.1.3}  The {\it cone} on a simplicial complex
$L$ is the join of $L$ with a single point, say $v$.  We will denote
it by $CL$ (or by $C_vL$ when we wish to distinguish the {\it cone
point} $v$).

\subsubsection{Suspensions}\label{4.1.4}  The {\it suspension} of $L$, 
denoted by $SL$, is the join of $L$ with a $0$--sphere $\BS^0$.

\subsubsection{Incidence relations and flags}\label{4.1.5}  A symmetric
and reflexive relation is an {\it incidence relation}.  Suppose $Q$ is
a set equipped with an incidence relation.  A {\it flag} in $Q$ is a
nonempty finite subset of pairwise related elements.  There is an
associated simplicial complex, $\flag (Q)$, the $i$--simplices of which
are flags of cardinality $i+1$.  (The vertex set of $\flag (Q)$ is $Q$
and two vertices are connected by an edge if and only if they are
incident.)

An important special case is where the incidence relation is given by
symmetrizing the partial order on a poset $P$.  A flag in $P$ is then
a nonempty finite totally ordered subset.  In this case, $\flag (P)$
is called the {\it derived complex} of $P$.  When $P$ is the poset of
cells of a regular $CW$ complex $X$, then $\flag (P)$ can be
identified with the barycentric subdivision of $X$.  As another
example, if $L$ is a simplicial complex, then $\flag ({\mathcal
S}(L))$ is the cone on the barycentric subdivision of $L$.  (The
vertex corresponding to $\emptyset$ is the cone point.)

\subsubsection{}\label{4.1.6}  
Given a poset $P$ and an element $x\in P$, define a
subposet by $P_{\le x} = \{ y\in P\mid y\le x\}$.  Subposets $P_{\ge
x}$, $P_{<x}$ and $P_{>x}$ are defined similarly.

\subsection{Links}\label{4.2}
If $\tau$ is a simplex of $L$, then $\Link(\tau,L)$, 
the {\it link} of $\tau$ in $L$, 
is the union of all simplices $\sigma$ such that
\begin{description}
\item[\rm(a)] intersection of $\sigma$ and $\tau$ is empty and
\item[\rm(b)] $\sigma$ and $\tau$ span a simplex of $L$.
\end{description}

  The subcomplex $\Link(\tau,L)$ is characterized by the condition that
$${\mathcal S}(L_\tau)\cong {\mathcal S}(L)_{\ge \tau}.$$

The {\it star} of $\tau$ in $L$, denoted $\st (\tau, L)$, is the union of all
simplices which intersect $\tau$. 

If $v$ is a vertex of $L$, then we will denote its link $\Link(v,L)$ by $L_v$. 
We have $\st (v, L)=C_vL_v$.  The {\it open
star} of $v$ is the complement of $L_v$ in $\st (v, L)$.  It is an open
subset of $L$.

\subsection{Generalized homology spheres and disks}\label{4.3}
A space $X$ is a {\it homology
$n$--manifold} over a ring $R$  if it has the same local homology groups, 
with coefficients in $R$, as does an $n$--manifold, ie, for all $x\in X$,
$$H_i(X, X-x; R) = \begin{cases}  0 & \text{ if $i\not= n$, }\\ 
R & \text{ if $i=n$.} \end{cases}$$ 

The definition of when a pair $(X,\partial X)$ is a {\it homology 
$n$--manifold with
boundary} over $R$ is similar.
It is well-known (cf \cite{Bre}) that a homology $n$--manifold 
over $R$ satisfies Poincar\'{e} duality over $R$. 
(In non-orientable case one have to use twisted coefficients. 
Also in general, for a finite 
dimensional locally compact space, possibly with a pathological topology, it 
is necessary to use Steenrod homology and \v{C}hech cohomology in order for 
this to be true.)

For the remainder of this paper it can be always assumed that the coefficients 
$R=\BQ$, the field of rational numbers.

\subsubsection{}\label{4.3.1} A simplicial complex $X$ is a  homology
$n$--manifold if and only if it is $n$--dimensional and for each $k$--simplex 
$\sigma$ in $X$, its link $\Link(\sigma,X)$ in $X$ has the same homology as 
$\BS^{n-k-1}$.  

\subsubsection{}\label{4.3.2} 
A simplicial complex $S$ is a {\it generalized homology $n$--sphere}
(abbreviated a $GHS^n$ or simply a $GHS$)  if it is a homology
$n$--manifold with the same homology as $\BS^n$.  A simplicial pair
$(D, \partial D)$ is a {\it generalized homology $n$--disk} (abbreviated
$GHD^n$) if it is a homology $n$--manifold with boundary and if
$$H_i(D, \partial D) = \begin{cases}  0 & \text{ if $i\not= n$, }\\ \BZ
& \text{ if $i=n$.} \end{cases}$$ 

\subsubsection{}\label{4.3.3}  
It follows from \ref{4.3.1} that an $n$--dimensional simplicial complex $X$ 
is a  homology $n$--manifold if and only if for each vertex $v$ of $X$, its 
link $X_v$ is a $GHS^{n-1}$. Similarly, $(X, \partial X)$ is a homology 
$n$--manifold with 
boundary if and only if for each vertex $v$ in $X- \partial X$, its link 
$X_v$ is a $GHS^{n-1}$ and for each $v \in \partial X $, the pair 
$(X_v, X_v \cap  \partial X)$ is a $GHD^{n-1}$.

\subsubsection{}\label{4.3.4} In particular, if $S$ is a $GHS^n$ and $v$ is 
a vertex of $S$, then its link $S_v$ is a $GHS^{n-1}$.
 
\subsubsection{}\label{4.3.5} 
 If $(D, \partial D)$ is a $GHD^n$, then it follows from Poincar\'{e} duality 
and the exact sequence of the pair that $D$ is  acyclic and that $\partial D$ 
has the same homology as does $\BS^{n-1}$. 

\subsubsection{}\label{4.3.6} 
We see from \ref{4.3.3} and \ref{4.3.5} that if a simplicial pair 
$(X, \partial X)$ is a  homology $n$--manifold with 
boundary, then $\partial X$ is a  homology $(n-1)$--manifold.

\subsubsection{}\label{4.3.7} 
If, for $i=1,2$, $S_i$ is a $GHS^{n_i}$, then it follows from the 
K\"{u}nneth Theorem and induction on dimension that the join 
$S_1 \ast S_2$ is a $GHS^{n_1+n_2+1}$. Similarly, if $S$ is a $GHS^n$ and 
$(D, \partial D)$ is a $GHD^m$, then $(S \ast D, S \ast \partial D)$  is 
a $GHD^{n+m+1}$.

\subsubsection{}\label{4.3.8} 
In particular, the suspension of a $GHS^n$ is a $GHS^{n+1}$ and 
the suspension of a $GHD^n$ is a $GHD^{n+1}$.  
 
\subsection{Flag complexes}\label{4.4}
Recall from the Introduction that a simplicial complex $L$ is a {\it
flag complex} if any nonempty finite set of vertices which are
pairwise connected by edges span a simplex in $L$.  In other words, 
$L$ is a flag complex if and only if whenever a subcomplex isomorphic
to the $1$--skeleton of a simplex is in $L$, then the entire simplex
lies in $L$.  (In \cite{G1} Gromov used the terminology that $L$
satisfies the ``no $\Delta$ condition'' for this property.)

\subsubsection{}\label{4.4.1}  
If $Q$ is a set with an incidence relation, then $\flag(Q)$ 
(defined in~\ref{4.1.5}) is a flag complex.  Conversely, any flag
complex arises from this construction.  (Indeed, given a flag complex
$L$, define two vertices in ${\mathcal S}_0(L)$ to be {\it incident}
if they are connected by an edge.  Then $L\cong\flag ({\mathcal
S}_0(L))$.)

\subsubsection{}\label{4.4.2}  
In particular, the barycentric subdivision of any
regular $CW$ complex is a flag complex.  Hence, the condition that $L$
be a flag complex imposes no restriction on its topological type:  it
can be any polyhedron.

\subsubsection{}\label{4.4.3}  
An $m$--gon (ie, a triangulation of a circle into $m$
edges) is a flag complex if and only if $m\ge 4$.

\subsubsection{}\label{4.4.4}  
Any full subcomplex of a flag complex is a flag complex.

\subsubsection{}\label{4.4.5}  
If $v$ is a vertex of a flag complex $L$, then its
link $L_v$ and its star $\st (v, L)$ are both full subcomplexes.
Hence, by 4.4.4, they are both flag complexes.

\subsubsection{Joins of flag complexes}\label{4.4.6}  The join of two flag
complexes is again a flag complex.  In particular, the cone on a flag
complex is a flag complex and the suspension of a flag complex is a
flag complex.

\subsubsection{Notation}\label{4.4.7}  For any set of vertices $T$ of $L$, 
let $N(T)$ be the union of all open stars of vertices in $T$.  We will
use $L-T$ to denote the complement of $N(T)$ in $L$.  In other words, 
$L-T$ is the full subcomplex of $L$ spanned by ${\mathcal S}_0(L)-T$.
For example, for any vertex $s$ of $L$, $L-s$ denotes the complement
of the open star of $s$ in $L$.  Similarly, if $A$ is any subcomplex
of $L$, then we will write $L-A$ for $L-{\mathcal S}_0(A)$.

\section{\label{5}Right-angled Coxeter groups}

\subsection{Definition of $\boldsymbol W_{\boldsymbol L}$}\label{5.1}
Suppose $L$ is a flag complex.  The $1$--skeleton of $L$ gives the data
for the presentation of a group $W_L$.  The set of generators in the
presentation is the vertex set ${\mathcal S}_0(L)$.  The edges of $L$
give relations, as follows:
\begin{align*}
s^2 & =  1 , \qquad\text{for all $s\in{\mathcal S}_0(L)$, }\\
(st)^2 & = 1 , \qquad\text{whenever $~\{ s, t\}$ spans an edge
in $L$.}
\end{align*}
The group $W_L$ is the {\it right-angled Coxeter group} associated to
$L$.  ${\mathcal S}_0(L)$, regarded as a subset of $W_L$, is the {\it
fundamental set of generators}.  The flag complex $L$ is called the {\it
nerve} of $(W_L, {\mathcal S}_0(L))$.

\subsection{Examples}\label{5.2}
We give some examples of this construction for various flag complexes
$L$.

\subsubsection{The empty set}\label{5.2.1}  If $L=\emptyset$, then
$W_\emptyset$ is the trivial group.

\subsubsection{A $0$--simplex}\label{5.2.2}  If $L$ is a single point $s$, 
then $W_s\cong\BZ_2$, the cyclic group of order $2$.

\subsubsection{Joins}\label{5.2.3}  By~\ref{4.1.2}, 
$W_{L_1\ast L_2} = W_{L_1}\times W_{L_2}$.

\subsubsection{Cones}\label{5.2.4}  By~\ref{5.2.2} and \ref{5.2.3}, 
$W_{CL} = \BZ_2\times W_L$.

\subsubsection{A $k$--simplex}\label{5.2.5}  If $\sigma$ is a $k$--simplex with
vertex set $\{ s_0, \dots , s_k\}$, then by~\ref{5.2.2} and~\ref{5.2.3}, 
$W_\sigma = W_{s_0}\times\cdots\times W_{s_k}\cong (\BZ_2)^{k+1}$.

\subsubsection{Disjoint unions}\label{5.2.6}  If $L$ is the disjoint union
of two flag complexes $L_1$ and $L_2$, then $W_L$ is the free product
of $W_{L_1}$ and $W_{L_2}$, ie, $W_{L_1\cup L_2} = W_{L_1}\ast
W_{L_2}$.

\subsubsection{Amalgamated products}\label{5.2.7}  More generally, if $L=
L_1\cup L_2$, $L_1\cap L_2 =A$, where $L_1$ and $L_2$ (and therefore, 
$A$) are full subcomplexes, then $W_L$ is the amalgamated product:
$$W_L=W_{L_1}\ast_{W_A} W_{L_2} .$$

\subsubsection{$k$ points}\label{5.2.8}  If $L$ is the disjoint union of $k$
points $s_1, \dots , s_k$, then $W_L$ is the free product
$W_{s_1}\ast\cdots\ast W_{s_k}$ ($\cong\BZ_2\ast\cdots\ast\BZ_2$).  In
particular, $W_{\BS^0}$ is the infinite dihedral group $D_\infty$.

\subsubsection{Suspensions}\label{5.2.9}  By~\ref{5.2.3} and~\ref{5.2.8}, 
$W_{SL} = D_\infty\times W_L$.

\subsection{Special subgroups}\label{5.3}
Let $A$ be a full subcomplex of $L$.  By \cite[Th\'e\-or\`eme 2, 
p. 20]{B}, $W_A$ can be identified with the subgroup of $W_L$ generated
by ${\mathcal S}_0 (A)$.  (N. B. Here it is important that $A$ be a
full subcomplex; for if two vertices of $A$ were connected by an edge
in $L$ which was not in $A$, then there would be a relation in $W_L$
not satisfied in $W_A$.)  Such a subgroup $W_A$ is called a {\it
special subgroup} of $W_L$.

\subsubsection{}\label{5.3.1}We note that a special subgroup $W_A$ is  finite
if and only if $A$ is a simplex of $L$ (or if $A=\emptyset$).  The
special subgroups of $W_L$ corresponding to the elements of ${\mathcal
S}(L)$ are sometimes called the {\it spherical} special subgroups.

\subsection{The poset of spherical cosets}\label{5.4}
A {\it spherical coset} in $W_L$ is a coset of the form $wW_\sigma$
for some $\sigma\in{\mathcal S}(L)$ and $w\in W_L$.  The set of all
spherical cosets will be denoted by $W_L{\mathcal S}(L)$, ie, 
$$W_L{\mathcal S}(L) = \bigcup_{\sigma\in{\mathcal S}(L)} W_L/W_\sigma.$$ 
It is partially ordered by inclusion of one coset in another.  The
group $W_L$ acts in an obvious way on the poset $W_L{\mathcal S}(L)$.
The quotient poset is ${\mathcal S}(L)$.

\section{\label{6}The complex $\Sigma_L$}

We retain the notation of the previous section:  $L$ is a finite flag
complex, $W_L$ is the associated right-angled Coxeter group and $W_L
{\mathcal S}(L)$ is the poset of spherical cosets.

\subsection{Definitions and basic properties}\label{6.1}
The space $\Sigma_L$ is defined as the geometric realization of the
poset $W_L{\mathcal S}(L)$.  (In other words, it is the simplicial
complex $\flag (W_L{\mathcal S}(L))$.)  Let $K_L$ denote the geometric
realization of ${\mathcal S} (L)$.  (By~\ref{4.1.5}, $K_L$ is the cone
on the barycentric subdivision $L$.)  The inclusion ${\mathcal
S}(L)\hookrightarrow W_L{\mathcal S}(L)$, defined by  $\sigma\mapsto
W_\sigma$, induces an inclusion $K_L\subset\Sigma_L$.  When regarded
in this way as a subset of $\Sigma_L$, $K_L$ is called the {\it
fundamental chamber}.

\subsubsection{The $W_L$--action}\label{6.1.1}  The natural $W_L$--action
on $W_L{\mathcal S}(L)$ induces a simplicial action on $\Sigma_L$.
The orbit space is $K_L$.  The action is proper (since each cell
stabilizer is a conjugate of a spherical special subgroup) and
cocompact (since ${\mathcal S}(L)$ is finite).

\subsubsection{Contractibility}\label{6.1.2}  It is proved in \cite{D1} that
$\Sigma_L$ is contractible.  In fact, $\Sigma_L$ is the universal
space for proper $W_L$--actions, in the sense of~\ref{1.4}.

\subsubsection{Special subcomplexes}\label{6.1.4}  Suppose $A$ is a full
subcomplex of $L$.  The inclusion $W_A\to W_L$ induces an inclusion of
posets $W_A{\mathcal S}(A)\to W_L{\mathcal S}(L)$ and hence, an
inclusion of $\Sigma_A$ as a subcomplex of $\Sigma_L$.  Such a
$\Sigma_A$ will be called a {\it special subcomplex} of $\Sigma_L$.
If $w\in W_L-W_A$, then $\Sigma_A$ and $w\Sigma_A$ are disjoint
subcomplexes.  It follows that the stabilizer of $\Sigma_A$ in
$\Sigma_L$ is $W_A$ and that
\subsubsection{}\label{6.1.5} 
\myequation{W_L\Sigma_A\cong W_L\times_{W_A}\Sigma_A,}

where $W_L\Sigma_A$ denotes the union of all translates of $\Sigma_A$
in $\Sigma_L$.

\subsection{Examples}\label{6.2}
We consider the above construction for the same flag complexes $L$ as
in~\ref{5.2}.

\subsubsection{The empty set}\label{6.2.1}  $\Sigma_\emptyset$ is a point.

\subsubsection{A $0$--simplex}\label{6.2.2}  If $L$ is a single point $s$, 
then $\Sigma_s$ can be identified with the interval $[-1, 1]$.  The
nontrivial element $s\in W_s$ ($W_s\cong\BZ_2$) acts as the reflection
$t\to -t$.

\subsubsection{Joins}\label{6.2.3}  By~\ref{4.1.2}, 
${\mathcal S}(L_1\ast L_2) \cong {\mathcal S}(L_1)\times{\mathcal
S}(L_2)$ and by~\ref{5.2.3}, $W_{L_1\ast L_2}=W_{L_1} \times W_{L_2}$.
It follows that
$$\Sigma_{L_1\ast L_2} = \Sigma_{L_1} \times \Sigma_{L_2}$$ with the
product action.

\subsubsection{Cones}\label{6.2.4}  
\myequation{\Sigma_{CL}=[-1, 1]\times\Sigma_L.}

\subsubsection{A $k$--simplex}\label{6.2.5}  If $\sigma$ is a $k$--simplex
with vertex set $\{ s_0, \dots , s_k\}$, then by~\ref{5.2.5}, $W_\sigma
=W_{s_0}\times \cdots\times W_{s_k}$ ($\cong(\BZ_s)^{k+1}$) and
by~\ref{6.2.2} and~\ref{6.2.3}, 
$$\Sigma_\sigma = \Sigma_{s_0}\times\cdots\times\Sigma_{s_k}~~(\cong
[-1, 1 ]^{k+1}) .$$

\subsubsection{Disjoint unions}\label{6.2.6}  If $L$ is the disjoint union
of $L_1$ and $L_2$, then $K_L$ is the one point union $K_{L_1}\vee
K_{L_2}$ (the common point corresponding to $\emptyset\in{\mathcal
S}(L_1)\cap{\mathcal S}(L_2)$).

\subsubsection{$k$ points}\label{6.2.7}  Suppose $L=P_k$, the disjoint union
of $k$ points.  Then $K_{P_k}$ is the cone on $k$ points  and, if
$k>1$, $\Sigma_{P_k}$ is the regular infinite tree where each vertex
has valence $k$.

\subsubsection{The $0$--sphere}\label{6.2.8}  In particular, $\Sigma_{\BS^0}$
can be identified with the real line $\BR$ cellulated as the union of
intervals of the form $[2m-1, 2m+1]$, $m\in\BZ$.  The action of the
infinite dihedral group $W_{\BS^0}$ is the standard one, generated by
the reflections across $0$ and $2$.

\subsubsection{Suspensions}\label{6.2.9}  By~\ref{6.2.3} and~\ref{6.2.8}, 
$\Sigma_{SL} = \BR\times\Sigma_L$.

\subsection{The cubical structure on 
$\boldsymbol \Sigma_{\boldsymbol L}$}\label{6.3}

\subsubsection{The case where $L$ is a simplex}\label{6.3.1}  Suppose
$\sigma$ is a $k$--simplex.  Then by~\ref{5.2.5}, $W_\sigma\cong
(\BZ_2)^{k+1}$ and by 6.2.5, $\Sigma_\sigma = [-1, 1]^{k+1}$.  The
group $W_\sigma$ acts simply transitively on the set of
$0$--dimensional faces ($=$ ``vertices'') of $[-1, 1]^{k+1}$.  Moreover, 
a set of such vertices is the vertex set of a face of $[-1, 1]^{k+1}$
if and only if it corresponds to the set of elements in a coset of the
form $w W_\tau$, for some $w\in W_\sigma$ and $\tau\in{\mathcal
S}(\sigma )$.  Hence, the poset of nonempty faces of $[-1, 1]^{k+1}$
($=\Sigma_\sigma$) is naturally identified with the poset
$W_\sigma{\mathcal S}(\sigma )$.

\subsubsection{The general case}\label{6.3.2}  Now suppose that $L$ is an
arbitrary flag complex.  For each $\sigma\in{\mathcal S}_k(L)$ and
$w\in W_L$, the subcomplex $w\Sigma_\sigma$ of $\Sigma_L$ is
homeomorphic to $[-1, 1]^{k+1}$.  This gives a decomposition of
$\Sigma_L$ into a family of subcomplexes, $\{w\Sigma_\sigma\}_{w
W_\sigma\in W_L{\mathcal S}(L)}$.  The family is indexed by the poset
of spherical cosets $W_L{\mathcal S}(L)$.  Each subcomplex is
homeomorphic to a cube.  Thus,  $\Sigma_L$ has the structure of a
regular $CW$ complex in which (a) the poset of cells is identified
with $W_L{\mathcal S}(L)$ and (b) the cell corresponding to $w
W_\sigma$ is a $(k+1)$--dimensional cube, where $k=\dim\sigma$.  As
before, there is a $0$--dimensional cube (vertex) for each element of
$W_L$ ($=W_L/W_\emptyset$) and a set of such $0$--cubes is the vertex
set of a $(k+1)$--cube, $w \Sigma_\sigma$, if and only if it is the set
of elements in the spherical coset $w W_\sigma$.

\subsubsection{The link of a vertex in $\Sigma_L$}\label{6.3.3}  With respect
to this cubical structure, the link of each vertex of $\Sigma_L$ is
$L$.  In other words, the poset of cubes of $\Sigma_L$ which properly
contain a given vertex is canonically identified with ${\mathcal
S}(L)_{>\emptyset}$.

\subsubsection{The orbihedral Euler characteristic 
of $\Sigma_L/W_L$}\label{6.3.4}   The $W_L$--orbits of cubical cells in
$\Sigma_L$ are bijective with ${\mathcal S}(L)$.  The dimension of a
cube in an orbit corresponding to $\sigma \in{\mathcal S}_k(L)$ is
$k+1$ and the order of its stabilizer is $2^{k+1}$.  Hence, 
by~\ref{1.3.1}, the orbihedral Euler characteristic is given by
$$\chi^{\orb}(\Sigma_L/W_L) = \sum_{\sigma\in{\mathcal
S}(L)}\left(-\frac{1}{2} \right)^{\dim\sigma +1}$$ or
$$\chi^{\orb}(\Sigma_L/W_L) = \sum^{\dim
L}_{k=-1}\left(-\frac{1}{2}\right)^{k +1}f_k(L) , $$ where $f_k(L)$ is
the number of elements in ${\mathcal S}_k(L)$.  We note that the right
hand side of the last equation is precisely the quantity $\kappa (L)$, 
mentioned in the Introduction, in connection with the Combinatorial
Gauss--Bonnet Theorem.  It is the local contribution to the Euler
characteristic coming from the link of a vertex in a piecewise
Euclidean, cubical cell complex.  (See \cite{CD}.)

Since $\Sigma_L$ is the universal space for proper $W_L$--actions, 
$\chi^{\orb} (\Sigma_L/W_L)$ is the Euler characteristic of $W_L$.

\subsubsection{}\label{6.3.5} Each special subcomplex of $\Sigma_L$ is also a
subcomplex in the cubical structure.

\subsubsection{}\label{6.3.6}  
For any $s\in{\mathcal S}_0(L)$, $\Sigma_s$ is an edge of
$\Sigma_L$.  Let $O(s)$ denote the union of the interiors of all cubes
of $\Sigma_L$ which have $\Sigma_s$ as a face (ie, $O(s)$ is the
open star of the interior of $\Sigma_s$).  For any subset $T$ of
${\mathcal S}_0(L)$, set
$$R(T)=\bigcup_{s\in T} W_L O(s) .$$ Thus, $R(T)$ is an open, 
$W_L$--stable subset of $\Sigma_L$.  Moreover, with notation as
in~\ref{4.4.7} and~\ref{6.1.5}, we have that
\subsubsection{}\label{6.3.7}
$\Sigma_L-R(T) = W_L\Sigma_{L-T} .$

\subsection{The commutator cover of 
$\boldsymbol \Sigma_{\boldsymbol L}/\boldsymbol W_{\boldsymbol L}$}\label{6.4}
In this section we will describe a finite cubical complex $P_L$ as a
subcomplex of a Euclidean cube.  It turns out that the universal cover
of $P_L$ can be identified with $\Sigma_L$.  This gives an
alternative, and perhaps more easily understandable method of
describing the cubical structure on $\Sigma_L$.

\subsubsection{The commutator subgroup}\label{6.4.1}  The abelianization
of $W_L$, denoted $W^{\ab}_L$, is obviously $(\BZ_2)^{{\mathcal
S}_0(L)}$, the direct product of cyclic groups of order two.  Let
$\varphi\colon W_L \to W^{\ab}_L$ be the natural epimorphism.  Its
kernel, denoted by $\Gamma_L$, is the commutator subgroup. 
Since any finite subgroup of $W_L$ is
contained in a conjugate of a finite special subgroup and since the
restriction of $\phi$ to any finite special subgroup is injective, $\Gamma
_L$  is torsion-free.  
 Hence, it acts freely on $\Sigma_L$.
The natural projection $\Sigma_L/\Gamma_L\to \Sigma_L/W_L$ is an
orbihedral covering in the sense of~\ref{1.2}; we call $\Sigma_L
/\Gamma_L$ the {\it commutator cover} of $\Sigma_L/W_L$.

\subsubsection{The complex $P_L$}\label{6.4.2}  Let $\square$ denote the
Euclidean cube $[-1, 1]^{{\mathcal S}_0(L)}$.  For each
$\sigma\in{\mathcal S}(L)$, let $\square_\sigma$ be the face of
$\square$ defined by
$$\square_\sigma = e\times [-1, 1]^{{\mathcal S}_0(\sigma)}$$ 
where $e$ is
the vertex of $[-1, 1]^{{\mathcal S}_0(L)-{\mathcal S}_0(\sigma )}$
with all coordinates equal to $1$.  The faces of $\square$ which are
parallel to $\square_\sigma$ have the form $f\times [-1, 1]^{{\mathcal
S}_0(\sigma )}$, where $f$ is some vertex of $[-1, 1]^{{\mathcal
S}_0(L)-{\mathcal S}_0(\sigma )}$.

Define $P_L$ to be the union of all faces of $\square$ which are
parallel to $\square_\sigma$, for some $\sigma\in{\mathcal S}(L)$.
Thus, $P_L$ is a subcomplex of $\square$.

Each generator $s$ of $W^{\ab}_L$ acts on $\square$ as reflection in
the $s^{\text{th}}$ coordinate.  Thus, $W^{\ab}_L$ acts on $\square$
as a finite reflection group.  The orbit space is $[0, 1]^{{\mathcal
S}_0(L)}$ and $P_L /W^{\ab}_L$ is the subcomplex consisting of all
faces of the form $e\times [0, 1]^{{\mathcal S}_0(\sigma )}$, for some
$\sigma\in{\mathcal S}(L)$.  (Moreover, this subcomplex can be
canonically identified with $K_L=\Sigma_L/W_L$.)

\subsubsection{Identification of $P_L$ with
$\Sigma_L/\Gamma_L$}\label{6.4.3} There is a natural
$\varphi$--equivariant map $p\colon\Sigma_L\to P_L$ which sends the
cube $w\Sigma_\sigma$ to $\varphi (w )\square_\sigma$.  It is obvious
that $p$ is a covering projection and that it induces an isomorphism
from $\Sigma_L/\Gamma_L$ onto $P_L$.  Henceforth, we identify these
two cubical complexes.

\subsection{The piecewise Euclidean metric on 
$\boldsymbol \Sigma_{\boldsymbol L}$}\label{6.5}
We review some material from \cite{G1} (which can also be found in
\cite{BH}, \cite{D2}, \cite{D3}, or \cite{M}).

Identify each $k$--dimensional cube in $\Sigma_L$ with the regular
Euclidean cube
of edge length $2$.  The length of a
piecewise linear curve in $\Sigma_L$ is then unambiguously defined.
The {\it distance} $d(x, y)$ between two points $x$ and $y$ in
$\Sigma_L$ is then defined to be the infimum of the lengths of
piecewise linear paths connecting them.  With this metric, $\Sigma_L$
becomes a {\it geodesic space}, that is, for any two points $x$ and
$y$ there is a path of length $d(x, y)$ between them.  Such a path is
called a {\it geodesic segment}.

\subsubsection{Nonpositive curvature}\label{6.5.1}  For a geodesic space
$X$ the concept of ``nonpositive curvature'' can be defined by
comparing distances on small triangles in $X$ (ie, configurations of
three geodesic segments in $X$) with distances on comparison triangles
in the Euclidean plane.  $X$ is {\it nonpositively curved} if Gromov's
$\text{CAT} (0)$--inequality (page 106 of \cite{G1}) holds for all
sufficiently small triangles in $X$.

\begin{Lemma}[Gromov]\label{6.5.2}  A cubical cell complex $X$ with
piecewise Euclidean metric defined as above is nonpositively curved if
and only if the link of each vertex is a flag complex.
\end{Lemma}
  This is
proved on page 123 of \cite{G1}.  The proof can also be found in
\cite{BH}, \cite{D3}, or~\cite{M}.

\subsubsection{}\label{6.5.3}
It follows from~\ref{6.3.3} and 
Gromov's Lemma that $\Sigma_L$
is nonpositively curved.  Since, by~\ref{6.1.2}, $\Sigma_L$ is 
contractible, this implies that it is a $\text{CAT}(0)$--space (ie, that the
$\text{CAT}(0)$ inequality holds for all triangles).

\subsubsection{}\label{6.5.4}It is not difficult to show that any special
subcomplex $\Sigma_A$ is a geodesically convex subspace of $\Sigma_L$. 
See  Proposition 1.7.1, page 514 of \cite{DJS} for details.

\subsection{Reflection groups on manifolds}\label{6.6}

\subsubsection{Classical reflection groups}\label{6.6.1}  Let $\BX^n$
stand for either Euclidean $n$--space $\BE^n$, hyperbolic $n$--space
$\BH^n$ or the $n$--sphere $\BS^n$.  A {\it classical reflection group}
$W$ is a discrete, cocompact group of isometries of $\BX^n$
generated by reflections.  Then $W$ is a Coxeter group.  (The theory
of general Coxeter groups arose from the study of this classical
situation.)

Suppose $W$ is a classical reflection group on $\BX^n$.  Choose
a component of the complement of the union of the reflecting
hyperplanes and call its closure $K$.  Then $K$ is a convex polytope.
Moreover, it is a fundamental domain for the $W$--action and the set of
reflections across the codimension-one faces of $K$ is a fundamental
set of generators for $W$.

Let $S$ be the simplicial complex dual to the boundary of $K$.  In the
spherical case, $S$ is the boundary of a simplex (and hence, not a
flag complex when the dimension of the simplex is greater than $1$).
In the Euclidean case, $S$ is the join of boundaries of simplices.

The condition that $W$ be right-angled means that the codimension-one
faces of $K$ are orthogonal whenever they intersect.  In the
right-angled Euclidean case, $\BX^n=\BE^n$, the only possibility is that
$K$ is a product of intervals, $S$ is the boundary of an
$n$--dimensional octahedron (an $n$--fold join of $0$--spheres) and $W=W_S$ is
an $n$--fold product of infinite dihedral groups.  If $K$ is the
regular $n$--cube $[0, 2]^n$ (which we may assume after conjugating by
an affine automorphism) then $\Sigma_S$ is isometric with $\BE^n$.
In the right-angled hyperbolic case, $\Sigma_S$ is equivariantly
homeomorphic to $\BH^n$ but not isometric to it.  (They are 
quasi-isometric.)  The cubical structure on $\Sigma_S$ is dual to the
tessellation of $\BH^n$ by the translates of $K$.

\subsubsection{An $m$--gon}\label{6.6.2}  Suppose $S$ is an $m$--gon, ie, 
a subdivision of the circle into $m$ edges.  To insure that $S$ is a
flag complex, we also assume $m\ge 4$.  Then $W_S$ is isomorphic to a
classical reflection group and $\Sigma_S$ is combinatorially dual to a
tessellation of the Euclidean plane (when $m=4$) by squares or to a
tessellation of the hyperbolic plane (when $m>4$) by right-angled
$m$--gons.

\subsubsection{Spheres}\label{6.6.3}  Suppose that $S$ is a triangulation
of $\BS^{n-1}$ as a flag complex.  Then, by~\ref{6.3.3}, $\Sigma_S$ is
an $n$--dimensional manifold (since a neighborhood of each vertex is
homeomorphic to the cone on $S$).  If $n>3$, then very few of these
triangulations correspond to classical reflection groups.  The
situation in dimension $3$ will be explained in Section~\ref{10}.

\subsubsection{Generalized homology spheres}\label{6.6.4}  Similarly, if
$S$ is a $GHS^{n-1}$, as defined in~\ref{4.3}, the, by~\ref{6.3.3}, 
$\Sigma_S$ is a polyhedral homology $n$--manifold.

\subsubsection{Generalized homology disks}\label{6.6.5}  If $D$ is a
triangulation of an $(n-1)$--disk as a flag complex and $\partial D$ is
a full subcomplex, then, by~\ref{6.3.3}, $\Sigma_D$ is an $n$--manifold
with boundary.  Its boundary is $W_D\Sigma_{\partial D}$.  Similarly, 
if $(D, \partial D)$ is a $GHD^{n-1}$, as defined in~\ref{4.3}, then
$\Sigma_D$ is a polyhedral homology $n$--manifold with boundary.

\section{\label{7}Properties of the $\ell^2$--homology of $\Sigma_L$}

From now on, {\it all simplicial complexes will be flag complexes and
all subcomplexes will be full subcomplexes}.  Given a finite flag
complex $L$, we have associated a group $W_L$, a geometric
$W_L$--complex $\Sigma_L$ and then, for each $i\in\BN$, a Hilbert
$W_L$--module, ${\h}_i(\Sigma_L)$.  Similarly, to each pair
$(L, A)$ we can associate the Hilbert $W_L$--module, ${\h}_i
(\Sigma_L, W_L\Sigma_A)$ (where, by~\ref{6.1.5}, $W_L\Sigma_A\cong
W_L\times_{W_A}\Sigma_A)$.

We introduce some useful notation which reflects this situation.

\subsection{Notation}\label{7.1}

\subsubsection{}\label{7.1.1} 
\myequation{ \H_i(L) = {\h}_i(\Sigma_L)}

\subsubsection{}\label{7.1.2} 
\myequation{\H_i(A) = {\h}_i(W_L\Sigma_A)}

\subsubsection{}
\label{7.1.3} 
\myequation{\H_i(L, A) = {\h}_i(\Sigma_L, W_L\Sigma_A)}

\subsubsection{}
\label{7.1.4} 
\myequation{\beta_i(A) = \dim_{W_L}(\H_i(A)) }

\subsubsection{}
\label{7.1.5} 
\myequation{\beta_i(L, A) = \dim_{W_L}(\H_i(L, A)) }

\subsubsection{}
\label{7.1.6} 
\myequation{\chi^{(2)}(L) = \sum (-1)^i\beta_i(L)} 

The notation in~\ref{7.1.2} and~\ref{7.1.4} will not lead to
confusion, since, by~\ref{2.4.5} and~\ref{6.1.5}, $\h_i(W_L\Sigma_A)$ 
is the induced representation from $\h_i(\Sigma_A)$ and, therefore, 
by~\ref{3.3.4}, 
$$b^{(2)}_i(W_L\Sigma_A;W_L)=b^{(2)}_i(\Sigma_A;W_A) .$$

\subsection{Basic algebraic topology}\label{7.2}  For the case at hand, we
rewrite some of the basic properties of reduced $\ell^2$--homology in
our new notation.  From~\ref{2.4.2}, we get the following.

\begin{Lemma}[\label{7.2.1}Exact sequence of the pair]   The sequence
$$\to \H_i(A)\to\H_i (L)\to\H_i(L, A)\to$$ is weakly exact.
\end{Lemma}

\begin{Lemma}[Excision]\label{7.2.2}  Given $(L, A)$ as above, let $T$ be a set
of vertices of $A$ such that the open star of any vertex in $T$ is
contained in the interior of $A$.  Then, with notation as
in~\ref{4.4.7}, 
$$\H_i(L, A)\cong \H_i(L-T, A-T) .$$
\end{Lemma}

\begin{proof}  This is immediate from~\ref{2.4.3} and~\ref{6.3.7}.
\end{proof}

\begin{Lemma}[Mayer--Vietoris sequences]\label{7.2.3}  Suppose $L=L_1\cup L_2$
and $A=L_1\cap L_2$, where $L_1$ and  $L_2$ (and therefore, $A$) are
full subcomplexes of $L$.
\begin{enumerate}
\item The Mayer--Vietoris sequence
$$\to\H_i(A)\to\H_i(L_1)\oplus\H_i(L_2)\to\H_i(L)\to$$ is weakly exact.

\item 
$\H_i(L, A)\cong \H_i(L_1, A)\oplus\H_i(L_2, A) .$
\end{enumerate}
\end{Lemma}

\begin{proof}  
Statement (1) follows from~\ref{2.4.4}.  For (2), use the following
relative version of the Mayer--Vietoris sequence, 
$$\to\H_i(A, A)\to\H_i(L_1, A)\oplus\H_i(L_2, A)\to\H_i(L, A)\to
\H_{i-1}(A, A)$$
and the fact that $\H_*(A, A)=0$.
\end{proof}

Using~\ref{5.2.3} and~\ref{6.2.3} the K\"unneth Formula, \ref{3.3.5}, 
translates to the following.

\begin{Lemma}[\label{7.2.4}The Betti numbers of a join]
$$\beta_k(L_1\ast L_2) = \sum_{i+j=k}\beta_i(L_1)\beta_j(L_2) .$$
\end{Lemma}

Using~\ref{6.3.4}, Atiyah's Formula, \ref{3.3.6}, translates as
follows.

\begin{Lemma}[Atiyah's Formula]\label{7.2.5}
$$\chi^{(2)}(L)=\sum^{\dim L}_{k=-1}\left(
-\frac{1}{2}\right)^{k+1}f_k(L) .$$
\end{Lemma}

\subsubsection{$0$--dimensional homology}\label{7.2.6}  If $L$ is nonempty
and not a simplex, then, by~\ref{2.5.1}, 
$$\b_0(L)=0.$$

\subsubsection{}\label{7.2.7}  Similarly, suppose $L$ is a pseudomanifold
of dimension $n-1$, as in~\ref{2.6}.  It then follows from~\ref{6.3.3}
that $\Sigma_L$ is an $n$-dimensional pseudomanifold and it can be
seen that each component of the complement of the codimension $2$
skeleton is infinite. Hence, by~\ref{2.6},
$$\b_n(L)=0.$$

\subsection{Examples}\label{7.3}
Next we calculate the Betti numbers, $\beta_i(L)$, for some of the
examples in~\ref{5.2} and~\ref{6.2}.

\subsubsection{The empty set}\label{7.3.1}  Since $W_\emptyset$ is trivial
and $\Sigma_\emptyset$ is a point, 
$$\beta_i(\emptyset )=\begin{cases} 1 &\text{if $i=0$, } \\ 0 &\text{if
$i\not= 0$.} \end{cases}$$

\subsubsection{A $k$--simplex}\label{7.3.2}  Given a $k$--simplex $\sigma$, 
$W_\sigma\cong (\BZ_2)^{k+1}$ and $\Sigma_\sigma$\break$ = [-1, 1]^{k+1}$.
Hence, 
$$\beta_i(\sigma ) =\begin{cases}  \left(\frac{1}{2}\right)^{k+1}
&\text{if $i=0$,}\\ 0 &\text{if $i\not= 0$.}
\end{cases}$$

\begin{Lemma}[The Betti numbers of a disjoint union]\label{7.3.4}
  Suppose that $L$ is the disjoint union of $L_1$ and $L_2$.  Then, 
for $i\ge 2$, 
$$\beta_i(L)=\beta_i(L_1)+\beta_i(L_2).$$ If neither $L_1$ nor $L_2$
is a simplex, then
$$\beta_1(L)=\beta_1(L_1)+\beta_1(L_2)+1.$$
\end{Lemma}

\begin{proof}  
This follows from the Mayer--Vietoris sequence, Lemma~\ref{7.2.3} (1), 
after noting that $L_1\cap L_2=\emptyset$ has nonzero Betti number, 
$\beta_0 (\emptyset )=1$.  The final sentence follows since if
$W_{L_1}$ and $W_{L_2}$ are both infinite, then, by~\ref{7.2.6}, 
$\beta_0(L_1)=\beta_0(L_2)=0$.
\end{proof}

\begin{Lemma}[The Betti numbers of $k$ points]\label{7.3.5}  Let $P_k$ denote
the disjoint union of $k$ points.  If $k\ge 2$, then
$$\beta_i (P_k) = \begin{cases}   \frac{k}{2} -1  &\text{if $i=1$, }\\
0              &\text{if $i\not= 1$.}
\end{cases}$$
In particular, 
$$\beta_i(\BS^0)=\beta_i(P_2)=0 \quad \text{   for all $i$. }$$
\end{Lemma}

\begin{proof}  Since $\Sigma_{P_k}$ is $1$--dimensional, $\beta_i(P_k)=0$ for
$i>1$.  Since $k\ge 2$, $\beta_0(P_k)=0$, by~\ref{7.2.6}.
By~\ref{6.3.4}, $\chi^{(2)}(P_k) =1-\frac{k}{2}$.  Hence, by Atiyah's
Formula~\ref{7.2.5}, $\beta_1(L)=-\chi^{(2)}(L)= \frac{k}{2}-1$.
\end{proof}

\begin{Lemma}[The Betti numbers of a suspension]\label{7.3.6}  $\beta_i(SL)=0$
for all $i$.
\end{Lemma}

\begin{proof}  This follows from Lemma~\ref{7.2.4} and~\ref{7.3.5}.
\end{proof}

\subsubsection{}\label{7.3.7}  
Suppose $L$ is the $m$--fold join, $L=P_{k_1}\ast\cdots
\ast P_{k_m}$ where each $k_j\ge 2$.  Then, by Lemmas~\ref{7.2.4} 
 and~\ref{7.3.5}, 
$$\beta_i(L)=\begin{cases}  \prod\left(\frac{k_j}{2} -1\right)
&\text{if $i=m$, }\\ 0 &\text{if $i\not= m$.}\end{cases}$$

\begin{Lemma}[The Betti numbers of a cone]\label{7.3.3} \

\begin{enumerate}

\item $\beta_i(CL)=\frac{1}{2}\beta_i(L)$.

\item $\beta_{i+1}(CL, L)=\frac{1}{2}\beta_i(L)$.

\item The sequence of the pair $(CL, L)$ breaks up into short exact
sequences:
$$0\to\H_{i+1}(CL, L)\to\H_i(L)\to\H_i(CL)\to 0 .$$
\end{enumerate}
\end{Lemma}

\begin{proof}  Although formulas (1) and (2) follow from the proof of (3), we
first give simple alternative arguments for them which illustrate the
above methods.  Since $\Sigma_{CL}=[-1, 1]\times\Sigma_L$, the
complexes $\Sigma_{CL}$ and $\Sigma_L$ are $W_L$--equivariantly
homotopy equivalent; hence, ${\h}_i(\Sigma_{CL})\cong{\mathcal
H}_i(\Sigma_L)$.  Since $W_{CL}=\BZ_2 \times W_L$, we have, 
by~\ref{3.2.9}, that $\beta_i(CL)=\frac{1}{2}\beta_i(L)$, proving (1).
To prove (2), let $-1$ and $+1$ denote the two points of $\BS^0$.
Then $SL=C_{-1}L\cup C_{+1}L$ is the union of two copies of the cone
on $L$ along $L$.  By excision, Lemma~\ref{7.2.2}, $\H_{i+1}
(C_{+1}L, L)\cong \H_{i+1}(SL, C_{-1}L)$ and by the exact sequence of
the pair, $\H_{i+1} (SL, C_{-1}L)\cong\H_i(C_{-1}L)$.  Hence, 
$\beta_{i+1}(CL, L)=\beta_{i+1}
(SL, CL)=\beta_i(CL)=\frac{1}{2}\beta_i(L)$, which proves (2).

In the exact sequence which we are considering in (3), $\Sigma_{CL}$
is the ambient space and $\H_i(L)$ means the reduced $\ell^2$--homology
of the subcomplex $\{\pm 1\}\times \Sigma_L$ in $\Sigma_{CL}$
($=[-1, 1]\times \Sigma_L$).  Thus, $\H_i(L)={\mathcal
H}_i(\Sigma_L)\oplus {\h}_i(\Sigma_L )$.  Let $j\colon\{\pm
1\}\times\Sigma_L\to\Sigma_{CL}$ be the inclusion.  A class of the
form $(\alpha , -\alpha )$ in the direct sum obviously goes to $0$ in
${\h}_i(CL)$, while $j_\ast$ maps the diagonal subspace of
elements of the form $(\alpha , \alpha )$ isomorphically onto
${\h}_i(CL)$.  Statement (3) follows (as do formulas (1) and
(2)).
\end{proof}

\subsection{Poincar\'e duality}\label{7.4}
If a flag complex $S$ is a generalized homology sphere with rational
coefficients, then $\Sigma_S$ is a polyhedral homology manifold with
rational coefficients.  
Hence, $\Sigma_S$
satisfies Poincar\'e duality, ~\ref{2.7.1}.  Similarly, if a pair 
$(D, \partial D)$ of
flag complexes is a generalized homology disk with rational
coefficients, then $\Sigma_D$ is a polyhedral homology manifold with
boundary with rational coefficients (its boundary being $W_D
\Sigma_{\partial D}$) and hence, it satisfies the relative version of
Poincar\'e duality.

\subsubsection{}\label{7.4.1}  
If $S$ is a $GHS^{n-1}$, then $\beta_i(S)=\beta_{n-i}(S)$.

\subsubsection{}\label{7.4.2}  
If $(D, \partial D)$ is a $GHD^{n-1}$, then $\beta_i
(D, \partial D)=\beta_{n-i}(D)$.

\subsubsection{}\label{7.4.3}  If $(D, \partial D)$ is a $GHD^{n-1}$, then the
homology and cohomology sequences of the pair $(D, \partial D)$ are
isomorphic under Poincar\'e duality in the sense that the following
diagram commutes up to sign, 
$$
\minCDarrowwidth 0.5cm
\begin{CD}
@>>> \H_{i+1}(D, \partial D)  @>>> \H_i(\partial D)  @>>>  \H_i(D) @>>>
\H_i(D, \partial D)  @>>>\\ && \Big\updownarrow\cong &&
\Big\updownarrow\cong && \Big\updownarrow\cong &&
\Big\updownarrow\cong  \\ {}  @>>>  \H^{n-i-1}(D)  @>>>
\H^{n-i-1}(\partial D)  @>>>  \H^{n-i} (D, \partial D)  @>>>
\H^{n-i}(D)  @>>>
\end{CD}
$$
where the vertical isomorphisms are given by Poincar\'e duality 
(cf Theorem 1.1.5 of \cite{Bro}). 
 From
this we deduce the following lemmas which we shall need in
Section~\ref{9}.

\begin{Lemma}\label{7.4.4} Suppose $(D, \partial D)$ is a $GHD^{2k}$ and that
$j_\ast\colon \H_k(\partial D)\to \H_k(D)$ is the map induced by the
inclusion.  Then
$$\dim_{W_L}(\Ker j_\ast )=\dim_{W_L}(\overline{\im
j_\ast})=\frac{1}{2} \beta_k(\partial D) .$$
\end{Lemma}

\begin{proof}  By~\ref{7.4.3}, the sequences
$$\H_{k+1}(D, \partial D)\xrightarrow{\partial_\ast}\H_k
(\partial D)\xrightarrow{j_\ast}\H_k(D)$$ and
$$\H^k(D)\xrightarrow{j^\ast}\H^k(\partial D) \xrightarrow{\partial^\ast}
\H^{k+1}(D, \partial D)$$
 are
isomorphic.  In other words, under Poincar\'e duality, the connecting
homomorphism $\partial_\ast$ is isomorphic to $j^\ast$.  Since
$j^\ast$ is the adjoint of $j_\ast$, \ref{3.2.8} implies that
$$\beta_k(\partial D) = \dim_{W_L}(\Ker j_\ast
)+\dim_{W_L}(\overline{\im \partial_\ast}) .$$ 
 By the exact sequence
of the pair, Lemma~\ref{7.2.1},  $\Ker j_\ast = \overline{\im
\partial_\ast}$, so \linebreak 
$\dim_{W_L}(\Ker j_\ast
)=\frac{1}{2}\beta_k(\partial D)$.  Then, by~\ref{3.2.7}, we also have
$ \dim_{W_L}(\overline{\im j_\ast})=\frac{1}{2} \beta_k(\partial D)$, 
which proves the lemma.
\end{proof}

\subsubsection{}\label{7.4.5}  
Suppose that $S=D_1\cup D_2$ and $S_0=D_1\cap D_2$.
Also suppose that $S$ is a $GHS^{n-1}$ and that $(D_1, S_0)$ and
$(D_2, S_0)$ are $GHD^{n-1}$'s.  By Lemma~\ref{7.2.3}(2), 
$\H^i(S, S_0)\cong\H^i(D_1, S_0)\oplus\H^i(D_2, S_0)$.  Similarly
to~\ref{7.4.3}, the homology Mayer--Vietoris sequence of $S=D_1\cup
D_2$ is isomorphic, via Poincar\'e duality, to the exact sequence of
the pair $(S, S_0)$ in cohomology.  In other words, the following
diagram commutes up to sign, 
$$
\minCDarrowwidth 0.5cm
\begin{CD} 
@>>> \H_{i+1}(S) @>>> \H_i(S_0) @>>> \H_i(D_1) \oplus\H_i(D_2) @>>>\\
&&\Big\updownarrow\cong   &&\Big\updownarrow\cong
&&\Big\updownarrow\cong  \\ @>>> \H^{n-i-1}(S) @>>> \H^{n-i-1}(S_0)
@>>> \H^{n-i}(D_1, S_0)\oplus\H^{n-i}(D_2, S_0) @>>>\end{CD}$$ where
the first row is the Mayer--Vietoris sequence, the second is the exact
sequence of the pair and the vertical isomorphisms are given by
Poincar\'e duality.  We record the special case of this where $n=2k+1$
and $i=k$ as the following lemma.

\begin{Lemma}\label{7.4.6} 
With hypotheses as in~\ref{7.4.5}, suppose $n=2k+1$.  Then
the map $i_\ast\colon \H_k(S_0)\to\H_k(S)$ induced by the inclusion is
dual  (under Poincar\'e duality) to the connecting homomorphism
$\partial_\ast\colon\H_{k +1}(S)\to\H_k(S_0)$ in the Mayer--Vietoris
sequence.
\end{Lemma}

\begin{proof}  
In this special case, the diagram in~\ref{7.4.5} becomes the following:
$$\begin{CD}  \H_{k+1}(S) @>\partial_\ast>> \H_k (S_0) @>>>
\H_k(D_1)\oplus\H_k(D_2)\\ \Big\updownarrow\cong
&&\Big\updownarrow\cong   &&\Big\updownarrow\cong \\ \H^k(S)
@>i^\ast>> \H^k(S_0) @>>> \H^{k+1}(D_1, S_0)\oplus\H^{k+1}(D_2, S_0)
\end{CD}$$\vspace{-5mm}
\end{proof}

\section{\label{8}Variations on Singer's Conjecture}

In this section, we will consider several conjectures, $\I(n)$, 
$\II(n)$, $\III(n)$, $\IV(n)$ and $\V(n)$, concerning the reduced
$\ell^2$--homology of $\Sigma_L$, where $L$ is either a generalized
homology sphere, usually denoted by $S$, or a generalized homology
disk, denoted by $D$.  Here the ``$n$'' refers to the dimension of
$\Sigma_L$, where $L=S$ or $D$, so that $\dim L=n-1$.

As usual all simplicial complexes are flag complexes and all
subcomplexes are full.

\subsection{Restatement of Singer's Conjecture}\label{8.1}

\begin{I(n)}  If $S$ is a $GHS^{n-1}$, then $\b_i(S)=0$ for all $i\not=
\frac{n}{2}$.
\end{I(n)}

\subsection{Singer's Conjecture for a disk}\label{8.2}

\begin{II(n)}  Suppose $(D, \partial D)$ is a $GHD^{n-1}$.
\begin{itemize}
\item If $n=2k$ is even, then $\b_i(D)=\b_i(D, \partial D)=0$ for all
$i\not= k$.

\item If $n=2k+1$ is odd, then
\begin{enumerate}
\item $\beta_i(D)=\b_{i+1}(D, \partial D)=0$ for all $i\not= k$, and
\item $\beta_k(D)=\beta_{k+1}(D, \partial
D)=\frac{1}{2}\beta_k(\partial D)$ and the following sequence of the
pair is weakly short exact, 
$$0\to\H_{k+1}(D, \partial D)\to\H_k(\partial D)\to\H_k(D)\to 0 .$$
\end{enumerate}
\end{itemize}
\end{II(n)}

\subsubsection{}\label{8.2.1} 
Given a $GHD$, $(D, \partial D )$, let $S$ denote
the $GHS$ formed by gluing on $C(\partial D)$ to $D$ along $\partial
D$.  If $v$ denotes the cone point, then $\partial D=S_v$ (the link of
$v$) and $C(\partial D)=CS_v$.  Conversely, given a $GHS$, call it
$S$, and a vertex $v$, we obtain a $GHD$, $D=S-v$ with $\partial
D=S_v$.

Next we consider some seemingly weaker statements in odd dimensions.

\subsection{A weak form of the conjecture}\label{8.3}
\begin{III(2k+1)}  Suppose $(D, S_v)$ is a $GHD^{2k}$ and that $S=D\cup
CS_v$ is as in~\ref{8.2.1}.  Then in the Mayer--Vietoris sequence, the
map, 
$$j_\ast\oplus h_\ast\colon\H_k(S_v)\to\H_k(D)\oplus\H_k (CS_v) , $$ is
a monomorphism.
\end{III(2k+1)}

\subsubsection{Remark}\label{8.3.1}  By Lemma~\ref{7.3.3}, 
$h_\ast\colon \H_k(S_v)\to \H_k(CS_v)$ is surjective and the von
Neumann dimension of its kernel is $\frac{1}{2}\beta_k(S_v)$.
Similarly, by Lemma 7.4.4, $\dim_{W_S}(\Ker j_\ast
)=\frac{1}{2}\beta_k(S_v)$.  So, it is not unreasonable to expect that
these subspaces intersect in general position: $\Ker j_\ast\cap \Ker
h_\ast =0$, in other words, that $\III(2k+1)$ is valid.

\subsubsection{}\label{8.3.2} 
By Lemma~\ref{7.4.6}, $\III(2k+1 )$ is equivalent to the
following.
\begin{III'(2k+1)}  Suppose $(D, S_v)$ is a $GHD^{2k}$ and that $S=D
\cup CS_v$ is as in~\ref{8.2.1}.  Then the map
$i_\ast\colon\H_k(S_v)\to \H_k(S)$, induced by the inclusion, is the
zero homomorphism.
\end{III'(2k+1)}

\subsection{A weak form of the  conjecture for a disk}\label{8.4}
\begin{IV(2k+1)}  Suppose $(D, \partial D)$ is a $GHD^{2k}$.  Then
$\b_{k+1}(D)=0$.
\end{IV(2k+1)}

\subsection{The strong form of the conjecture}\label{8.5}  
The formulation of the last conjecture is key to our approach.

\begin{V(n)}  Suppose $S$ is a $GHS^{n-1}$ and that $A$ is any full
subcomplex.
\begin{itemize}
\item If $n=2k$ is even, then $\b_i(S, A)=0$ for all $i>k$.
\item If $n=2k+1$ is odd, then $\b_i(A)=0$ for all $i>k$.
\end{itemize}
\end{V(n)}

\subsection{Joins}\label{8.6}  It follows from the K\"unneth Formula, 
Lemma~\ref{7.2.4}, that the above conjectures are compatible with the
operation of taking joins.  For example, let $\J$ stand for $\I$, 
$\III$ or $\V$.  If $S_1$ and $S_2$ are $GHS$'s of dimension $n_1-1$
and $n_2-1$ for which $\J(n_1)$ and $\J(n_2)$ hold, then $\J(n_1+n_2)$
holds for $S_1\ast S_2$ (which is a $GHS^{n_1+n_2-1}$ by~\ref{4.3.7}).  
Similarly, let $\J$ stand for $\II$ or
$\IV$.  If $S_1$ is a $GHS^{n_1-1}$ for which $\I(n_1)$ holds and
$D_2$ is $GHD^{n_2-1}$ for which $\J(n_2)$ holds, then $\J(n_1+n_2)$
holds for $S_1\ast D_2$ (which is a $GHD^{n_1+n_2-1}$ by 4.3.7).

\subsection{An $\boldsymbol m$--gon}\label{8.7}  
Suppose $S$ is an $m$--gon, $m\ge 4$.
Then by~\ref{7.2.6}, $\b_0(S)=0$ and by~\ref{7.2.7} (or
by~\ref{7.3.1})  $\b_2 (S)=0$.  So, $\I(2)$ holds.  Similarly, 
$\b_0(\BS^0)=\b_1(\BS^0)=0$, so $\I(1)$ holds.

\subsection{Some implications}\label{8.8}
Next we list some obvious implications among\break these
conjectures.

\subsubsection{}\label{8.8.1}Let $\J$ stand for $\I$, $\II$, $\III$, 
$\III'$, $\IV$, or $\V$. Then $\J(n)\implies \J(n-2)$.

\begin{proof}  Suppose $L^{n-3}$ is a $GHS$ or $GHD$ for which $\J(n-2)$ fails.
Let $S_1$ be a $5$--gon.  Since $\beta_1(S_1)\not= 0$,  the K\"unneth
Formula~\ref{7.2.4} shows that $\J(n)$ also fails for $S_1\ast
L^{n-3}$.
\end{proof}

\subsubsection{}\label{8.8.2} $\II(n)\implies \I(n-1)$.
\begin{proof}  
Let $S$ be a $GHS^{n-2}$.  By Lemma~\ref{7.3.3}, if $\II(n)$ holds
for $(CS, S)$, then $\I(n-1)$ holds for $S$.
\end{proof}

\subsubsection{}\label{8.8.3} $[\I(n-1)$ and $\I(n)]\implies \II(n)$.
\begin{proof}  Suppose $\I(n-1)$ and $\I(n)$ hold, that $(D, \partial D)$ is a
$GHD^{n-1}$ and, as in~\ref{8.2.1}, that $S=D\cup C(\partial D)$.  If
$n=2k$, then since $\I(2k-1)$ holds for $\partial D$, $\b_i(\partial
D)=0$ for all $i$.  By $\I(2k)$, $\b_i(S)=0$ for $i\not= k$.  The
Mayer--Vietoris sequence then yields that $\b_i(D)=0$ for $i\not= k$
and that $\H_k (D)\oplus \H_k(C(\partial D))\cong\H_k(S)$.  So, 
$\II(2k)$ holds for $(D, \partial D)$.  If $n=2k+1$, then by
$\I(2k+1)$, $\b_i(S)=0$ for all $i$.  Hence, the Mayer--Vietoris
sequence yields, $\H_i(\partial D)\cong \H_i(D)\oplus\H_i(C(\partial
D))$.  By $\I(2k)$, $\b_k(\partial D)=0$ for $i\not= k$.  It then
follows from~\ref{7.4.2} and Lemma~\ref{7.4.4} that $\II(2k+1)$ holds
for $(D, \partial D)$.
\end{proof}

\subsubsection{}\label{8.8.4} $\II(2k)\implies \I(2k)$.
\begin{proof}  Let $S$ be a $GHS^{2k-1}$ and $v$ a vertex of $S$.  Write
$S=D\cup CS_v$, as in~\ref{8.2.1}.  Assume $\II(2k)$ holds.
By~\ref{8.8.2}, $\I(2k-1)$ holds for $S_v$, ie, $\b_i (S_v)=0$ for
all $i$.  From the sequence of the pair $(S, S_v)$ we get:
$\H_i(S)\cong\H_i(S, S_v)\cong\H_i(D, S_v)\oplus\H_i(CS_v, S_v)$.
Since $\II(2k)$ holds, the last two
terms are nonzero only in the middle dimension.  Hence, 
$\II(2k)\implies \I(2k)$.
\end{proof}

The next two implications, \ref{8.8.5} and \ref{8.8.6}, are immediate.

\subsubsection{}\label{8.8.5} $\I(2k+1)\implies \III(2k+1 )$.

\subsubsection{}\label{8.8.6} $\II(2k+1)\implies \IV(2k+1 )$.

\subsubsection{}\label{8.8.7} $\V(n)\implies \I(n )$.

\begin{proof}  This follows from Poincar\'e duality.  (If $n=2k$, take
$A=\emptyset$ to get $\b_i(S)=0$ for $i>k$.  If $n=2k+1$, take $A=S$, 
to get $\b_i(S)=0$ for $i>k$.)
\end{proof}

\subsubsection{}\label{8.8.8} $\V(n)\implies \II(n )$.

\begin{proof}  We proceed as in~\ref{8.8.3}.  Given $(D, \partial D)$, set $S=
D\cup C(\partial D)$, as before.  Assume $\V(n)$ holds.  If $n=2k$, 
then $\b_i(S)=0$ for $i\not= k$ (by~\ref{8.8.7}).  For $i>k$, by
$\V(2k)$, we have that
$$0=\H_i(S, \partial D)\cong\H_i(D, \partial D)\oplus\H_i (C(\partial
D), \partial D).$$ So, for $i>k$, $\b_i(D, \partial D)=0$ and
$\b_{i-1}(\partial D)=0$ (by Lemma~\ref{7.3.3}(2)).  By~\ref{7.4.1}, 
the second equation implies that $\b_i(\partial D)=0$ for all $i$.  It
then follows from the exact sequence of the pair and~\ref{7.4.2}, that
$\II(2k)$ holds for $(D, \partial D)$.  If $n=2k+1$, then $\b_i(S)=0$
for all $i$ (by~\ref{8.8.7}).  The sequence of the pair $(S, \partial
D)$ then gives:
$$\H_i(D, \partial D)\oplus\H_i(C(\partial D), \partial D)\cong\H_{i-1}
(\partial D).$$ By $\V(2k+1)$, $\b_j(\partial D)=0$ for $j>k$; hence, 
by 7.4.1, it vanishes for $j\not= k$.  Therefore, $\b_i(D, \partial
D)=0$ for $i\not= k+1$ and by~\ref{7.4.2}, $\b_i(D)=0$ for $i\not= k$.
So, $\II(2k+1)$ holds for $(D, \partial D)$.
\end{proof}

\begin{Lemma}\label{8.8.9}Statement $\V(2k)$ implies that for any full
subcomplex $A$ of $S$ (a $GHS^{2k-1}$), we have
$$\b_i(A)=0 \quad \text{  for all $i>k$. }$$
\end{Lemma}

\begin{proof}  Assume $\V(2k)$ holds.  By~\ref{8.8.7}, $\b_i(S)=0$ for $i\not=
k$.  Hence, in the exact sequence of the pair, 
$$\H_{i+1}(S, A)\to\H_i(A)\to\H_i(S) , $$ the first and third terms
vanish for all $i>k$.
\end{proof}

\subsection{A conjecture for groups of finite type}\label{8.9}Recent
work by Bestvina, Kapovich and Kleiner in \cite{BKK} shows that
if  $A=P_3*\dots *P_3$ is an $k$--fold join of $3$ points with itself,
then $W_A$ cannot act properly on a contractible $(2k-1)$--manifold. 
Their argument uses the well-known fact that $A$ does
not embed in $\BS^{2k-2}$.  We note that this well-known fact follows
from Conjecture $\V(2k-1)$ since, by~\ref{7.2.4}, $b_k^{(2)}(W_A) \neq 0$
which would contradict Conjecture $\V(2k-1)$ if $A$ were a full
subcomplex of some flag triangulation of $\BS^{2k-2}$.
 (Here the $\ell^2$--Betti numbers $b^{(2)}_i(W_A)$ are as defined
in~\ref{3.3.7}.)  These remarks suggest the following generalization
of Singer's Conjecture.
 
\begin{Conjecture}\label{8.9.1}Suppose that a discrete group $G$ acts
properly on a contractible $n$--manifold.  Then
$$b^{(2)}_i(G)=0 \quad \text{   for $i>\frac{n}{2}$. }$$
\end{Conjecture}  (In the case where $G$ does
not act cocompactly on its universal space $\underline{EG}$, define
its $\ell^2$--Betti numbers as in \cite{CG2}.)

\section{\label{9}Inductive Arguments}

We describe a partially successful program for proving conjecture
$\V(n)$.  The idea is to use a double induction:  first, induction on
the dimension $n$ and second, depending on the parity of $n$, 
induction either on the number of vertices of $A$ or on the number of
vertices in $S-A$.

\subsection{Notation}\label{9.1}  We set up some notation for the
induction on the number of vertices.  Suppose $A$ and $B$ are full
subcomplexes of $S$, the vertex sets of which differ by only one
element, say $v$.  In other words, $B=A-v$, for some $v\in{\mathcal
S}_0(A)$.  Let $A_v$ and $S_v$ denote the link of $v$ in $A$ and $S$, 
respectively.  Thus, $A=B\cup CA_v$ and $CA_v\cap B=A_v$.  We note
that $S_v$ is a $GHS$ of one less dimension than $S$ and that $A_v$ is
a full subcomplex of $S_v$.

\subsection{Induction on the number of vertices}\label{9.2}

\begin{Lemma}\label{9.2.1}$\V(2k-1)\implies \V(2k)$.
\end{Lemma}

\begin{proof}  Suppose $\V(2k-1)$ holds.  Let $(S, A)$ be as in $\V(2k)$ and
let $B=A-v$.  Assume, by induction on the number of vertices in $S-A$, 
that $\V(2k)$ holds for $(S, A)$.  (The case $A=S$ being trivial.)  We
want to prove it also holds for $(S, B)$, ie, that $\b_i(S, B)=0$ for
$i>k$.  Consider the exact sequence of the triple $(S, A, B)$:
$$\to\H_i(A, B)\to\H_i(S, B)\to\H_i(S, A)\to .$$ 
Suppose $i>k$.  By
inductive hypothesis, $\b_i(S, A)=0$.  By excision, Lemma \ref{7.2.2}, 
$\beta_i(A, B)=\beta_i(CA_v, A_v)$.  By Lemma~\ref{7.3.3}(2), $\beta_i
(CA_v, A_v)=\frac{1}{2}\beta_{i-1}(A_v)$.  Since $\V(2k-1)$ holds for
$(S_v, A_v)$ and since $i-1>k-1$, $\beta_{i-1}(A_v)=0$.  So, 
$0=\beta_i(CA_v, A_v) =\beta_i(A, B)$.  Consequently, $\b_i(S, B)=0$.
\end{proof}

Essentially the same argument proves the following lemma (which we
will need in Section~\ref{11.4}).

\begin{Lemma}\label{9.2.2}Assume that $\V(2k)$ holds.  Suppose that a flag
complex $L$ is a polyhedral homology manifold of dimension $2k$ and
that $A$ is a full subcomplex.  Then $\b_i(L, A)=0$ for $i>k+1$.
\end{Lemma}

\begin{proof}  We proceed as in the previous proof.  If $B=A-v$, then
$\beta_i(A, B)=\beta_i(CA_v, A_v)=\frac{1}{2}\beta_{i-1} (A_v)$.  Since
we are assuming $\V(2k)$ holds, Lemma~\ref{8.8.9} implies that
$\beta_{i-1} (A_v)=0$ for $i>k+1$.  Hence, if we assume by induction
that the lemma holds for $(L, A)$, then it also holds for $(L, B)$.
\end{proof}

\begin{Lemma}\label{9.2.3}
$\left[\V(2k)\text{ and }\III(2k+1)\right]\implies \V(2k+1) .$
\end{Lemma}

\begin{proof} 
 Assume $\V(2k)$ and $\III(2k+1)$ hold.  Let $(S, A)$ be as in $\V(2k
+1)$ and let $B=A-v$.  Assume, by induction on the number of vertices
in $B$, that $\V(2k+1)$ holds for $B$.  (The case $B=\emptyset$ being
trivial.)  We want to prove that it also holds for $A$, ie, that
$\b_i(A)=0$ for $i>k$.

First suppose that $i>k+1$.  Consider the Mayer--Vietoris sequence for
$A=B\cup CA_v$:
$$\H_i(B)\oplus\H_i(CA_v)\to\H_i(A)\to\H_{i-1} (A_v) .$$ By $\V(2k)$
and Lemma~\ref{8.8.9}, $\b_{i-1}(A_v)=0$ (since $i-1>k$) and hence, 
$\b_i(CA_v)=0$ (by Lemma~\ref{7.3.3}(1)).  By inductive hypothesis, 
$\b_i(B)=0$, and consequently, $\b_i(A)=0$.

For $i=k+1$, we compare the Mayer--Vietoris sequence of $A=B\cup CA_v$
with that of $S=D\cup CS_v$ (where $D=S-v$):
$$\begin{CD}  & & & &\H_{k+1}(S_v, A_v)\\ & & & & @VVV\\ 0 @>>>
\H_{k+1}(A) @>>> \H_k(A_v) @>{j'_\ast \oplus h'_\ast}>>
\H_k(B)\oplus\H_k(CA_v)\\ & & & & @V{f_\ast}VV  @VVV \\ & & & &
\H_k(S_v) @>{j_\ast\oplus  h_\ast}>> \H_k (D)\oplus
\H_k(CS_v)\end{CD}$$ By $\V(2k)$, $\b_{k+1}(S_v, A_v)=0$; hence, 
$f_\ast$ is injective.  By $\III(2k+1)$, $j_\ast\oplus h_\ast$ is
injective.  Hence,  $j'_\ast\oplus h'_\ast$ is injective and therefore, 
$\b_{k+1}(A)=0$.
\end{proof}

\subsection{Induction on dimension}\label{9.3}  Our main result is the
following.

\begin{Theorem}\label{9.3.1} Statement $\III(2k-1)$ implies that $\V(n)$ holds
for all $n \leq 2k$.
\end{Theorem}

\begin{proof}  
By~\ref{8.2.1}, $\III(2k-1)$ implies $\III(2l -1)$, for all $l\le
k$.  Suppose, by induction on $n$, that $\V(n-1)$ holds for some $n\le
2k$.  If $n-1$ is odd, then by Lemma~\ref{9.2.1}, $\V(n-1)$ implies
$\V(n)$.  If $n-1$ is even, then by Lemma~\ref{9.2.3}, $\V(n-1)$ and
$\III(n)$ imply $\V(n)$.
\end{proof}

\subsubsection{}\label{9.3.2} An $n$--dimensional simplicial complex $L$ has
{\it spherical links in codimensions} $\le m$ if for each
$\sigma\in{\mathcal S}(L)$ of dimension $n-i$, with $i\le m$, its link
$\Link(\sigma , L)$ is a $GHS^{i-1}$.

When $\dim L<m$, this condition means that $L$ is a $GHS$ (take
$\sigma = \emptyset$).  When $\dim L\ge m$, it means that the
complement of its codimension--$(m+1)$ skeleton is a homology manifold.
For example, for $m=1$, it means that $L$ is a pseudomanifold.

We note that the condition is inherited by links of vertices:  if $L$
has spherical links in codimensions $\le m$, then so does $L_v$ for
any vertex $v$.

\begin{Theorem}\label{9.3.3} Assume that $\III(2l +1)$ is true.  Let $L$ be
an $(n-1)$--dimensional flag complex, $n\ge 2l +1$, with spherical
links in codimensions $\le 2l +1$.  Then for any full subcomplex $A$
of $L$, 
$$\beta_{n-i}(A)=0 \quad \text{ for $i\le l$, }$$
and if $n\not= 2l +1$, then
$$\beta_{n-i}(L, A)=0 \quad \text{  for $i\le l$. }$$
\end{Theorem}

\begin{proof}  The proof is by induction on $n$, starting at $n=2l +1$.  If
$n=2l +1$, then $L$ is a $GHS^{2l}$ and the result follows from
Theorem~\ref{9.3.1}.

If $n> 2l +1$, then for any vertex $v \in A$, we have, by inductive hypothesis
that $\beta_{(n-1)-i}(A_v)=0$ for $i\le l$.  The proof of
Lemma~\ref{9.2.1} then shows that $\beta_{n-i}(L, A)=0$ and the first
part of the proof of Lemma~\ref{9.2.3} shows that $\beta_{n-i}(A)=0$
for $i\le l$.
\end{proof}

\begin{Corollary}\label{9.3.4}Assume that $\III(2l +1)$ is true.  Suppose
that a flag complex $S$ is a $GHS^{n-1}$, $n\ge 2l+1$. Then
$$\beta_{i}(S)=\beta_{n-i}(S)=0 \quad \text{  for $i\le l$. }$$
\end{Corollary}

\begin{proof}  This follows from the previous theorem (taking $L=A=S$) and
Poincar\'e duality.
\end{proof}

Theorem~\ref{9.3.3} suggests the following generalization of Singer's
Conjecture.

\begin{Conjecture}\label{9.3.5} Suppose that $X$ is a contractible, geometric
$G$--com\-plex of dimension $n$ and that $X$ has spherical links in
codimensions $\le 2l +1$, where $2l +1\le n$.  
Then 
$$b^{(2)}_{n-i}(X;G)=0 \quad \text{ 
for $i\le l$. }$$
\end{Conjecture}

\section{\label{10}The conjecture in dimension $3$}

\subsection{Review of previous results}\label{10.1}  In \cite{LL} Lott and
L\"uck proved Singer's Conjecture for any closed, irreducible
$3$--manifold with infinite fundamental group for which
Thurston's Geometrization Conjecture holds.  In other words, for such
a $3$--manifold, the reduced $\ell^2$--homology of its universal cover
vanishes.  As we shall see in~\ref{10.1.5}, below, the Geometrization
Conjecture holds for the $3$--dimensional orbifolds which we are
interested in.  Hence, conjecture $\I(3)$, from Section~\ref{8.1}, is
true.

The calculation in \cite{LL} depends on the following two facts, 
stated as~\ref{10.1.1} and~\ref{10.1.2}, below.

\subsubsection{}\label{10.1.1}  
Suppose $M$ is the compact $3$--manifold formed by
gluing together two compact $3$--manifolds $M_1$ and $M_2$ along one or
more boundary components which are incompressible tori.  If the
reduced $\ell^2$--homology of their universal covers, $\widetilde{M_1}$
and $\widetilde{M_2}$, vanishes, then so does the reduced
$\ell^2$--homology of $\widetilde{M}$.  (This follows the
Mayer--Vietoris sequence, ~\ref{2.4.4}, and the vanishing of the reduced
$\ell^2$--homology of the universal cover of $T^2$.)

\subsubsection{Theorem 5.14 of \cite{LL}\label{10.1.2}}  Let $M^3$ be a
compact $3$--manifold with boundary such that its interior is
homeomorphic to a complete hyperbolic manifold of finite volume.  Then
the reduced $\ell^2$--homology of its universal cover vanishes.

We note that if $M^3$ is a compact $3$--manifold formed by chopping off
the cusps of a complete hyperbolic $3$--manifold of finite volume, then
each boundary component of $M^3$ is a $2$--torus (or possibly a Klein
bottle in the nonorientable case).  Hence, the result of Lott--L\"uck
follows from 10.1.1, ~\ref{10.1.2} and a similar result for aspherical
Seifert fiber spaces.

The result in~\ref{10.1.2} is, in turn, a consequence of the next two
facts, stated as~\ref{10.1.3} and~\ref{10.1.4}, below.

\subsubsection{}\label{10.1.3}  
The reduced $\ell^2$--homology of any odd-dimensional
hyperbolic space, $\BH^{2k+1}$, vanishes.  (This is proved in
\cite{Do2}.)

The next result is stated on page 226 of \cite{G3}.  It is proved by
Cheeger and Gromov in \cite{CG1}.

\subsubsection{Bounded geometry}\label{10.1.4}  Suppose $X$ is a complete
contractible Riemannian manifold with uniformly bounded geometry, 
ie, its sectional curvature is bounded and its injectivity radius is
bounded away from $0$.  Let $\Gamma$ be a discrete group of isometries
of $X$ with $\text{Vol}  X/\Gamma<\infty$.  Then
$$b^{(2)}_k(\Gamma )=\dim_\Gamma ({\h}_k(X)) .$$ (Here
$b^{(2)}_k(\Gamma )$ is the $\ell^2$--Betti number of $\Gamma$ defined
in 3.3.7.)

Thus, ~\ref{10.1.2}  follows from~\ref{10.1.4} and~\ref{10.1.3} in the
case where $X=\BH^3$.

\subsubsection{Haken manifolds}\label{10.1.5}  Thurston proved that the
Geometrization Conjecture holds for Haken $3$--manifolds.  Suppose that
$S$ is a triangulation of the $2$--sphere as a flag complex and that
$M^3 =P_S$, the commutator cover of $\Sigma_S/W_S$ considered
in~\ref{6.4}.  Then $M^3$ is obviously Haken.  Indeed, for any vertex
$s$ of $S$, the special subcomplex $\Sigma_{S_s}$ is geodesically
convex; hence, its image in $M^3$ is an incompressible surface. 
(See the argument in \ref{14.1.6}, below.)  Therefore, the reduced
$\ell^2$--homology of $\Sigma_S$ (the universal cover of $M^3$)
vanishes.

In this special case, Thurston's Theorem is basically a consequence of
Andreev's Theorem, which was proved several years earlier in
\cite{An},\cite{An1}.  
We explain Andreev's Theorem in Section~\ref{10.3}, below.
However, we first need to develop some material about triangulations
of the $2$--sphere.

\subsection{Triangulations of $\BS^2$}\label{10.2}  Let $S$ be a
triangulation of $\BS^2$ as a flag complex.

\subsubsection{}\label{10.2.1}  The {\it valence} of a vertex $s$ of $S$ is the number
of vertices in its link.  In what follows we shall be concerned with
the vertices of valence $4$.

\subsubsection{}\label{10.2.2}  Let $C$ be a circuit of length $4$ in the $1$--skeleton
of $S$.  Then $C$ is an {\it empty $4$--circuit} if (a) $C$ is not the
link of a vertex and (b) $C$ is not the boundary of the union of two
adjacent $2$--simplices.  Since $S$ is a flag complex, it follows from
(b) that any empty $4$--circuit $C$ is a full subcomplex.

\begin{Lemma}\label{10.2.3}Suppose that
\begin{itemize}
\item{\rm(i)} $S$ has no empty $4$--circuits and
\item{\rm(ii)} $S$ is not the suspension of a $4$-- or $5$--gon.
\end{itemize}
Then no two valence $4$ vertices of $S$ are connected by an edge.
\end{Lemma}

\begin{proof}  Suppose that $s_1$ and $s_2$ are valence $4$ vertices which
are connected by an edge.  Then the star of that edge is the
configuration pictured in the figure below.

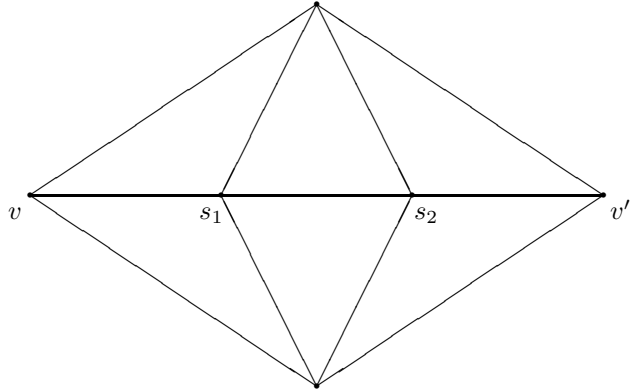
\setlength{\unitlength}{0.25in}
\begin{picture}(16, 10)(-10.3, -5) \small
\put(-6, 0){\line(1, 0){12}} \put(-6, 0){\line(3, 2){6}}
\put(-6, 0){\line(3, -2){6}} \put(6, 0){\line(-3, 2){6}}
\put(6, 0){\line(-3, -2){6}} \put(0, 4){\line(1, -2){2}}
\put(0, 4){\line(-1, -2){2}} \put(0, -4){\line(1, 2){2}}
\put(0, -4){\line(-1, 2){2}}

\put(-6.5, -0.5){$v$} \put(-2.5, -0.5){$s_1$} \put(2, -0.5){$s_2$}
\put(6.1, -0.5){$v'$}

\put(-6, 0){\circle*{.1}} \put(-2, 0){\circle*{.1}}
\put(2, 0){\circle*{.1}} \put(6, 0){\circle*{.1}}
\put(0, 4){\circle*{.1}} \put(0, -4){\circle*{.1}}

\end{picture}

The indicated vertices $v$ and $v'$ cannot coincide, since if they did $S$ 
would
contain an empty $3$--circuit and hence, not be a flag complex.
Similarly,
the top and bottom vertices cannot be connected by an edge, since $S$ would
again contain empty $3$--circuits.
 Let $C$ be the boundary of the star in
the figure.  If $C$ is the boundary of two adjacent $2$--simplices, then $S$
is the suspension of a $4$--gon.  If $C$ is the link of a missing vertex,
then $S$ is the suspension of a $5$--gon.  Otherwise, $C$ is an empty
$4$--circuit, contradicting (i).
\end{proof}

\begin{Lemma}\label{10.2.4}Let $T$ be a set of valence $4$ vertices of
$S$, no two of which are connected by an edge.  Then
$\beta_i(S)=\beta_i (S-T)$ for all $i$.
\end{Lemma}

\begin{proof}  By \ref{7.3.7} $\H_\ast (S_s)$ vanishes for any $s\in T$.  
Hence, it
follows from the Mayer--Vietoris sequence, Lemma~\ref{7.2.3}(1), that
we can adjoin $CS_s$ to $S-T$ without changing $\beta_i$.
\end{proof}

\subsubsection{}\label{10.2.5}  For $j=1, 2$, suppose that $S_j$ is a flag
triangulation of $\BS^2$ and that $s_j$ is a vertex of valence $4$ in
$S_j$.  Choose an identification of the link of $s_1$ with that of
$s_2$.  (They are both $4$--gons.)  Define a new triangulation
$S_1\square S_2$ of $\BS^2$ by gluing together the $2$--disks $S_1-s_1$
and $S_2-s_2$ along their boundaries.

\subsubsection{}\label{10.2.6}  
Conversely, suppose $C$ is an empty $4$--circuit in
$S$.  Then $C$ separates $S$ into two $2$--disks, $D_1$ and $D_2$.  Let
$S_1$ and $S_2$ denote the result of capping off $D_1$ and $D_2$, 
respectively (where ``capping off'' means adjoining a cone on the
boundary).  Then $S=S_1\square S_2$.

The next lemma is a version of~\ref{10.1.1}.

\begin{Lemma}\label{10.2.7}$\beta_1 (S_1\square
S_2)=\beta_1(S_1)+\beta_1(S_2)$.  Thus, $\H_\ast$ vanishes for
$S_1\square S_2$ if and only if it vanishes for both $S_1$ and $S_2$.
\end{Lemma}

\begin{proof}  This follows from the Mayer--Vietoris sequence as before.
\end{proof}

\subsubsection{}\label{10.2.8}  
Suppose $S$ satisfies the conditions of Lemma~\ref{10.2.3} and let $T$  
denote the set of valence $4$ vertices of $S$. 
Consider a cellulation $[S-T]$ of $\BS^2$ obtained by replacing stars
of  vertices of $T$ by square cells. By Lemma~\ref{10.2.3}, $[S-T]$ is
a well-defined  $2$--complex homeomorphic to $\BS^2$ with triangular
and square faces.  In fact, it is easy to see that, under the
assumptions of  Lemma~\ref{10.2.3}, this complex is a {\it cell}
complex in a strict sense that any nonempty intersections of two cells
is a cell. It is a classical fact that any such complex 
is combinatorially dual to the boundary complex of a convex polytope,
which we will denote $K_{[S-T]}$.

\subsection{Andreev's Theorem}\label{10.3}  In \cite{An1} Andreev determined
which convex polytopes could occur as fundamental chambers of
classical reflection groups on $\BH^3$.  More precisely, given a
convex polytope with assigned dihedral angles in $\left( 0,
\frac{\pi}{2}\right]$ on the edges, he gave necessary and sufficient
conditions for it to be realized as a (possibly ideal) convex polytope
in $\BH^3$.  A special case of his result is the following.

\begin{Theorem}[Andreev]\label{10.3.1}  Suppose that $S$ is a flag
triangulation of $\BS^2$ and that
\begin{itemize}
\item{\rm(i)} $S$ has no empty $4$--circuits, and

\item{\rm(ii)} $S$ is not the suspension of a $4$-- or $5$--gon.
\end{itemize}
Let $T$  denote the set of valence $4$ vertices of $S$ and let
$K_{[S-T]}$  be the dual of the cellulation $[S-T]$ of $\BS^2$ obtained by
replacing stars of  vertices of $T$ by square cells.
 
Then $K_{[S-T]}$ can be realized as an ideal, right-angled convex
polytope in $\BH^3$.  (The ideal vertices correspond to the square faces
of $[S-T]$,  ie, to the vertices of valence $4$ in $S$.)  The
resulting classical reflection group is the right-angled Coxeter group
$W_{S-T}$.
\end{Theorem}

\begin{proof} By \ref{10.2.8},  $K_{[S-T]}$
is combinatorially equivalent to the boundary complex of a convex
polytope with vertices of valence $3$ and $4$ only. In Theorem~2
of~\cite{An1}, Andreev lists $6$ conditions
$\mathfrak{m}0$--$\mathfrak{m}5$ on assigned angles for such a
polytope to be realized in $\BH^3$.

The conditions $\mathfrak{m}0$ and $\mathfrak{m}1$ are immediate under
our hypothesis, since all angles are $\frac{\pi}{2}$. The remaining
conditions refer to certain configurations of faces of the polytope,
and turn out to be vacuous in our case, since these configurations
never appear under our hypothesis.

Indeed, since $S$ is a flag triangulation, it follows that  $[S-T]$
does not contain empty $3$--circuits, and therefore, $K_{[S-T]}$ does
not contain triangular prismatic elements and cannot be a triangular
prism. Similarly, since $S$ does not contain empty $4$--circuits, every
$4$--circuit in $[S-T]$ is a boundary of either two adjacent triangles
or of a square cell, and therefore, $K_{[S-T]}$ does not contain
quadrangular prismatic elements. Thus, we have verified conditions
$\mathfrak{m}2$, $\mathfrak{m}3$ and $\mathfrak{m}4$.

To verify condition $\mathfrak{m}5$ we note that if two faces of
$K_{[S-T]}$ intersect at a vertex, but are not adjacent, then this
vertex has to have valence $4$. So this vertex corresponds to a square
cell of $[S-T]$, and the two faces correspond to opposite corners of
the square.  The configuration in condition  $\mathfrak{m}5$ has a
third face, adjacent to both previous two, so the corresponding vertex
in $[S-T]$ is connected to these corners.  In $S$ this square is
subdivided by the diagonals and,  since $S$ does not contain empty
$4$--circuits, one of the remaining corners of the square in $S$ must
be connected to the vertex corresponding to the third face. This means
that $S$ contains a configuration pictured in Lemma~\ref{10.3.2},
which according to that lemma is impossible.
\end{proof}

\subsubsection{Remark}\label{10.3.2}  Thurston gives a proof of Andreev's
Theorem in \cite{T}.   
Hypothesis (ii) does
not occur in his statement of the result.  The reason is that
Thurston's statement is in terms of finding a collection of
half-spaces in $\BH^3$ with nonempty intersection such that their
supporting planes intersect in the prescribed combinatorial pattern
with the prescribed dihedral angles.  When all the dihedral angles are
strictly less than $\frac{\pi}{2}$, he shows that the intersection of
half-spaces is a (possibly ideal) polytope.  However, when some of the
angles $=\frac{\pi}{2}$, the intersection can degenerate to a lower
dimensional set.  In the case of interest, all the angles are
$\frac{\pi}{2}$.  It is easy to see that when the intersection is a
planar set, $S$ must be the suspension of a $5$--gon and similarly,
when it is $0$-- or $1$--dimensional, that $S$ is the suspension of a
$4$--gon.

\subsection{$\I(3)$ is true}\label{10.4}

\begin{Theorem}\label{10.4.1}Let $S$ be a triangulation of the $2$--sphere
as a flag complex.  Then $$\b_i (S)=0 \quad \text{  for all $i$. }$$
\end{Theorem}

\begin{proof}  If $S$ is the suspension of a $4$-- or $5$--gon, then the theorem
follows from Lemma~\ref{7.3.6}.  If $S$ is not the suspension of a
$4$--gon or a $5$--gon and if it has no empty $4$--circuits, then
by~\ref{10.1.3}, ~\ref{10.1.4} and~\ref{10.3.1}, $\H_i(S-T)$ vanishes
for all $i$, where $T$ denotes the set of valence $4$ vertices.
Hence, by Lemmas~\ref{10.2.3} and~\ref{10.2.4}, $\H_i(S)$ also
vanishes.

In every other case, $S$ has an empty $4$--circuit which we can use to
decompose $S$ as, $S=S_1\square S_2$, as in~\ref{10.2.6}.  Since $S_1$
and $S_2$ each have fewer vertices than does $S$, this process must
eventually terminate.  So, the theorem follows from Lemma~\ref{10.2.7}.
\end{proof}

\section{\label{11}Some consequences}

\subsection{$\V(3)$ and $\V(4)$ are true}\label{11.1}  Since $\I(3 )$ is
true, Theorem~\ref{9.3.1} (together with~\ref{8.8.5}) yields the
following.

\begin{Theorem}\label{11.1.1}
Statement $\V(n)$ (from~\ref{8.5}) is true for $n\le 4$.
\end{Theorem}

\subsection{$\boldsymbol 4$--dimensional consequences}\label{11.2}  
Since $\V(4)$
implies $\I(4)$ (by~\ref{8.8.7}), Sin\-ger's Conjecture holds for
$\Sigma_S$, where $S$ is any flag triangulation of a rational homology
$3$--sphere, ie, $\beta_1(S)=\beta_3(S)=0$.  By Atiyah's Formula, 
this implies that $\chi^{\orb} (\Sigma_S/W_S)=\beta_2(S)\ge 0$, and
hence, by~\ref{6.3.4}, that the Flag Complex Conjecture~\ref{0.2} is
true in dimension $3$.  We restate this as follows.

\begin{Theorem}[The Flag Complex Conjecture in dimension $3$]\label{11.2.1}
Let $S$ be any triangulation of a rational homology $3$--sphere as a
flag complex.  Then
$$\sum^3_{i=-1}\left(-\frac{1}{2}\right)^{i+1}f_i(S)\ge 0 , $$ where
$f_i(S)$ denotes the number of $i$--simplices in $S$ and where
$f_{-1}=1$.
\end{Theorem}

As explained in \cite{CD}, this implies the following $4$--dimensional
result.

\begin{Theorem}\label{11.2.2}
The Euler Characteristic Conjecture~\ref{0.1} holds for
all \linebreak nonpositively curved, piecewise Euclidean $4$--manifolds
which are cellulated by regular Euclidean cubes.  In other words, for
any such $4$--manifold $M^4$, $$\chi (M^4)\ge 0.$$
\end{Theorem}

In fact, one only need require $M^4$ to be a rational homology
$4$--manifold (rather than a $4$--manifold).

\subsection{Higher dimensional consequences}\label{11.3}  From Theorem
9.3.3 and Corollary~\ref{9.3.4}, we get the following.

\begin{Theorem}\label{11.3.1} 
Suppose $L$ is an $(n-1)$--dimensional flag complex, 
$n\ge 3$, with spherical links in codimensions $\le 3$.  Then for any
full subcomplex $A$ of $L$
$$\beta_n(A)=\beta_{n-1}(A)=0$$ 
and
$$\beta_n (L, A)=\beta_{n-1}(L, A)=0.$$
\end{Theorem}

\begin{Theorem}\label{11.3.2}Suppose $S$ is a $GHS^{n-1}$, $n\ge 3$.  Then
$$\beta_1 (S)=\beta_{n-1}(S)=0.$$
\end{Theorem}

\subsection{$\boldsymbol 3$--dimensional consequences}\label{11.4}  
We restate $\V(3)$
as follows.

\begin{Theorem}\label{11.4.1}Let $A$ be a finite flag complex of dimension $\le
2$.  Suppose $A$ is planar (ie, it can be embedded as a subcomplex
of the $2$--sphere).  Then $$\beta_2(A)=0.$$
\end{Theorem}

\begin{proof} By Lemma~\ref{7.3.4}, we may assume that $A$ is connected.
 Suppose that $A$ is piecewise linearly embedded in $\BS^2$. By
introducing a new vertex in the interior of each complementary region
and then coning off the boundary of each region, we obtain a flag
triangulation $S$ of the $2$--sphere with $A$ embedded as a full
subcomplex.  By $\V(3)$, $\beta_2(A)=0$.
\end{proof}

\subsubsection{Example}\label{11.4.2}  
The contrapositive of Theorem~\ref{11.4.1} states that if
$\beta_2(A)$\break$\not= 0$, then $A$ is not planar.  Kuratowski's graph
$K_{3, 3}$ is defined to be $P_3\ast P_3$, the join of $3$ points with
itself.  By Lemma~\ref{7.3.5},  $\beta_2(K_{3, 3})=\frac{1}{4}$.   So,
as suggested in \ref{8.9}, we have a complicated proof of the
classical fact that $K_{3, 3}$ is not planar.

Statement $\V(3)$ implies the following generalization of $\II(3)$.

\begin{Proposition}\label{11.4.3} Suppose $A$ is a flag triangulation of a
$2$--sphere with $g+1$ holes.  Let $S_0, \dots , S_g$ be the boundary
components of $A$ and set $$\alpha =\frac{1}{2} (\beta_1(S_0)+\cdots
+\beta_1(S_g)).$$ Then
$$\beta_i(A)=\begin{cases}   \alpha &\text{if $i=1$, }\\  0
&\text{if $i\not= 1$.}
\end{cases}$$
If, in addition, $\partial A$ is a full subcomplex of $A$, then
$$\beta_i(A, \partial A)=\begin{cases}  g+\alpha &\text{if $i=2$, }\\
0 &\text{if $i\not= 2$.}
\end{cases}$$
\end{Proposition}

\begin{proof}  
As in Theorem~\ref{11.3.1}, embed $A$ in the flag triangulation $S$ of
$\BS^2$ obtained by introducing new vertex $s_i$ for each boundary
component $S_i$.  Since $\H_i(S)$ vanishes,
$$\H_i(A)\cong\H_{i+1}(S, A)\cong \H_{i+1}(C_{s_0} S_0\cup\dots\cup
C_{s_g} S_g, S_0\cup\dots\cup S_g).$$ A simple calculation using
Lemma~\ref{7.3.3} and~\ref{7.3.4} gives that this is nonzero only for
$i+1=2$ and that
$$\beta_2(C_{s_0}S_0\cup\dots\cup C_{s_g}S_g, S_0\cup \dots\cup
S_g)=\frac{1}{2}(\beta_1(S_0)+\cdots +\beta_1(S_g)).$$ The first
formula follows.

To prove the second, consider the pair $(S, C_{s_0}S_0\cup\dots\cup
C_{s_g}S_g )$.  By excision, its homology is isomorphic to that of
$(A, \partial A)$.  Hence,
$$\beta_{i+1}(A, \partial A)=\beta_i(C_{s_0}S_0\cup\dots\cup C_{s_g}
S_g)$$ and by Lemma~\ref{7.3.3} and~\ref{7.3.4}, the second term is
nonzero only for $i=1$, in which case,
$$\beta_1(C_{s_0}S_0\cup\dots\cup C_{s_g}S_g)=g+\alpha.$$
\end{proof}

\subsection{Surfaces of higher genus}\label{11.5}  Suppose $L_g$ is a
triangulation of a closed orientable surface of genus $g$ as a flag
complex.  In \cite{Ak}, Akita points out that $\chi^{\orb}
(\Sigma_{L_g}/W_{L_g})=g$.  This, together with the calculation in
Proposition~\ref{11.4.3}, makes the following generalization of
$\I(3)$ a very plausible conjecture.

\begin{Conjecture}\label{11.5.1}$\beta_i(L_g)=0$ for $i\not=2$ and
$\beta_2(L_g)=g$.
\end{Conjecture}

Akita also proves in~\cite{Ak} that if $A$ is a  $1$--dimensional flag
complex (ie,  if it is a simplicial graph without any circuits of
length $3$) and if $A$ embeds in an orientable surface of genus $g$,
then  $\chi^{\orb} (\Sigma_A/W_A) \le g$.  If the above conjecture
holds, then the following analog of Theorem~\ref{11.4.1} gives a
stronger result.

\begin{Proposition}\label{11.5.2} 
Assume Conjecture~\ref{11.5.1}. If a finite flag complex $A$ can be
embedded as a subcomplex of an orientable surface of genus $g$, then
$\beta_2(A)\le g$.
\end{Proposition}

\begin{proof}  
As in the proof of Theorem~\ref{11.4.1}, we can assume that $A$ is a
full subcomplex of some flag triangulation $L$ of the orientable
surface of genus $g$.  By Lemma~\ref{9.2.2}, $\b_3(L, A)=0$; hence,
the map $\H_2(A)\to \H_2(L)$ is injective.  Since we are assuming
$\beta_2(L)=g$, the result follows.
\end{proof}

\subsubsection{Example}\label{11.5.2.1}
We repeat an example from~\cite{Ak}. Let $K_{m, n}$ denote the join of
$m$ points and $n$ points (a complete bipartite graph).
By~\ref{7.3.7},
$$\b_2(K_{m, n})=(\frac{m}{2} -1 )(\frac{n}{2} -1
)=\frac{(m-2)(n-2)}{4}.$$ By~\cite[Theorem 4.5.3]{GT} the minimal
genus of a surface in which $K_{m, n}$ embeds is the least integer
$\ge$ this number.

\section{\label{12}Reflection type covers}

In this section we use the fact that $\l ^2$--Betti numbers are
multiplicative with respect to finite coverings (cf~\ref{3.3.3}).  In
particular, in~\ref{12.3}, we use this to show that Conjecture $\III(2k+1)$ is
implied by Conjecture $\IV (2k+1)$ (that $\beta _{k+1}(D)=0$ whenever
$D$ is a $GHD^{2k}$).

\subsection{Reflection subgroups}\label{12.1}  Let $L$ be a finite flag
complex.  To simplify notation we write $W$, $\Sigma$ and $K$ for
$W_L$, $\Sigma_L$ and $K_L$, respectively.  In this subsection we will
state some basic facts about $W$ and $\Sigma$. Most of the proofs will
be left as exercises for the reader.  (They are all straightforward
adaptations of standard arguments from the theory of classical
reflection groups, for example, as explained in \cite{B}.)

\subsubsection{Reflections, walls, half-spaces}\label{12.1.1}  An element
of $W$ is a {\it reflection} if it is conjugate to a fundamental
generator, ie, to an element of ${\mathcal S}_0(L)$.  Given a
reflection $r$, the fixed set of $r$ on $\Sigma$ is denoted by $\Sigma
(r)$ and called the {\it wall} associated to $r$.  Each wall separates
$\Sigma$ into two pieces, called the {\it half-spaces} bounded by the
wall.  To be more explicit, for each reflection $r$, let $P_r=\{ w\in
W| \l(rw)>\l(w)\}$ (where $\l(~)$ denotes word length) and let $H(r)$
denote the union of the chambers $wK$, with $w\in P_r$.  Then $H(r)$
is the half-space bounded by $\Sigma (r)$ which contains the
fundamental chamber $K$.  The other half-space is $rH(r)$.

\subsubsection{Convexity and half-spaces}\label{12.1.2}  Each half-space
is geodesically convex.  (The proof uses the fact that there is a
distance decreasing retraction $\Sigma\to H(r)$ called the ``folding
map''.)  If $C$ is any convex union of chambers, then it is the
intersection of the half-spaces which contain it.

\subsubsection{Supporting walls}\label{12.1.3}  Suppose $C$ is a convex
union of chambers.  A wall $\Sigma (r)$ is a {\it supporting wall} of
$C$ if (a)  $C$ is contained in one of the half-spaces bounded by
$\Sigma (r)$ and (b)  the intersection $C\cap\Sigma (r)$ is nonempty
and is not contained in any other wall.  Let $\Supp(C)$ denote the set
of reflections $r$ such that $\Sigma (r)$ is a supporting wall of $C$.

\subsubsection{The subgroup generated by $\Supp(C)$}\label{12.1.4}  
For each $r \in \Supp(C)$ denote $C \cap \Sigma (r)$ by $C_r$ and call
it the {\it mirror} of $C$ associated to $r$. Let $G= \langle \Supp(C)
\rangle$ be the subgroup of $W$ generated by $\Supp(C)$.  Next we want
to give a  standard argument which shows that $G$ is a Coxeter group
and that $C$ is a fundamental domain for the $G$--action on $\Sigma$.
Let $\hat G$ be the group defined by the following presentation: there
is a generator $\hat r$ for each $r \in \Supp(C)$ and there are
relations,  $\hat r^2 = 1$, for each $r \in \Supp(C)$ and $(\hat r_1
\hat r_2)^2 = 1$, whenever $C_{r_1} \cap C_{r_2} \neq \emptyset$.
Thus, $\hat G$ is a right-angled Coxeter group.   Let $\theta: \hat G
\to G$ be the epimorphism defined by $\theta (\hat r) = r$.  Let
$\mathcal U(\hat G, C) = (\hat G \times C)/\sim$, where $\sim$ denotes
the equivalence relation generated by $(\hat g, x) \sim (\hat g \hat
r, x)$ whenever $x \in C_r$.  Let  $[\hat g, x]$ denote the image of
$(\hat g, x)$ in $\mathcal U(\hat G, C)$.  The group $G$ acts
naturally on $\mathcal U(\hat G, C)$.  For each $x \in C$, let $G_x$
(resp.  $\hat G_x$) denote the subgroup of $G$ (resp. $\hat G$)
generated by the reflections across the mirrors of $C$ which contain
$x$ and let $U_x$ be an open neighborhood of $x$ in $C$ which
intersects only those mirrors which contain $x$.  Let $\mathcal U(\hat
G_x, U_x)$ denote the image of $\hat G_x \times U_x$ in $\mathcal
U(\hat G, C)$.  Then $G_xU_x$ is an open neighborhood of $x$ in
$\Sigma$, $\mathcal U(\hat G_x, U_x)$ is an open neighborhood of
$[1,x]$ in $\mathcal U(\hat G, C)$ and both $G_x$ and $\hat G_x$ are
isomorphic to $(\BZ_2)^m$ where $m$ is the number of mirrors containing
$x$.  Let $f: \mathcal U(\hat G, C) \to \Sigma$ denotes the
$\theta$--equivariant map defined by $f([\hat g, x]) = \theta (\hat
g)x$. Using the fact that $W$ is right-angled,  it can be seen that
$U_x$ is a fundamental domain for  the $G_x$--action on $G_xU_x$.  It
follows from this that  $f$ maps $\mathcal U(\hat G_x, U_x)$
homeomorphically onto $G_xU_x$ and  consequently, that $f$ is a
covering projection.  Since $\Sigma$ is simply connected, this implies
that $f$ is a homeomorphism and that $\theta$ is an isomorphism.
Thus, $G$ is a right-angled Coxeter group, $C$ is a fundamental domain
and  $\Supp(C)$ is a fundamental set of generators.

The nerve of $\langle\Supp(C)\rangle$ (cf~\ref{5.1}) is the flag
complex $L(C)$ which can be defined as follows.  The vertex set of
$L(C)$ is $\Supp(C)$ and two distinct vertices $r_1$ and $r_2$ span an
edge if and only if $(r_1r_2)^2=1$ (a flag complex is determined by
its $1$-skeleton).  Thus, $\langle \Supp(C)\rangle\cong W_{L(C)}$.

\subsubsection{}\label{12.1.5}  
Suppose $W_A$ is a special subgroup of $W$.  Then $W_A K$ is a convex
union of chambers.  The corresponding subgroup $W_{L(W_A K)}$ can be
identified with the kernel of the homomorphism $\varphi_A\colon W\to
W_A$, defined by specifying its values on the generating set
${\mathcal S}_0(L)$ as follows:
$$\varphi_A (s)=\begin{cases} s &\text{if $s\in{\mathcal S}_0(A)$, }\\
1 &\text{if $ s\not\in{\mathcal S}_0(A)$.}\end{cases}$$ We note that
this kernel is of finite index in $W$ if and only if $A$ is a simplex.

\subsubsection{Doubling along a vertex}\label{12.1.6}  Suppose $\sigma$
is a simplex of $L$.  Denote the corresponding flag complex
$L(W_\sigma K)$ by $d_\sigma L$.  Thus, $W_{d_\sigma L}$ is a normal
subgroup of index $2^{\dim\sigma +1}$ in $W$.  The special case where
$\sigma$ is a vertex $v$ will be denoted $d_vL$ and called the {\it
double of $L$ along $v$}.

\subsubsection{Description of $d_vL$}\label{12.1.7}  For each vertex $s$ of
$L-v$, we get two supporting walls of $W_vK$, namely, $\Sigma (s)$ and
$v\Sigma (s)=\Sigma (vsv^{-1})$.  When $s\in L_v$, $vsv^{-1}=s$ and
these two walls coincide.  Hence, $d_vL$ is formed by taking two
copies of $L-v$ and gluing them together along the subcomplex $L_v$.

\subsubsection{Iterated doubles}\label{12.1.8}  Suppose $s_1$ and $s_2$
are two vertices of $L$ which are not connected by an edge.  For each
positive integer $N$, let $F_N$ be the set of the first $N$ elements
in the list $1, s_1, s_1s_2, s_1s_2s_1, s_1s_2s_1s_2, \dots$  ($F_N$
is a subset of the infinite dihedral group generated by $s_1$ and
$s_2$.)  Then $F_NK$ is a convex union of chambers.  The corresponding
flag complex $d^NL$ is the $N$-fold {\it iterated double} of $L$ along
$(s_1, s_2)$. By  \ref{12.1.4}, $F_NK$ is a fundamental domain for the
$W_{d^NL}$-action; hence, the subgroup $W_{d^NL}$ is of index $N$ in
$W$.

\subsubsection{}\label{12.1.9}
Suppose that $L=S$, a $GHS^{n-1}$ and that $C$ is a convex union of a
finite number of chambers in $\Sigma_S$.  Then $C$ is contractible
(since it is $\text{CAT}(0)$) and hence, a generalized homology
$n$-disk.  It follows that the flag complex $L(C)$ (which is ``dual''
to the boundary of $C$) is also a $GHS^{n-1}$.

\subsection{Inequalities}\label{12.2}  In this subsection we return to
the situation of statement $\III(2k+1)$ in~\ref{8.3}:  $(D, S_v)$ is a
$GHD^{2k}$ and $S=D\cup CS_v$ is the generalized homology $2k$-sphere
obtained by adjoining the cone on the boundary.  (The cone point is
$v$.)  Set
$$\alpha_{k+1} = \dim_{W_S}(\overline\im
(i_\ast\colon\H_{k+1}(D)\to\H_{k+1} (S))) .$$ By excision and
Lemma~\ref{7.3.3} (1), $\beta_{k+1}(S, D)=\beta_{k+1} (CS_v, S_v)=
\frac{1}{2}\beta_k(S_v)$.  Hence, the sequence of the pair $(S, D)$
gives the following inequality.

\subsubsection{}\label{12.2.1}
\myequation{0 \le \beta_{k+1}(S)-\alpha_{k+1} \le
\frac{1}{2}\beta_k(S_v).}

Next suppose $\beta_{k+1}(S_v)=0$.  (For example, this holds if
$\I(2k)$ holds for the link $S_v$.)  Then since $\beta_{k+2}(S,
D)=\frac{1}{2}\beta_{k+1} (S_v)=0$, the map
$i_\ast\colon\H_{k+1}(D)\to\H_{k+1}(S)$ is injective and
$\alpha_{k+1}=\beta_{k+1}(D)$.  So, ~\ref{12.2.1}, can be rewritten as:

\subsubsection{}\label{12.2.2}
\myequation{0 \le \beta_{k+1}(S)-\beta_{k+1}(D) \le
\frac{1}{2}\beta_k(S_v) .}

The next lemma shows that this inequality can be improved by a factor
of $2$.

\begin{Lemma}\label{12.2.3}Suppose, as above, that $(D, S_v)$ is a $GHD^{2k}$
and that \linebreak  $\beta_{k+1}(S_v)=0$.  Then
$$\beta_{k+1}(S)-\beta_{k+1}(D)\le \frac{1}{4}\beta_k(S_v) .$$
\end{Lemma}

\begin{proof}  
By~\ref{12.1.9}, the double of $S$, $d_vS$, is also a $GHS^{2k}$.
By~\ref{12.1.7}, $d_vS$ is the union of two copies of $D$ glued along
$S_v$.  So, we have a Mayer-Vietoris sequence,
$$
0 \to \H_{k+1}(D)\oplus\H_{k+1}(D) \to \H_{k+1}(d_vS)
\xrightarrow{\partial}  \H_k(S_v) \to \H_k(D)\oplus\H_k(D).
$$
By Lemma~\ref{7.4.4}, the kernel of the map $\H_k(S_v)\to\H_k(D)$ into
either factor has dimension $\frac{1}{2}\beta_k(S_v)$.  Thus, the
kernel of the map $\H_k(S_v)\to\H_k(D)\oplus\H_k(D)$ has dimension
$\le \frac{1}{2}\beta_k(S_v)$.  Hence,
$$\beta_{k+1}(d_vS)\le 2\beta_{k+1}(D)+\frac{1}{2}\beta_k(S_v) .$$
Substituting in $2\beta_{k+1}(S)$ for $\beta_{k+1}(d_vS)$
(by~\ref{12.1.6} and~\ref{3.3.3}), we get the desired inequality.
\end{proof}

\subsection{$\IV(2k+1) \implies \III(2k+1)$}\label{12.3}  Suppose that
$\IV(2k+1)$ is true (ie, that $\beta_{k+1}(D')=0$ for any
generalized homology $2k$-disk $D'$).  Then the
inequality~\ref{12.2.1} becomes,

\subsubsection{}\label{12.3.1}
\myequation{\beta_{k+1}(S) \le \frac{1}{2}\beta_k(S_v).}

As we shall see below, this inequality forces $\beta_{k+1}(S)=0$.
Since $\H_{k+1}(S)$ is the previous term for the map in the
Mayer-Vietoris sequence which is under consideration in $\III(2k+1)$,
the next lemma shows that $\IV(2k+1)$ implies $\III(2k+1)$.

\begin{Lemma}\label{12.3.2} As in $\III(2k+1)$, let $(D, S_v)$ be a $GHD^{2k}$
and let $S=D\cup CS_v$ be the $GHS^{2k}$ formed by adjoining a cone on
the boundary.  Assume $\IV(2k+1)$ is true.  Then $$\beta_{k+1}(S)=0.$$
\end{Lemma}

\begin{proof}  {\bf Case 1}\qua  Suppose $D-S_v$ is not a simplex.  Then we can
find vertices $s_1, s_2$ in $D-S_v$ which are not connected by an
edge.  Let $\widetilde{S}$ be the $N$--fold iterated double $d^NS$
along $(s_1, s_2)$, as defined in~\ref{12.1.8}.  Then $v$ has $N$
preimages in $\widetilde{S}$ and the link of each is isomorphic to
$S_v$.  Choose one, say $v_1$, and set
$\widetilde{D}=\widetilde{S}-v_1$.  Since we are assuming $\IV(2k+1)$,
we have, by~\ref{12.3.1}, that $\beta_{k+1}(\widetilde{S})\le
\frac{1}{2}\beta_k(S_{v_1})$.  By~\ref{3.3.3}, $\beta_{k+1}
(\widetilde{S}) = N\beta_{k+1}(S)$.  Hence,
$$\beta_{k+1}(S)\le \frac{1}{2N}\beta_k(S_v) .$$ Since this holds for
any $N$, $\beta_{k+1}(S)=0$.

\medskip

{\bf Case 2}\qua  $D-S_v$ is a simplex $\sigma$.  If $\dim\sigma =0$,
then $S$ is a suspension and we are done by Lemma~\ref{7.3.6}.  If
$\dim\sigma >0$, then let $S'=d_\sigma S$ (defined in~\ref{12.1.6}).
By~\ref{3.3.3}, $\beta_{k+1}(S')= m\beta_{k+1}(S)$ where
$m=2^{\dim\sigma +1}$.  Moreover, there are $m$ preimages of $v$ in
$S'$, no two of which are connected by an edge and such that the link
of each is isomorphic to $S_v$.  Choose one of these preimages, say
$v_1$, and set $D'=S'-v_1$.  Since $D'$ contains $m-1$ preimages of
$v$, $m-1\ge 2$, we can apply Case 1 to $(D', S_v)$ to conclude that
$0=\beta_{k+1} (S')=m\beta_{k+1}(S)$.
\end{proof}

\subsection{Atiyah's Conjecture}\label{12.4} In \cite{A} Atiyah 
conjectured that $\ell^2$--Betti numbers of any geometric $G$--complex
$X$ are rational numbers. A refinement of this states that if $m$
denotes  the least common multiple of the orders of the finite
subgroups of $G$, then  $mb^{(2)}_i(X, G)$ is an integer. An
equivalent form (see~\cite{E}) of  this conjecture is  the following.

\begin{Conjecture}[Atiyah]\label{12.4.1} Suppose 
$\phi\colon (\BZ G)^p \to (\BZ G)^q$ is a homomorphism of free  $\BZ
G$--modules and that $\Hat{\phi}\colon {\l^2(G)}^p \to {\l^2(G)}^q$  is
the induced map of Hilbert $G$--modules.  As above, let  $m$ denote
the least common multiple of the orders of the finite subgroups of
$G$,    and suppose that $m$ is finite. Then
$$m\dim_G(\Ker\Hat{\phi})\in \BN.$$
\end{Conjecture}
\subsubsection{}\label{12.4.2}
The above conjecture implies that in the (weakly) exact sequence of
any  pair of geometric $G$--complexes, the von Neumann dimension of the
kernel or  image of any map  is a nonnegative rational number with
denominator  dividing $m$.

\subsubsection{}\label{12.4.3} If $W_L$ is a right-angled Coxeter group, then 
$m=2^{\dim L +1}$.

\subsubsection{}\label{12.4.4} Taken together with Atiyah's Conjecture, 
Lemma~\ref{12.2.3} provides some convincing evidence for the truth of
$\III(2k+1)$. Let $S$, $S_v$ and $D$ be as above and assume $\I(2k)$
holds for $S_v$.  Then the largest possible denominator for
$\b_k(S_v)$  ($=(-1)^k \chi^{(2)}(S_v)$) is $2^{2k}$. If
Conjecture~\ref{12.4.1}  holds for $W_S$, then the largest  possible
denominator for $\b_{k+1}(S)-\b_{k+1}(D)$ is $2^{2k+1}$.  So, if
$\b_k(S_v)$ has the smallest possible nonzero value, namely
$(1/2)^{2k}$,  and if Conjecture~\ref{12.4.1} is true, then
Lemma~\ref{12.2.3} implies that $\b_{k+1}(S)=\b_{k+1}(D)$.  This
implies that in the Mayer--Vietoris sequence for $S=D\cup CS_v$,  the
map $\H_{k+1}(D)\oplus\H_{k+1}(CS_v) \to \H_{k+1}(S)$ is surjective.
and hence, that the map $\H_{k}(S_v) \to \H_{k}(D)\oplus\H_{k}(CS_v)$
is injective, ie, that $\III(2k+1)$ holds for the pair $(S, S_v)$.

\section{Inclusions of walls}\label{13}

Let $S$ be a flag triangulation of a generalized homology sphere of
dimension $2k$. Suppose that $s$ is a vertex of $S$ and that $S_s$
denote its link in $S$. By~\ref{8.3.2} and Theorem~\ref{9.3.1}, our
conjecture has been reduced to $\III'(2k+1)$, which asserts that the
map $i_*\colon \H_k(S_s) \to \H_k(S)$, induced by inclusion, is zero.

In this section we shall make a series of observations about this
problem. Our eventual point is made in~\ref{13.3}:  Conjecture  $\III'(2k+1)$ is
essentially equivalent to a certain estimate on the rate of growth of
the norms of  $k$--dimensional homology classes in the hypersurface
$\Sigma _{S_s}$ as they are pushed onto an ``equidistant
hypersurface''.

To simplify notation set $\Sigma=\Sigma_S$ and $W=W_S$. We recall
(from \ref{7.1}) that $\H_k(S_s)$ stands for $\h_k(W\Sigma_{S_s}) $,
where $\Sigma_{S_s}$ is the special subcomplex corresponding to
$S_s$. We also note that $\Sigma_{CS_s}$ can be identified with
$\Sigma_{S_s} \times [-1, 1]$, where the wall $\Sigma(s)$ (defined
in~\ref{12.1.1}) corresponds to $\Sigma_{S_s}\times 0$. In particular,
$\Sigma_{S_s}$ is $W_{S_s}$--equivariantly homeomorphic to the wall
$\Sigma(s)$.

\subsection{Reduction to a single wall}\label{13.1} 
Since $W\Sigma _{S_s}$ is the disjoint union of copies of $\Sigma(s)$,
one for each coset of $W_{S_s}$ in $W$, the Hilbert space
$\h_k(W\Sigma _{S_s})$ is an orthogonal sum of copies of
$\h_k(\Sigma(s))$. Hence, to prove that $i_*\colon \h_k(W\Sigma
_{S_s})$\break $ \to \h_k(\Sigma)$  is the zero map, it is necessary and
sufficient to show that its restriction to one summand,
$\h_k(\Sigma(s))$, is zero.

\subsection{The map into unreduced homology}\label{13.2} 
The map $i_*\colon \h_k(\Sigma(s))\to \h_k(\Sigma)$ factors as a
composition $p\circ \widehat{i_*}$, where the map
$\widehat{i_*}\colon \h_k(\Sigma(s)) \to H^{(2)}_k(\Sigma)$ is
induced by the  inclusion of the harmonic $k$--cycles into the
$\ell^2$--cycles on $\Sigma$ and where the map $p\colon
H^{(2)}_k(\Sigma) \to \h_k(\Sigma)$ is projection onto the harmonic
cycles. (Recall, from 2.3.3, that $H^{(2)}_k(~)$ denotes unreduced
$\ell^2$--homology.)

\begin{Lemma}\label{13.2.1} If 
$i_*\colon \h_k(\Sigma(s)) \to \h_k(\Sigma)$ is the zero map, then the
map $\widehat{i_*}\colon \h_k(\Sigma(s)) \to H^{(2)}_k(\Sigma)$ is
injective.
\end{Lemma}

\begin{proof}
Suppose that $x$ is a  harmonic $k$--cycle in $\h_k(\Sigma(s))$ such
that $\widehat{i_*}({x})=0$ in $H^{(2)}_k(\Sigma)$. In other words,
$x=d(y)$ for some $(k+1)$--chain $y$ in $C_{k+1}(\Sigma)$.  ( We
identify $x$ with its image under the inclusion of chains  $
C_k(\Sigma(s)) \hookrightarrow C_k(\Sigma)$.)  The wall $\Sigma(s)$
divides $\Sigma$ into two half-spaces; let us call them $\Sigma_+$ and
$\Sigma_-$.  We first claim that we can find a $(k+1)$--chain $y' \in
C_{k+1}(\Sigma)$ so that $x=d(y')$ and so that $y'$ is supported on
only one half-subspace, say $\Sigma_+$. To see this, first  write
$y=y_+ + y_-$, where  $y_+$ (respectively, $y_-$) is supported on
$\Sigma_+$ (respectively, $\Sigma_-$) and therefore, $d(y_+)$ and
$d(y_-)$ are both supported on $\Sigma(s)$. Then set  $y'=y_+ +
sy_-$. Since $s$ fixes $\Sigma(s)$,
$$d(y')=d(y_+)+sd(y_-)=d(y_+ + y_-)=d(y)=x.$$ Set $z=y' - sy'$. Then
$$d(z)=d(y')-sd(y')=0, $$ so $z$ is a $(k+1)$--cycle in
$C_{k+1}(\Sigma)$. Let $\overline{z}$ denote its image in reduced
$\ell^2$--homology $\h_{k+1}(\Sigma)$.

Consider  the Mayer--Vietoris sequence of $\Sigma=\Sigma_+\cup
\Sigma_-$ in unreduced  $\ell^2$--homology. Let $\partial \colon
H^{(2)}_{k+1}(\Sigma) \to H^{(2)}_k(\Sigma(s))$ be the connecting
homomorphism and let  $\partial_* \colon \h_{k+1}(\Sigma) \to
\h_k(\Sigma(s))$ be the induced map of quotients. It follows from the
definition of $\partial$, that $\partial([z])=[x]$ in unreduced
homology and therefore,  $\partial_*(\overline{z})=x$, since $x$ is
harmonic.  On the other hand, just as in \ref{7.4.5}, the map
$i_*\colon \h_k(\Sigma(s)) \to \h_k(\Sigma)$ is isomorphic, under
Poincar\'{e} duality to $\partial^* \colon \h_{k}(\Sigma(s)) \to
\h_{k+1}(\Sigma)$. Hence, if $i_*$ is the zero map, then so is
$\partial_*$. Therefore, our hypothesis implies that $x=
\partial_*(\overline{z})=0$ and consequently, that $\widehat{i_*}$ is
injective.
\end{proof}

\subsubsection{}\label{13.2.2} An alternative proof of this lemma can be
constructed as follows. If  $i_*\colon \h_k(\Sigma(s)) \to
\h_k(\Sigma)$ is zero, then, by Theorem~\ref{9.3.1},
$\h_{k+1}(\Sigma)=0$. As before, suppose $\widehat{i_*}(x)=0$ and
define $z$ as before. If $x\not =0$, then we can find a $u\in
\h_k(\Sigma(s))$ with nonzero intersection number with $x$. It is
geometrically clear that the intersection number of $u$ and
$\overline{z}$ in $\Sigma$ is the same. Hence, $u\not=0$,
contradicting $\h_k(\Sigma(s))=0$.

\subsection{Equidistance hypersurfaces}\label{13.3}
There is an infinite sequence of disjoint $W_{S_s}$--stable
hypersurfaces  $\Sigma_1,  \Sigma_2, \dots, \Sigma_n, \dots$ in
$\Sigma$ and  $W_{S_s}$--equivariant homotopy equivalences $p_n\colon
\Sigma_n \to \Sigma(s)$.

To define these, let $X_n$ denote the union of all cells in the
half-space $\Sigma_+$  of combinatorial distance $\le n$ from
$\Sigma(s)$. (This definition is intentionally vague; there are
several possible definitions of ``combinatorial distance'' and at
least two possible cell structures on $\Sigma$ --- one is the
cellulation by cubes and the other is the dual cellulation by
chambers.) The boundary of $X_n$ has two components, one is
$\Sigma(s)$, the other is denoted by $\Sigma_n$.

There is a $W_{S_s}$--equivariant retraction of $X_n$ onto
$\Sigma(s)$. Its restriction to $\Sigma_n$ is $p_n$.  Since $X_n$ is
$W_{S_s}$--cocompact and $W_{S_s}$--homotopy equivalent to $\Sigma(s)$
we have that $i_*\circ {p_n}_*={i_n}_*$ where $i_n\colon \Sigma_n \to
\Sigma$ denotes the inclusion.

\subsubsection{}\label{13.3.1} Given a cycle $x\in C_k(\Sigma(s))$ we can
find  a cycle $x_n\in C_k(\Sigma_n)$ with  ${p_n}_*(x_n)=x$ in reduced
homology.  So, $x_n$ will be homologous to $x$ in $\Sigma$. Therefore,
any linear combination $a_1x_1+\dots+a_nx_n $ with $a_1+\dots+a_n=1$
is also homologous to $x$  in $\Sigma$. Let $y_n$ denote a linear
combination of the above form which has the minimal norm.  Since the
cycles $x_i$ are supported on disjoint sets, they are mutually
orthogonal in  $C_k(\Sigma)$. Then an  easy inductive argument shows
that the norm of $y_n$ is given by
$$\frac{1}{\norm{y_n}^2}=\frac{1}{\norm{x_1}^2}+\dots+\frac{1}{\norm{x_n}^2}.$$
Hence, if the series $\sum_{n=1}^\infty \frac{1}{\norm{x_n}^2}$ is
divergent, then $\lim_{n\to \infty} \norm{y_n}=0$ and therefore,
$i_*(x)=0$ (since there would be a sequence of cycles representing
$i_*(x)$, with norms going to $0$).  For example,  this argument works
if $\norm{x_n}^2$ grows sublinearly in $n$.

\subsubsection{}\label{13.3.2} The cell structure on $\Sigma_n$ can be
obtained from that of $\Sigma_{n-1}$ by a subdivision process which
can be described by a regular procedure which depends only on the
initial data. It follows that $x_{n-1}$ can be pushed to $\Sigma_n$ by
a process which replaces each $k$--cell in $x_{n-1}$ by a $k$--chain in
$\Sigma_n$ and that the maximum norm of this chain can be bounded
above by a constant $D$ which is independent of $n$. This gives
$\norm{x_n}\le D\norm{x_{n-1}}$ and hence, that $\norm{x_n}\le
D^n\norm{x}$, an estimate that is  much worse than what we want. On
the other hand, there are many possible choices for the ``pushing
procedure'' of associating a $k$--chain to each $k$--cell. Roughly, the
hope is that one can show that there are at least $D$ such choices of
disjoint $k$--chains. We could then choose $x_n$ to be the average of
$D$ disjoint pushes of $x_{n-1}$, obtaining  $\norm{x_n}\le \norm{x}$,
the desired result.

\subsection{Equidistance hypersurfaces in hyperbolic
space}\label{13.4} In the case of hyperbolic $(2k+1)$--space
$\BH^{2k+1}$, the above argument can be made precise.  Let $\BH^{2k}$
be a totally geodesic hyperplane in $\BH^{2k+1}$.  We claim  that the
map $\h_k(\BH^{2k}) \to \h_k(\BH^{2k+1})$, induced by inclusion, is
the zero map. This is of course true, since $\h_*(\BH^{2k+1})=0$ by
\cite{Do2}, but our proof below  does not depend on that, and, in
fact, can be used to give an alternative proof of  Singer's Conjecture
for hyperbolic space.  Our argument  uses $L^2$--de~Rham cohomology
theory and is dual to the argument in~\ref{13.3.1}--\ref{13.3.2}. We
will show that the map $\h_k(\BH^{2k+1}) \to \h_k(\BH^{2k})$, induced
by restriction of forms, is the zero map.
 
Let $N_t$ be the hypersurface in $\BH^{2k+1}$ consisting of the points
of (oriented) distance $t$ from $\BH^{2k}$. Let $p_t\colon N_t\to
\BH^{2k}$ be the projection which takes a point in $N_t$ to the
closest point in $\BH^{2k}$. Then $p_t$ is a homothety. Let
$\phi_t\colon \BH^{2k} \to N_t$ be its inverse. Also, let  $i\colon
\BH^{2k} \to \BH^{2k+1}$ and  $i_t\colon N_t \to \BH^{2k+1}$ be the
inclusions. Thus, $i$ and $i_t\circ\phi_t$ are properly homotopic.

Let $\omega$ be a closed $L^2$-$k$--form on $\BH^{2k+1}$.  We claim
that the restriction $i^*(\omega)$ of $\omega$ to $\BH^{2k}$
represents the zero class in reduced  $L^2$--cohomology. Indeed,
suppose  $[i^*(\omega)]\neq 0$. Then   $\norm{i^*(\omega)}\geq
\norm{[i^*(\omega)]}\geq 0$, where  $\norm{[i^*(\omega)]}$ denotes the
norm of the harmonic representative of the cohomology class
$[i^*(\omega)]$.  Since $\phi_t$ is a conformal diffeomorphism, it
follows  that it preserves norms of middle-dimensional forms:
$\norm{\phi^*_t(i^*_t(\omega)}=\norm{i^*_t(\omega)}$.  Since $i$ and
$i_t\circ\phi_t$ are properly homotopic, $[\phi^*_t(i^*_t(\omega)]=
[i^*(\omega)]$, so it follows that $\norm{i^*_t(\omega)}\ge
\norm{[i^*(\omega)]}$.  Now, since $i^*_t(\omega)$ is just a
restriction of $\omega$, we have a pointwise inequality
$\norm{\omega}_x \ge \norm{i^*_t(\omega)}_x$. Therefore, using
Fubini's Theorem, we obtain
\begin{multline*}
\norm{\omega}^2=\int_{\BH^{2k+1}} \norm{\omega}^2_x\, dV= \int_\BR
\int_{N_t} \norm{\omega}^2_x \, dA \, dt \ge \\  \int_\BR \int_{N_t}
\norm{i^*_t(\omega)}^2_x \, dA \, dt  =\int_\BR
\norm{i^*_t(\omega)}^2 \, dt \ge  \int_\BR   \norm{[i^*(\omega)]}^2 \,
dt= \infty.
\end{multline*}
This contradicts our assumption that $\omega$ is $L_2$--form and
completes the proof.

\section{Virtual fibrations over the circle}\label{14}

In this section we discuss some ideas for another possible attack on
the conjecture that $\I(2k+1)$ is true.

\subsection{Another conjecture}\label{14.1} 
As in \ref{1.1} suppose that $X$ is a simply connected geometric
$G$--complex. The orbihedron $X/G$ {\em virtually fibers over $\BS^1$}
if there is a subgroup $\Gamma$ of finite index in $G$ such that
$\Gamma$ acts freely on X and such that $X/\Gamma$ fibers over $\BS^1$.

\begin{Theorem}[L{\"u}ck \cite{L1}]\label{14.1.1} 
If a finite complex fibers over $\BS^1$, then the reduced
$\l^2$--homology of its  universal cover vanishes in all dimensions.
\end{Theorem}

\subsubsection{}\label{14.1.2}
It follows that if $X$ is simply connected and if $X/G$ virtually
fibers over $\BS^1$, then $\h_i(X)$=0 for all $i$.

\subsubsection{}\label{14.1.3}
There is an obvious obstruction for  $X/G$ to virtually fiber over
$\BS^1$: its orbihedral Euler characteristic must vanish. We note,
however, that if $X/G$ is a closed odd-dimensional orbifold, then this
obstruction always vanishes.

In the late 1970's Thurston asked if the following conjecture were
true.

\begin{Conjecture}[Thurston]\label{14.1.4}
Let $M^3$ be a closed irreducible $3$--mani\-fold (or developable
orbifold) with  infinite fundamental group. Then $M^3$  virtually
fibers over $\BS^1$.
\end{Conjecture}

From now on we suppose that $S$ is a smooth triangulation of
$\BS^{n-1}$ as a flag complex. As usual, to simplify notation, we set
$W=W_S$, $K=K_S$ and $\Sigma=\Sigma_S$. Then $\Sigma/W$ can be given
the structure of a smooth $n$--dimensional orbifold. In the following
conjecture we shall also assume that $n$ is odd.

\begin{Conjecture}\label{14.1.5}
Suppose $S$ is a smooth triangulation of $\BS^{2k}$ as a flag
complex. Then  $\Sigma/W$ virtually fibers over $\BS^1$.
\end{Conjecture}

\subsubsection{}\label{14.1.6}
By L{\"u}ck's Theorem, this conjecture implies $\I(2k+1)$ (at least in
the case where $S$ is a sphere rather than just a generalized homology
sphere).

There are reasons for believing that, in odd dimensions, the orbifolds
$\Sigma/W$ should be viewed as being analogous to $3$--manifolds. One
such reason is the following. Suppose that $M^n$ is a manifold
covering of  $\Sigma/W$ corresponding to a normal torsion-free
subgroup $\Gamma$ of finite index in $W$.  Let $\Sigma (r)$ be a wall
of $\Sigma$ and let $N^{n-1}$ denote its image in $M^n$.  Since
$\Sigma(r)$ is the fixed point set of an isometric reflection on
$\Sigma$ it is a geodesically convex subset.   A well-known argument
of \cite{Mil} then shows that $N^{n-1}$ is a totally geodesic
hypersurface in $M^n$.  (This argument goes as follows.  Let $H$
denote the centralizer of $r$ in $\Gamma$.  Suppose that for some
$\gamma \in  \Gamma$ there is a point $x$ in $\Sigma (r) \cap \gamma
\Sigma (r)$.  The element $r\gamma r\gamma ^{-1}$ fixes $x$.  Since
$\Gamma$ is normal, $r\gamma r\gamma ^{-1} \in \Gamma$ and since
$\Gamma$ is torsion-free, $r\gamma r\gamma ^{-1} = 1$.  Consequently,
$\gamma \in H$.  Therefore, $\Sigma (r)/H =N$ and $\Sigma (r) \to N$
is a covering projection.)  In particular, $N$ is aspherical, $\pi
_1(N)=H$ and the induced homomorphism $\pi_1(N)\to \pi_1(M)$ is
injective. The fact that $M$ has many such ``incompressible
hypersurfaces'' $N$ means that $M$ is a higher-dimensional analog of a
Haken 3--manifold.

We will discuss below two ideas for attacking
Conjecture~\ref{14.1.5}. The first idea  is to find a nowhere-zero,
closed $1$--form on M. It is discussed in subsections~\ref{14.2}
and~\ref{14.3}. The second idea is to find an incompressible
hypersurface $F$ in $M$ to serve as the fiber of a fibration over
$\BS^1$. It is described in subsections~\ref{14.4} to~\ref{14.7}.

\subsection{Nowhere-zero closed $1$--forms}\label{14.2}
It is well-known that a smooth closed manifold $M^n$ fibers over
$\BS^1$ if and only if it admits a nowhere-zero closed
$1$--form. Indeed, if $p\colon M^n\to \BS^1$ is a smooth submersion,
then $p^*(d\theta)$ is such a $1$--form.  Conversely, suppose $\omega$
is a  nowhere-zero closed $1$--form on $M^n$. After adding a closed
$1$--form (of small pointwise norm) to $\omega$ we may assume that, in
addition to being nowhere-zero, $\omega$ has rational periods. In
other words, we may assume that its cohomology class $[\omega]$
actually lies in $H^1(M^n;\BQ)$. Then after replacing $\omega$ by a
suitable multiple, we may assume that it has integral periods, ie,
that $[\omega]$  lies in the image of $H^1(M^n;\BZ)$ in
$H^1(M^n;\BR)$. The cohomology class $[\omega]$ then defines a
homomorphism $\phi_{[\omega]}\colon \pi_1(M^n) \to \BZ$. Finally,
after choosing a basepoint, integration of $\omega$ along paths yields
a submersion  $p_{\omega}\colon M^n\to \BR/\BZ=\BS^1$, with
$p_{\omega}^*(d\theta)=\omega$.

\subsubsection{Remark}\label{14.2.1} 
In~\cite{Fa} Farber gives a direct argument for showing that if $M^n$ 
admits a nowhere-zero closed $1$--form, then the reduced $\l^2$--homology 
of its universal cover vanishes. This gives another proof of L{\"u}ck's 
Theorem in the smooth case.

\subsubsection{Example}\label{14.2.3}
Suppose $S$ is the boundary of an $n$--dimensional octahedron. Then $K$
is an $n$--cube, $W=(D_\infty)^n$ and $\Sigma=\BR^n$.  The commutator
cover of $\Sigma/W$ is an $n$--torus $T^n$, which, of course, fibers
over $\BS^1$. (Any nonzero element of  $H^1(T^n;\BR)$ can be
represented by a linear $1$--form, which is never zero.) So, in this
case, $\Sigma/W$ virtually fibers over $\BS^1$.

\subsection{Double branched covers}\label{14.3}
Suppose $\pi\colon M^n_1 \to  M^n_2$ is a smooth bran\-ched covering, 
branched along a codimension--$2$ submanifold $B^{n-2}\subset
M^n_2$. At any point  $x\not \in \pi^{-1}(B)$, the differential
$d\pi_x \colon T_x(M^n_1) \to T_x(M^n_2)$ is  an isomorphism. On the
other hand,  for  $x \in \pi^{-1}(B)$,   $d\pi_x$ maps
$T_x(\pi^{-1}(B))$ isomorphically  onto $T_x(B)$ and it maps  a
complementary $2$--plane to $0$.

Now suppose that $M^n_2$ admits a nowhere-zero closed $1$--form
$\omega$ such that the restriction of $\omega$ to the tangent bundle of $B$
 is also nowhere-zero. Then
$\pi^*(\omega)$ is a nowhere-zero closed $1$--form on $M^n_1$. Another
way to say essentially the same thing is the following. If $p \colon
M^n_2 \to \BS^1$ is a smooth fibration and if the branch set $B$ is
never tangent to the fibers, then $p\circ\pi \colon M^n_1 \to \BS^1$
is also a smooth fibration.

\subsubsection{An example of Thurston}\label{14.3.1}
The following example of Thurston is explained in~\cite{Su}.  Let $S$
and $S'$ be the triangulations of $\BS^2$ as a boundary of an
octahedron and an icosahedron, respectively.  Then $K_S$ is a cube and
$K_{S'}$ is a dodecahedron. Drawing the  dodecahedron as below, shows
that there is a
map of orbifolds from $K_{S'}$ to $K_S$.

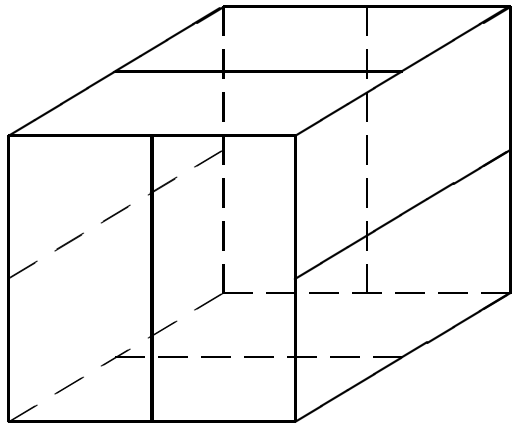
\setlength{\unitlength}{0.25in}
\begin{picture}(16, 11)(-5, -1)
\thicklines \put(0, 0){\line(1, 0){6}}  \put(0, 0){\line(0, 1){6}}
\put(6, 0){\line(0, 1){6}}  \put(0, 6){\line(1, 0){6}}

\put(6, 0){\line(5, 3){4.5}}  \put(6, 6){\line(5, 3){4.5}}
\put(0, 6){\line(5, 3){4.5}} \put(4.5, 8.7){\line(1, 0){6}}
\put(10.5, 2.7){\line(0, 1){6}}

\put(3, 0){\line(0, 1){6}} \put(2.25, 7.35){\line(1, 0){6}}
\put(6, 3){\line(5, 3){4.5}} \thinlines
\multiput(0, 0)(.975, .585){5}{\line(5, 3){.6}}
\multiput(4.5, 2.7)(0, .9){7}{\line(0, 1){.6}}
\multiput(4.5, 2.7)(.9, 0){7}{\line(1, 0){.6}}
\multiput(0, 3)(.975, .585){5}{\line(5, 3){.6}}
\multiput(7.5, 2.7)(0, .9){7}{\line(0, 1){.6}}
\multiput(2.25, 1.35)(.9, 0){7}{\line(1, 0){.6}}
\end{picture}

Take an $8$--fold cover of $K_S$ to get the $3$--torus $T^3$. The
induced covering of $K_{S'}$ is still topologically a torus, but
orbifold singularities remain along six circles (which are the inverse
images of the edges which were introduced in the subdivision of the
cube above). These circles are parallel to the coordinate circles
$T^3$. Let $\pi\colon M^3 \to T^3$ be the $2$--fold branched covering
of $T^3$, branched along these circles. (By Andreev's Theorem $M^3$ 
can be given the structure of the hyperbolic $3$--manifold.) Now
choose a linear fibration $p\colon T^3\to \BS^1$ such that the
coordinate circles are transverse to the fibers. Then  $p\circ\pi
\colon M^3 \to \BS^1$ is a fibration. In other words, 
$\Sigma_{S'}/W_{S'}$ virtually fibers over $\BS^1$.

\subsubsection{}\label{14.3.2}
In the above example we started by subdividing each codimension-one
face of  of a cube. The process of dividing a  codimension-one face of
$K$ in two is dual to the ``edge subdivision process'' of Section 5.3
in~\cite{CD}, applied to an appropriate edge of $S$. So, the method of
the above example shows that we can construct many further examples of
right-angled Coxeter orbifolds (in any dimension $n$) which virtually
fiber over $\BS^1$ by using the following three steps.
\begin{enumerate}
\item Start with a triangulation $S$ so that some finite cover $M^n$
of $\Sigma/W$ fibers over $\BS^1$ (or admits a nowhere-zero closed
$1$--form $\omega$).

\item Subdivide an edge of $S$ to obtain a new triangulation $S'$ and
a double branched cover $\pi\colon M' \to M$.

\item If necessary, perturb the $1$--form $\omega$ on $M$ so that it is
transverse to the branch set. If this is possible, then
$\pi^*(\omega)$ will be  a nowhere-zero closed $1$--form on $M'$.
\end{enumerate}

\subsection{Finding a fiber}\label{14.4}
Suppose $F^{n-1}$ is a hypersurface in $M^n$. Let $\widehat M(F)$
denote the result of cutting open $M$ along $F$. If  $\widehat M(F)$
is homeomorphic to $F\times [0, 1]$, then $M$ fibers over $\BS^1$ with
fiber $F$. (Proof: M is obtained by gluing together two ends of
$F\times [0, 1]$, ie, it is the mapping torus of a homeomorphism.)

There are some obvious necessary conditions on $F$ for it to be a
fiber. For example, its Betti numbers must be fairly large. Indeed, by
the Wang sequence, we have the  inequality: $b_i(F)+b_{i-1}(F)\ge
b_i(M)$. Furthermore, if $M$ is aspherical, then $F$ must be also
aspherical and the induced homomorphism $\pi_1(F) \to \pi_1(M)$ must
be an injection onto a normal subgroup.

Suppose that $p\colon M \to \BS^1$ is a fibration with fiber $F$ and
with fibration $1$--form $\omega=p^*(d\theta)$. Assume that $M$ is
oriented. Give $F$ the induced orientation and let $[F]$ denote the
image of its fundamental class in $H_{n-1}(M^n;\BR)$. Then $[F]$ is
Poincar{\'e} dual to $[\omega]$.

If $\omega$ is any nowhere-zero $1$--form, then $\Ker \omega$ is an
oriented $(n-1)$--dimensional subbundle of the tangent bundle of
$M^n$. Let  $e_\omega \in H^{n-1}(M^n;\BZ)$ denote the Euler class of
this bundle. When $\omega$ is the $1$--form of a fibration, we clearly
have that $e_\omega([F])=\chi(F)$.

\subsection{Cutting and pasting hypersurfaces}\label{14.5}
In this subsection we shall discuss a procedure for amalgamating
oriented hypersurfaces in an arbitrary closed oriented $n$--manifold
$M^n$.

\subsubsection{Switching sheets}\label{14.5.1}
Suppose $N_1^{n-1}$ and $N_2^{n-1}$ are oriented hypersurfaces in
$M^n$ intersecting transversely. As in~\cite{Th1} we can associate to
this situation a new oriented hypersurface $N_{12}$ in $M$ with the
following two properties:
\begin{enumerate}
\item\label{Class} $[N_{12}]=[N_{1}]+[N_{2}]\qquad$ and

\item\label{Euler} $\chi(N_{12})=\chi(N_{1})+\chi(N_{2})$.
\end{enumerate}

The procedure can be described as follows. A neighborhood of $N_1 \cap
N_2$ in  $N_1 \cup N_2$ is homeomorphic to the product of   $N_1 \cap
N_2$ with the cone over $4$ points. Replace this neighborhood by two
copies of $(N_1 \cap N_2)\times [-1, 1]$ by gluing each side of $N_1
\cap N_2$ in $N_1$ to the appropriate side of $N_1 \cap N_2$ in
$N_2$, as indicated in the picture below.

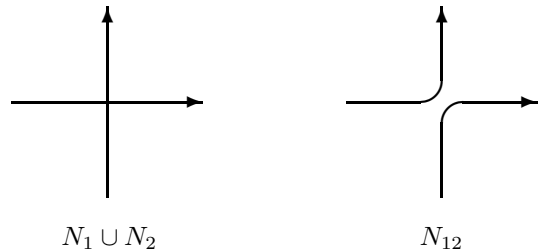
\setlength{\unitlength}{0.25in}
\begin{picture}(9, 6)(-4, -1) \small
\thicklines 
\put(0, 2){\vector(1, 0){4}}  
\put(2, 0){\vector(0, 1){4}} 
\put(1, -1){$N_1\cup N_2$}
\put(7, 0){\begin{picture}(2, 3)
\thicklines 
\put(2.5, 1.5){\oval(1, 1)[tl]}
\put(1.5, 2.5){\oval(1, 1)[br]}

\put(0, 2){\line(1, 0){1.5}} 
\put(2.5, 2){\vector(1, 0){1.5}} 
\put(2, 0){\line(0, 1){1.5}} 
\put(2, 2.5){\vector(0, 1){1.5}}
  
\put(1.5, -1){$N_{12}$}
\end{picture}}
\end{picture}

After a small perturbation we may assume that the result is still
embedded in $M$.

We note that in this procedure each point of $N_1 \cap N_2$ is
replaced by two points in $N_{12}$.

Since $N_{12}$ is just a small perturbation of the $(n-1)$--cycle
$N_1+N_2$, property~\ref{Class} holds.  If $n$ is even, then
property~\ref{Euler} is automatic (since the Euler characteristic of
an odd-dimensional manifold vanishes). If $n$ is odd, then the
codimension--$2$ submanifold $N_1 \cap N_2$ has Euler characteristic
0. Hence, $\chi(N_1\cup N_2)=\chi(N_{1})+\chi(N_{2})=\chi(N_{12})$ so
property~\ref{Euler} holds.

\subsubsection{Iterating this procedure}\label{14.5.2}
Next suppose $N_1$, \dots, $N_m$ are oriented hypersurfaces in general
position in $M$. We can assume that $N_3$ and $N_{12}$ are in general
position. Define $N_{123}$ to be the result of applying the switching
sheets procedure to $N_{12}$ and $N_3$. Continuing in this fashion, we
eventually arrive to an embedded hypersurface $N_{12\dots m}$. As
before, 
\begin{enumerate}
\item\label{Class1} $[N_{12\dots m}]=\sum [N_{i}]\qquad$ and

\item\label{Euler1} $\chi(N_{12\dots m})=\sum \chi(N_{i})$.
\end{enumerate}

We note that in this procedure each point of $N_1 \cup \dots \cup N_m$
which lies in a $j$--fold intersection is blown up into $j$ points in
$N_{12\dots m}$.

In the next subsection we shall give another description of this
procedure which is independent of the ordering of the $N_i$. To give
this description it suffices to consider the local model.

\subsection{The local model}\label{14.6}
For each $i$, $1\le i \le n$, let $P_i$ denote the coordinate
hyperplane in $\BR^n$ defined by $x_i=0$. Orient $P_i$ by requiring
its unit normal vector $e_i$ to be positively oriented.

Let $\lambda=(\lambda_1, \dots , \lambda_n)$ be any function from
$\{1, \dots , n\}$ to $\{-1, 0, +1\}$. Let $z(\lambda)$ denote the number
of zeroes in $(\lambda_1, \dots , \lambda_n)$ and $n(\lambda)$ the
number of $(-1)$'s.

The {\em quadrant} $Q_\lambda$ corresponding to $\lambda$ is the
subset of $\BR^n$ defined by
$$Q_\lambda=\left\{(x_1, \dots , x_n)\in \BR^n\ | \lambda_i x_i \ge
    0 \text{ if $\lambda_i\not=0$; }x_j=0 \text{ if $\lambda_j=0$}
  \right\}.$$
 Clearly, $Q_\lambda$ is isomorphic to the cone on a
simplex of dimension $n-z(\lambda)-1$. It is a manifold with corners
of codimension $z(\lambda)$ in $\BR^n$.

Each hyperplane $P_i$ is divided into $2^{n-1}$ $(n-1)$--dimensional
quadrants. Thus, in total there are $n2^{n-1}$ $(n-1)$--dimensional
quadrants. We shall now reassemble these into $n$ different
sheets. For $l=0, \dots , n-1$, let $\E(l)$ denote the set of functions
$\lambda$ with  $z(\lambda)=1$ and $n(\lambda)=l$.  Define

\subsubsection{}\label{14.6.1} 
$\displaystyle P(l)=\bigcup_{\lambda \in \E(l)} Q_\lambda$.

As we shall see below, $P(l)$ is a piecewise differentiable
submanifold of $\BR^n$ which is homeomorphic to $\BR^{n-1}$. Moreover, 
$P_{12\dots n}$ can be identified with the disjoint union of the
$P(l)$.

The whole arrangement of the $Q_\lambda$ is isomorphic to the cone
over the  triangulation $O^{n-1}$ of $\BS^{n-1}$ as the boundary of the
standard $n$--dimensional octahedron. The vertex set of $O^{n-1}$ is
$\{\pm e_i\}_{1\le i\le n}$ and the simplex corresponding to $\lambda$
is the spherical $(n-z(\lambda)-1)$--simplex spanned by $\{\lambda_j
e_j\}_{\lambda_j\not =0}$.

\begin{Lemma}\label{14.6.2}
Let $O^{n-1}$ be the boundary of the $n$--dimensional octahedron. For
$l=0, \dots , n-1$, let $B(l)$ denote the union of the
$(n-1)$--simplices in $O^{n-1}$ with $n(\lambda) \le l$. Then $B(l)$ is
a topological ball.
\end{Lemma}
\begin{proof}
$B(0)$ is an $n-1$--simplex and $B(l)$ collapses onto $B(l-1)$.
\end{proof}

\subsubsection{}\label{14.6.3}
Clearly, $\partial B(l)$ is the union of the   $(n-2)$--simplices in
$O^{n-1}$ with $n(\lambda)= l$. Hence, $P(l)$ (defined
in~\ref{14.6.1}) is homeomorphic to the cone on $\partial B(l)$. By
the above lemma,  $\partial B(l)$ is a triangulation of
$\BS^{n-2}$. Thus,  $P(l)$ is homeomorphic to $\BR^{n-1}$.

\subsubsection{}\label{14.6.4} The $P(l)$ are not disjointly embedded (since 
$\bigcup P(l)=\bigcup P_i$). To remedy this we alter the embedding of
$P(l)$ in $\BR^n$ by a small isotopy as follows. Choose a decreasing
sequence of real numbers $\mu_0>\mu_1>\dots >\mu_{n-1}$. Let $e$ be
the vector $(1, 1, \dots , 1)$ in  $\BR^n$. Finally let $P'(l)$ be the
subset of $\BR^n$ defined by
$$P'(l)=P(l)+\mu_l e.$$ The $P'(l)$ are now disjointly embedded. This
gives the desired local description of $P_{12\dots n}$: it is the
union of the $P'(l)$.

\subsection{Potential fibers}\label{14.7}
We return to our consideration of the orbifold $\Sigma/W$. Let
$\Gamma$ be a normal, torsion-free subgroup of finite index in $W$ and
set $M^n=\Sigma/\Gamma$.  It follows from the assumptions that
$\Gamma$ is  normal and  torsion-free that, given any wall of
$\Sigma$, its image in $M^n$ is an embedded hypersurface
$N^{n-1}$. We call such an $N$ a {\em standard hypersurface} in $M$. A
standard hypersurface is totally geodesic in the nonpositively curved
cubical structure on $M$.

By passing to a deeper subgroup of finite index if necessary, we may
assume that $M$ is orientable and that each  standard hypersurface is
orientable.

If $W$ does not split as product with an infinite dihedral group, then
a standard hypersurface can never be the fiber of a fibration over
$\BS^1$. However, the cutting and pasting procedure of~\ref{14.5}
applied to various collections of oriented  standard hypersurfaces
$\{N_1, \dots, N_m \}$ yields a good source of candidates for fibers.

\subsubsection{Fundamental domains for the commutator subgroup}\label{14.7.1} 
In this subsection $\Gamma$ is the commutator subgroup of $W$ and
$M^n$ is the commutator cover of $\Sigma/W$ (see~\ref{6.4}). The
quotient group $W/\Gamma$ is $\left( \BZ_2\right)^{{\mathcal S}_0(S)}$ where
${\mathcal S}_0(S)$ denotes the vertex set of $S$.  Order elements of 
${\mathcal S}_0(S)$:
$s_1$, $s_2$, \dots , $s_p$. Next we shall inductively define   an
increasing sequence, $D(0)\subset \dots \subset D(p)$, such that each
$D(i)$ is a convex union of chambers and such that for all $j>i$, 
$\Sigma(s_j)$ is a supporting wall of $D(i)$. (See~\ref{12.1}.)

Put $D(0)=K$. Assuming $D(i)$ has been defined for some $i<p$, define
$D(i+1)$  to be the double of  $D(i)$ along $\Sigma(s_{i+1})$. Set
$D=D(p)$. We claim that $D$ has the following properties:
\begin{enumerate}
\item\label{Convex} $D$ is a convex union of chambers.
\item\label{Fund} $D$ is a fundamental domain for the $\Gamma$--action
on $\Sigma$.
\item\label{Pairs} The codimension-one faces of $D$ are identified in
pairs by 
elements of~$\Gamma$.
\item\label{Standard} The image of the boundary $\partial D$ of $D$ in 
$M$ is a union of
standard 
hypersurfaces $\{N_1, \dots, N_m \}$.
\end{enumerate}

Property~\ref{Convex} is immediate. To see property~\ref{Fund} first
observe that any chamber in $D$ has the form $s_{i_t}\cdots s_{i_1}K$
where $(i_1, \dots , i_t)$ is an increasing sequence of integers
(possibly the empty sequence) in $[1, p]$. Since the group elements
corresponding to such sequence map bijectively onto the quotient group
$W/\Gamma=(\BZ_2)^p$, we see that $D$ is a fundamental domain.

Similarly, any supporting wall of $D$ can be written in the form 
$s_{i_t}\cdots s_{i_1}\Sigma(s_j)$, where $(i_1, \dots , i_t)$ is 
a nonempty increasing
sequence of integers in $[1, p]$ and $j\ge i_t$, and where
$s_{i_t}\cdots s_{i_1}$ does not commute with $s_j$. In particular, 
consider the supporting walls $w\Sigma(s_j)$  and $s_j w\Sigma(s_j)$
where $w=s_{i_t}\cdots s_{i_1}$ and $j>i_t$. The element $s_j w s_j
w^{-1}$ takes the first wall to the second and this element lies in
the commutator subgroup $\Gamma$. Property~\ref{Pairs} follows. The
image of $\partial D$ in $M$ is the same as the image of the union of
the  supporting walls of $D$ in $M$. Hence, property~\ref{Standard}
holds.

\subsubsection{}\label{14.7.2} 
Once again, $\Gamma$ is an arbitrary normal, torsion-free subgroup of
finite index in $W$ such that $M$ and the standard hypersurfaces are
orientable. We further suppose that $\Gamma$ has a fundamental domain
satisfying properties~\ref{Convex} through~\ref{Standard}
in~~\ref{14.7.1}. (By~\ref{14.7.1} such $\Gamma$ exist.) Since $D$ is 
convex, it
is a disk. Let $N_1, \dots , N_m$ be the standard hypersurfaces coming
from the supporting walls of $D$. We shall now describe an attempt
to fiber $M$ over $\BS^1$ which sometimes works.

Let $R$ be a regular neighborhood of $\bigcup N_i$ in $M$. Since $M-R$
can be identified with the complement of a collared  neighborhood of
$\partial D$ in $D$, it is a disk. Thus, $M$ is formed by attaching a
$n$--disk to  $\bigcup N_i$. Since $M$ is orientable, the attaching map
$\partial D\to \bigcup N_i$ is trivial on homology. Hence, for $j<n$, 
$H_j(M)\cong H_j(\bigcup N_i)$. Using Mayer--Vietoris sequences it is
easy to see that we have an injection from $\bigoplus
H_{n-1}(N_i)$~($\cong \BZ^m$) into $ H_{n-1}(\bigcup N_i)$. (In
favorable circumstances this injection is an isomorphism.)

Now choose an orientation for each $N_i$ and a positively oriented
normal vector field $v_i$ on $N_i$. Pulling this back to $D$ we obtain
a normal vector on each of its codimension-one faces. If two faces are
identified by an element of $\Gamma$, then the vectors point in
opposite directions. Call a codimension-one face {\em positive}
(respectively, {\em negative}) if the normal vector is outward pointing
(respectively, inward pointing). Let $D_+$ (respectively  $D_-$) denote 
the union of
the positive (respectively, negative)  codimension-one faces. Thus, each
choice of the orientations of the $N_i$ leads to a partition of
$\partial D$ into positive and negative regions.

\begin{Proposition}\label{14.7.3}
With hypothesis as above, suppose it is possible to choose
orientations for $N_i$ such that the positive region $D_+$ is a disk
(of codimension $0$) in $\partial D$. Then  $N_{12\dots m}$ is the
fiber of a fibration of $M$ over $\BS^1$.
\end{Proposition}

\begin{proof} Let $F= N_{12\dots m}$ and let $\widetilde{F}$ denote the
inverse image of $F$ in $\Sigma$. Then $D_+$ (or $D_-$) can be
regarded as a fundamental domain for the $\Gamma$--action on
$\widetilde{F}$. Let $\widehat{M}(F)$ denote the result of cutting $M$
open along $F$.  Take a component of the complement of
$\widetilde{F}$ in $\Sigma$ and let $\widetilde{M}$ be its
closure. Then $\widetilde{M}$ is a covering space of $\widehat{M}(F)$
with group of covering transformations $\Gamma'$. Furthermore, $D$ is
a fundamental domain for the $\Gamma'$--action on $\widetilde{M}$. The
only points of $D$ that are identified under $\Gamma'$ lie on the
common boundary of $D_+$ and $D_-$. It follows that $\widehat{M}(F)$
is a quotient space of $D$ by an equivalence relation $\sim$. Since by
hypothesis, D is homeomorphic to $D_+ \times [0, 1]$, it follows that
the quotient space is homeomorphic to $(D_+/\sim) \times [0, 1]$, where
$(D_+/\sim)=F$. The proposition follows.
\end{proof}

\end{document}